\documentclass[12pt]{amsart}
\usepackage{amsfonts,amssymb,amsmath,amscd,amsthm, amstext,color}
\usepackage[colorlinks, citecolor=blue,pagebackref,hypertexnames=false]{hyperref}
\usepackage{enumerate}
\usepackage{epstopdf}
\usepackage{fancyhdr}
\usepackage{geometry}
\usepackage{graphicx}
\usepackage{mathrsfs}
\usepackage{txfonts}
\usepackage{tikz}
\usetikzlibrary{arrows.meta, positioning, decorations.pathreplacing}
\usepackage{url}

\usepackage{comment}
\allowdisplaybreaks

\newtheorem{theorem}{Theorem}[section]
\newtheorem{thmA}{Theorem}

\newtheorem{proposition}[theorem]{Proposition}
\newtheorem{lemma}[theorem]{Lemma}
\newtheorem{corollary}[theorem]{Corollary}
\newtheorem{remark}[theorem]{Remark}

\def\supp {{\rm supp}\,}

\def \f {\frac}

\title[ Chernoff--Hoeffding inequalities]{ Chernoff--Hoeffding inequalities,  almost orthogonality and bounded multi-parameter square functions}

\thanks{{\it {\rm 2010} Mathematics Subject Classification:} 42B20, 42B25, 60G48.  }
\thanks{{\it Key words and phrase:} Chernoff--Hoeffding type inequality; square functions;
bi-parameter setting; Shilov boundary}
 
 \author{Ji Li}
\address{Ji Li, School of Mathematical and Physical Sciences, Macquarie University, NSW 2109, Australia}
\email{ji.li@mq.edu.au}

\author{Jill Pipher}
\address{Department of Mathematics, Brown University, Providence, RI 02906 USA}
\email{jpipher@math.brown.edu}

\author{Liangchuan Wu}
\address{Liangchuan Wu,  School of Mathematical Science, Anhui University, Hefei, 230601,  P.R.~China and School of Mathematical and Physical Sciences, Macquarie University, NSW 2109, Australia}
\email{wuliangchuan@ahu.edu.cn}
 \date{}
 
\begin{document}

\begin{abstract}

 We establish multi-parameter Hoeffding inequalities for the summation of atomic functions possessing suitable almost-orthogonality
 via introducing a new square function and by constructing an iterative algorithm on scales of this new square function. 
The main result yields new multi-parameter Hoeffding inequalities, and generalize those for dyadic martingales. In addition, it applies to broader settings such as the Shilov boundary of tensor product domains in $\mathbb C^{2n}$, $n\geq2$, studied by Nagel and Stein (2004).

As an application, we obtain the sharp order of local integrability for our generalized multi-parameter square functions, extending the classical result of Chang--Wilson--Wolff (1985) in the one-parameter setting and the second author’s work on the bi-disc (1986).  
\end{abstract}
 
\maketitle

\tableofcontents

\section{Introduction and Statement of Main Results}\label{sec:1}
\setcounter{equation}{0}
 
\subsection{Background and aim}

A fundamental question in probability is the likelihood of a random variable being close to its expectation. Some of the most widely used tools for this question are various {\it concentration  inequalities}. Among these, the well-known Chernoff--Hoeffding bound (also known as the  Chernoff--Hoeffding inequality or Hoeffding's inequality)  is probably one of the most widely employed. It states that the sum $X=\sum_{j=1}^N X_j$ of random variables $X_j$ independently distributed in $[0,1]$,   is highly concentrated around its expectation $  E[X]$ (\cite{Ch52, Hoe63}): for all $t>0$,
\begin{equation}\label{eqn:C-H}
    \mathbf{P}\,\big(X-  E[X]>t\big), \  \   \mathbf{P}\,\big(X-  E[X]<-t\big)\leq e^{-2t^2/N}.
\end{equation}
 One may view \eqref{eqn:C-H} as a concentrated version of the central limit theorem. Notably, this result has an analog in terms of dyadic martingales (see \cite{CWW}). Specifically,
 let $\mathcal F_n$ be the $\sigma$-field generated by the  dyadic intervals of length $2^{-n}$ in $[0,1]$, and suppose
 that the sequence $\{f_n\}_{n\geq 1}$ is a dyadic martingale: $f_n$ is measurable with 
 respect to $\mathcal F_n$ and $ E(f_{n+1}| \mathcal F_n)=f_n$ for all $n$. Assume $f_0=d_0=0$ and set $d_n(x)=f_n(x)-f_{n-1}(x)$, $n\geq1$. Then for any $N\in \mathbb N_+$ and for any $\alpha>0$, we have the following inequality due to H. Rubin,
\begin{equation}\label{eqn:ESI}
    \int_{[0,1]} \exp\Big(  \alpha F_N(x)-\frac{\alpha^2}{2} \mathcal S^2(F_N)(x)\Big)\, {\rm d}x\leq 1,
\end{equation}
where $F_N(x):=\sum_{n=1}^N d_n(x)$ and the square function $\mathcal S$ is defined as 
$
    \mathcal S(F_N)(x):=\big(\sum_{n=1}^N d_n^2(x)\big)^{1\over2}.
$

The higher dimensional version of \eqref{eqn:ESI}  is due to Chang--Wilson--Wolff \cite[Theorem 3.1]{CWW}: 
\begin{equation}\label{eqn:ESI higher}
    \int_{[0,1]^d} \exp\Big(  \alpha F_N(x)-\frac{\alpha^2}{2} \mathcal S^2_{\mathbb R^{\mathrm d}}(F_N)(x)\Big)\, {\rm d}x\lesssim 1,
\end{equation}

where the square function $\mathcal S_{\mathbb R^{\mathrm d}}$ is defined as follows:
\begin{align}\label{eqn: square}
    \mathcal S^2_{\mathbb R^{\mathrm d}}(F_N)(x) = \sum_{n=1}^N\left \| \mathsf 1_{Q_{n-1}(x)}(\cdot)\, d_n(\cdot) \right\|_{L^\infty}^2,\quad x\in\mathbb R^{\mathrm d}
\end{align}  
with $Q_{n-1}(x)$ the unique dyadic cube in $[0,1]^d$ of side length $2^{-(n-1)}$ containing $x$.

Versions of inequality \eqref{eqn:ESI higher} hold on
any fixed dyadic cube $Q$ in $\mathbb R^{\mathrm d}$, ${\mathrm d}>1$, with $F_N$ and $\mathcal S$ appropriately localized to $Q$.

By a standard argument via iteration, the proof of \eqref{eqn:ESI} (as well as \eqref{eqn:ESI higher}) relies indispensably on the cancellation $\int_{I_n} d_n =0$ for each $n$ and dyadic interval of sidelength $2^{-n}$,  from the perspective of  Hoeffding's lemma: $ E[e^{\alpha (X-  E[X])}]\leq e^{\frac{\alpha^2 (b-a)^2}{8}} $ for any random variable $X$ with $a\leq X\leq b$ and any $\alpha>0$. Additionally, the mutual independence of random variables $X_j$ in \eqref{eqn:C-H} corresponds to the orthogonality $\int d_n(x) d_m(x)\,dx =0$ for $n\neq  m$ of martingale differences in \eqref{eqn:ESI}, this orthogonality is a natural consequence of the cancellation of $d_n$'s, and implies 
$\int  \big\| \sum_n d_n(x)\big \|^2  {\rm d}x=\int \sum |d_n(x)|^2 {\rm d}x$.

These fundamental results \eqref{eqn:C-H}---\eqref{eqn:ESI higher} have various applications in different fields, e.g., sharp estimates for square functions (Chang--Wilson--Wolff \cite{CWW}),  random data Cauchy theory for supercritical wave equations (Burq--Tzvetkov \cite{BT1, BT2}), {and the sparse dimensionality reduction in a variety of disciplines (Bourgain--Dirksen--Nelson \cite{BDN15}).} We also refer to \cite{BOP1, BOP3, B,  GHS, IT, Pe, T95, T96, Wi89, Wi08} and references therein for extensive applications in the past few decades.

In the bi-parameter setting $[0,1]\times [0,1]$, the second author \cite{Pi86} considered the double dyadic martingales $\{f_{n,m}\}$ in $[0,1]\times [0,1]$, where $f_{n,m}$ is a dyadic martingale in the first component for each fixed $m$, and also a dyadic martingale in the second component for each fixed $n$. In particular, we also assume $f_{n,0}=f_{0,m}=0$ for each $n$ and $m$. Let $d_{n,m}=f_{n,m}-f_{n-1,m}-f_{n,m-1}+f_{n-1,m-1}$, $n,m\geq1$. Then
 from \eqref{eqn:ESI}, one obtains that for any fixed $x_2\in[0,1]$ and any $\alpha>0$,
 \begin{align}\label{S1 bi-para}
\int_0^1\exp\bigg( \alpha F_{N,N}(x_1,x_2) - {\alpha^2\over2} \mathcal S^2_1(F_{N,N})(x_1,x_2) \bigg) \, {\rm d}x_1\leq1,
\end{align}  
and for any fixed $x_1\in[0,1]$,  there exist $c,C>0$ such that for any $\alpha>0$,
\begin{equation}\label{S bi-para}
   	\int_{0}^1 \exp\left(  \alpha \mathcal S_1(F_{N,N})(x_1,x_2) -c \alpha^2 \mathcal S^2(F_{N,N})(x_1,x_2)\right){\rm d} x_2\leq C,
\end{equation}
where $F_{N,N} = \sum_{n=1}^N \sum_{m=1}^N d_{n,m}$,
$$\mathcal S^2_1(F_{N,N})(x_1,x_2) = \sum_{n=1}^N\Big(\sum_{m=1}^N d_{n,m}\Big)^2(x_1,x_2)\qquad   \text{and}  \qquad 
    \mathcal S^2(F_{N,N})(x_1,x_2) = \sum_{n=1}^N\sum_{m=1}^N d_{n,m}^2(x_1,x_2). 
$$

Based on \eqref{S1 bi-para} and \eqref{S bi-para}, a bi-parameter version of Rubin's Lemma was established in \cite[Section 2]{Pi86} as follows: there exist positive constants $c$ and $C$ such that for each $N\in \mathbb N_+$, such that  
\begin{equation}\label{eqn:double-ESI intro}
   \int_{[0,1]^2} \exp\Bigg(\frac{c\sum_{n=1}^{N} \sum_{m=1}^{N} d_{n,m} }{\big\| \big(\sum_{n=1}^{N} \sum_{m=1}^{N} d_{n,m}^2\big)^{1/2} \big\|_{L^\infty}}\Bigg)\, {\rm d}x_1 {\rm d}x_2 \leq C.
\end{equation}

A natural {\bf question} arises as inspired by the one-parameter result \eqref{eqn:ESI higher}  (Chang--Wilson--Wolff \cite{CWW}) in  $\mathbb R^{\mathrm d}$ and by the sharp lower bound of multi-parameter square functions:

\smallskip

{\bf  Can we establish \eqref{eqn:double-ESI intro} in higher dimensional spaces $\mathbb R^{\mathrm d_1}\times \mathbb R^{\mathrm d_2}$ and more general settings\footnote{There are many fundamental tensor product spaces arising from several complex variables beyond $\mathbb R^{\mathrm d_1}\times \mathbb R^{\mathrm d_2}$, one of which is the well-known Shilov boundary $\mathscr M_1\times\cdots\times \mathscr  M_n$ of tensor product domains $\Omega_1\times \cdots\times \Omega_n$ in $\mathbb C^{2n}$, $n\geq2$, studied by Nagel and Stein \cite{NS04}. See also the decoupled boundaries by Nagel and Stein \cite{NS06}.} ?}
\medskip

One may think that it is straightforward to first extend both \eqref{S1 bi-para} and \eqref{S bi-para} to higher dimensional space, i.e., similar to \eqref{eqn:ESI higher} as in \cite{CWW}.
However, the above argument  (\cite{Pi86}) can not be applied to 
$\mathbb R^{\mathrm d_1}\times \mathbb R^{\mathrm d_2}$ directly. To see this, we note that for any fixed $x_2\in \mathbb R^{\mathrm d_2}$,  the extension of \eqref{S1 bi-para} in $\mathbb R^{\mathrm d_1}$ (as a consequence of \eqref{eqn:ESI higher} by \cite{CWW})
is of the form 
\begin{equation}\label{eqn: 1.8}
    \int_{Q_1} \exp\Big(  \alpha F_{N,N}(x_1,x_2)-\frac{\alpha^2}{2} \mathcal S^2_{\mathbb R^{\mathrm d_1},1}(F_{N,N})(x_1,x_2)\Big) \,{\rm d}x_1\leq |Q_1|,\qquad \forall\alpha>0,
\end{equation}
where $Q_1$ is any fixed dyadic cube in $\mathbb R^{\mathrm d_1}$, and the square function $\mathcal S_{\mathbb R^{\mathrm d_1},1}$ is defined as follows:
\begin{align}\label{eqn: 1.9}
    \mathcal S^2_{\mathbb R^{\mathrm d_1},1}(F_{N,N})(x_1,x_2) = \sum_{n=1}^N \bigg\| \mathsf 1_{Q_{n-1}(x_1)}(\cdot) \sum_{m=1}^N d_{n,m}(\cdot, x_2)\bigg\|_{L^\infty(\mathbb R^{\mathrm d_1})}^2
\end{align}  
with $Q_{n-1}(x_1)$ is the unique dyadic cube in $Q_1$ of side length $2^{-(n-1)}\ell(Q_1)$ containing $x_1$. It is obvious that this $L^\infty(\mathbb R^{\mathrm d_1})$-norm 
over $\sum_{m=1}^N d_{n,m}$ destroys the martingale structure (and hence the cancellation, orthogonality or almost orthogonality) in the second component in $\mathbb R^{\mathrm d_2}$, thus, the next extension of \eqref{S bi-para} to $\mathbb R^{\mathrm d_2}$ does not work either.

\smallskip

{\bf The aim of this paper is to address the above question} by establishing a multi-parameter analogue of \eqref{eqn:double-ESI intro} in general product spaces of homogeneous type, which covers not only the bi-disc and $\mathbb R^{\mathrm d_1}\times \mathbb R^{\mathrm d_2}$, but also the Shilov boundary $\mathscr M_1\times\cdots\times  \mathscr M_n$, $n\geq2$, defined by Nagel and Stein \cite{NS04}. 

We provide a different approach via introducing a new square function and constructing an iterative algorithm on scales of the square functions first in the one-parameter setting and then in the bi-parameter setting. This overcomes the above obstacles in \eqref{eqn: 1.8} arising from the structure of \eqref{eqn: 1.9}. Moreover, our approach bypasses the dependence on precise cancellation, and hence allows for a wider class of square functions beyond the  dyadic martingales (see examples in Section \ref{sec:5}). {Additionally, note that the concentration of measure phenomenon often occurs without independence, which motivates the exploration of  Hoeffding's inequalities in certain dependent settings (see \cite[Chapter 3]{DP} and \cite[Chapter 5]{V} for examples). Our framework, which does not rely on precise cancellation, is also partly inspired by this perspective. Indeed, our result not only encompasses the classical Hoeffding's inequalities related to martingale differences,  but also sharpens them with a new,  suitable square function under weaker orthogonality conditions, which can be regarded as a form of limited dependence in probabilistic terms, and perhaps has the potential to become an useful tool for a variety of related problems. }

Some direct applications include the sharp $L^p$ lower bound of Lusin area function on  $\mathbb R^{\mathrm d_1}\times \mathbb R^{\mathrm d_2}$ and on $\mathscr M_1\times\cdots\times \mathscr  M_n$ studied by Nagel and Stein \cite{NS04}.

\subsection{Statement of main results}   
Let $(X,\rho,\mu)$ be a space of homogeneous type in the sense of Coifman and Weiss and let $\mathcal D=\{\mathcal D_k\}_{k\in \mathbb{Z}}$ be a system of dyadic cubes on $X$
(full definition will be provided in Section \ref{sec:2}).
Our first step is the following one-parameter Hoeffding type inequality in $(X,\rho,\mu)$, 
where we only require almost-orthogonality, instead of specializing to dyadic martingales.

 \begin{thmA}  \label{thm:newCH-dyadic}
Let  $(X,\rho,\mu)$ be a quasi-metric measure space satisfying the volume doubling condition.
For any given $Q\in \mathcal D$, let $\big \{\mathfrak{a}_P:\, P\in \mathcal D,\, P\subseteq Q\big  \}$ be ``atomic functions'' satisfying
\begin{itemize}
	\item[\textbf{(P1)} ] $\mathfrak{a}_P$ is supported on $P$;
	\smallskip
	\item[\textbf{(P2)} ]  $\mathfrak{a}_P$ is constant on each $P'\in Ch(P)$, here $Ch(P)$ represents the set of the children of $P$.
	\smallskip
	\item[\textbf{(P3)} ]  \textbf{{\rm(}almost orthogonality{\rm)}} there exists  $C_1>0$ such that for any dyadic cube $R\subseteq Q$ and any sub-collection $\mathcal G_R\subseteq \{P\in \mathcal D:\, P\subseteq R\}$,
\begin{equation}\label{eqn:P3}
	\bigg(\frac{1}{\mu(R)}  \int_R  \Big|   \sum_{P\in \mathcal G_R} \mathfrak{a}_P(x)\Big|^2 {\rm d}\mu(x) \bigg)^{1/2} \leq C_1\,  \bigg\|\,  \bigg(\sum_{ P\in \mathcal G_R}   {\mathfrak{a}_P^2(x) }\bigg)^{1/2}\,\bigg\|_{L^\infty(R)} .
\end{equation}
\end{itemize}

Then there exist positive constants $c_1$, $c_2$, $c_3$ and $c_4$ independent of $Q$ such that for any $\alpha>0$ and $N\in \mathbb N_+$, 
\begin{equation}\label{eqn:newCH-dyadic}
       \int_Q \exp\left(\alpha \left|F_N(x)\right|  
       \right)  {\rm d}\mu(x)\leq c_1 \int_Q \exp\left( c_2\alpha^2 S^2(F_N)(x)
       \right)  {\rm d}\mu(x),
\end{equation}
and that  for any $N\in \mathbb N_+$, 
\begin{equation}\label{eqn:newCWW-one-parameter}
       \int_Q \exp\left(    c_3    \frac { |F_N(x) |^2 }   {\|S(F_N)\|_ {L^\infty(Q) }^2  }
       \right)  {\rm d}\mu(x)\leq c_4  \mu(Q),	
\end{equation}
where 
$$
    F_{N}(x)=\sum_{\substack{P\in  Ch^{(n)}(Q)\\ 0\leq n\leq N}}   \mathfrak{a}_{P} ( x) \quad    {\rm and}\quad    
     S(F_N)(x)=\bigg(\sum_{ \substack{P\in  Ch^{(n)}(Q)\\ 0\leq n\leq N}}  {\mathfrak{a}_P^2 (x) }\bigg)^{1/2}.
$$
\end{thmA}

\smallskip

Our notation ($F_N, S(F_N)$, etc.) suppresses the dependence on $Q$ of these quantities, as we shall always assume that $F$ and its square function
are suitably localized to the cube $Q$ over which the integration takes place.

Based on Theorem \ref{thm:newCH-dyadic}, we will further establish the bi-parameter Hoeffding type inequalities (and one may refer to Proposition \ref{prop:multi-newCH-dyadic}  for the multi-parameter version). 

 \begin{thmA}   \label{thm:double-newCH-dyadic}
 For   $j=1,2,$ let  $(X_j,\rho_j,\mu_j)$ be a quasi-metric measure space satisfying the volume doubling condition and let $\mathcal D^{[j]}$ be a system of dyadic cubes on $X_j$. For any given $Q_1\times Q_2\in \mathcal D^{[1]}\times \mathcal D^{[2]}$, let $\left \{\mathfrak a_{P_1\times P_2}:\, P_1\times P_2\in \mathcal D^{[1]}\times \mathcal D^{[2]}, \,  P_1\subseteq Q_1,\, P_2\subseteq Q_2\right\}$ be ``atomic functions'' satisfying
 \begin{itemize}
	\item[\textsf{(P1)} ] $\mathfrak{a}_{P_1\times P_2}$ is supported on $P_1\times P_2$.
	\smallskip
	\item[\textsf{(P2)} ]  $\mathfrak{a}_{P_1\times P_2}$ is constant on each $P'_1\times P'_2\in Ch(P_1\times P_2):=Ch(P_1)\times Ch(P_2)$.	
	\smallskip
	\item[\textsf{(P3)} ]  \textbf{{\rm(}almost orthogonality{\rm)}} there exists  $C_2>0$ such that for any dyadic rectangle $R_1\times   R_2\subseteq Q_1\times Q_2$ and any subcollection $\mathcal G_{R_1\times R_2}\subseteq \left\{P_1\times P_2\in \mathcal D^{[1]}\times \mathcal D^{[2]} :\, P_1\times P_2 \subseteq R_1\times R_2\right\}$,
\begin{align}\label{eqn:double-P3-revise-I}
	& \bigg(\int_{R_1} \bigg |   \sum_{P_1}  \sum_{P_2:\, P_1\times P_2\in \mathcal G_{R_1\times R_2}} \mathfrak a_{P_1\times P_2} (x_1,x_2)\bigg |^2  {\rm d}\mu_1(x_1)\bigg)^{1/2}\\
	& \leq\, C_2   \bigg(\int_{R_1}   \sum_{P_1}  \bigg | \sum_{P_2:\, P_1\times P_2\in \mathcal G_{R_1\times R_2}} \mathfrak a_{P_1\times P_2} (x_1,x_2)\bigg |^2  {\rm d}\mu_1(x_1)\bigg)^{1/2},\qquad\qquad \forall{\rm\  fixed\ } x_2\in R_2;\ \  {\rm and}\nonumber
\end{align}
\begin{align}\label{eqn:double-P3-revise-II}
	&\bigg(\int_{R_2} \bigg |    \sum_{P_2:\, P_1\times P_2\in \mathcal G_{R_1\times R_2}} \mathfrak a_{P_1\times P_2} (x_1,x_2) \bigg |^2  {\rm d}\mu_2(x_2)\bigg)^{1/2}\\
	&\leq\,  C_2   \bigg(\int_{R_2}     \sum_{P_2:\, P_1\times P_2\in \mathcal G_{R_1\times R_2}} \left|\mathfrak a_{P_1\times P_2} (x_1,x_2)\right|^2  {\rm d}\mu_2(x_2)\bigg)^{1/2}, \qquad \forall{\rm\  fixed\ }P_1 {\rm\ and\ } x_1\in P_1.\nonumber
\end{align}
\end{itemize}

Then there exist positive constants $c'$ and $C'$ independent of $Q_1$ and $Q_2$ such that for any $N\in \mathbb N_+$,
\begin{equation}\label{eqn:double-newCH-dyadic}
 \int_{Q_1\times Q_2} \exp\Bigg(c' \frac{   | F_{N,N}( x_1, x_2)|   }{\|S(F_{N,N})\|_{L^\infty(Q_1\times Q_2)}}\Bigg)  \,{\rm d}\mu_1(  x_1) \, {\rm d}\mu_2( x_2) \leq C' \mu_1(Q_1)\,\mu_2(Q_2),
\end{equation}
where 
$$
     F_{N,N}(  x_1, x_2)=\sum_{j=1}^2\sum_{\substack{P_j\in  Ch^{(n_j)}(Q_j)\\ 0\leq n_j\leq N}}  \mathfrak{a}_{P_1\times P_2} ( x_1,  x_2)\quad   {\rm and}\quad 
    S(F_{N,N})( x_1,  x_2) =\Bigg(\sum_{j=1}^2\sum_{\substack{P_j\in  Ch^{(n_j)}(Q_j)\\ 0\leq n_j\leq N}}    { \mathfrak{a}^2_{P_1\times P_2}  (x_1,x_2)     }\Bigg)^{1/2}. 
$$
 \end{thmA}

\begin{remark}
{ One key new idea in establishing \eqref{eqn:newCH-dyadic} is to apply the almost
orthogonality \textbf{(P3)} to establish a  John--Nirenberg type estimate associated to our new square function $S(F_N)$ as a bootstrap, and further exploit suitable classifications/rearrangements of dyadic cubes and stopping-time argument to construct an iterative algorithm on scales of $S(F_N)$. One new approach focuses on the definition of our square functions without using the classical $L^\infty$ norm of atom functions as in \eqref{eqn: square}. 

}
\end{remark}

\begin{remark}
Note that for any $f\in L^p(X,{\rm d}\mu)$ with $1<p<\infty$,
 every martingale difference  $\Delta_P(f)$ of $f$ associated to $P\in \mathcal D$
 satisfies \textbf{(P1)}, \textbf{(P2)} and the cancellation property $\int_P \Delta_P(f) (x)\, {\rm d}\mu(x)=0$, where 
 $$
     \Delta_P(f)(x):=\sum_{R\in Ch(P)} f_R \mathsf 1_R(x) -f_P \mathsf 1_P(x)
 $$
and $f_P$ is the average of $f$ over $P$. These conditions ensure that \textbf{(P3)} holds with the constant  $C_1=1$.

With the new square function $S(F_N)$ in Theorem \ref{thm:newCH-dyadic},  our results \eqref{eqn:newCH-dyadic}  and \eqref{eqn:newCWW-one-parameter} are new even  when $\mathfrak a_P$ is the dyadic martingale  $\Delta_P(f)$ associated to $P\in \mathcal D$ for a given function $f\in  L^2(\mathbb R^{\mathrm d})$.  The classical iterative argument on generations of cubes {\rm(}see \cite[Theorem 3.1]{CWW}  
for instance{\rm)} only gives the following argument  
\begin{equation}\label{eqn:classical-dyadic-Hoeffding}
	\int_Q \exp\Bigg (\alpha  \sum_{\substack{P\in  Ch^{(n)}(Q)\\ 0\leq n\leq N}}   \Delta_{P}(f) ( x)  -\frac{\alpha^2}{2} \sum_{\substack{P\in  Ch^{(n)}(Q)\\ 0\leq n\leq N}}   \|\Delta_{P}(f)\|_{L^\infty(P)}^2 \mathsf 1_P(x)  \Bigg)  \, {\rm d}x\leq |Q|,
\end{equation}
which is associated to the usual dyadic square function $\mathcal S_{\mathbb R^{\mathrm d}}(F_N)(x) $ as in \eqref{eqn: square}.  Details are in the appendix.  As a consequence of \eqref{eqn:classical-dyadic-Hoeffding}, we have the classical Chang--Wilson--Wolff result for dyadic martingales  that
\begin{align}\label{CWW classical}
     \int_Q \exp\left(    c_3'    \frac { |F_N(x) |^2 }   {\|\mathcal S_{\mathbb R^{\mathrm d}}(F_N)\|_ {L^\infty(Q) }^2  }
       \right)  {\rm d}\mu(x)\leq c_4'  \mu(Q).
\end{align}
We note that this result \eqref{eqn:newCWW-one-parameter} improves the classical Chang--Wilson--Wolff result \eqref{CWW classical}, since
it involves a smaller square function, i.e.
$
   S(F_N)(x)\leq \mathcal S_{\mathbb R^{\mathrm d}}(F_N)(x).
$


We also highlight that the definition of $\mathcal S_{\mathbb R^{\mathrm d}}(F_N)(x)$ is not suitable for  establishing the bi-parameter version of \eqref{eqn:classical-dyadic-Hoeffding} in $\mathbb R^{\mathrm d_1}\times \mathbb R^{\mathrm d_2}$ with double dyadic martingales.
 
Similarly, the bi-parameter square function $S(F_{N,N})$ in Theorem \ref{thm:double-newCH-dyadic} has the right form to permit an induction to higher parameters - see Section \ref{sec:4} for explanations.
\end{remark}

\begin{remark}
We observe that to prove the one-parameter Chang--Wilson--Wolff type result \eqref{eqn:newCWW-one-parameter}, it suffices to apply
a rough version of  Hoeffding's inequality:
$$
      \int_Q \exp\left(\alpha \left|F_N(x)\right|  
       \right)  {\rm d}\mu(x)\leq c_1   \exp\left( c_2\alpha^2 \|S(F_N)\|_{L^\infty(Q)}^2\right)
      \mu(Q).
$$
However, for the bi-parameter and multi-parameter analogues of \eqref{eqn:newCWW-one-parameter}, i.e.,  Theorem \ref{thm:double-newCH-dyadic} and Proposition \ref{prop:multi-newCH-dyadic}, the corresponding Hoeffding's inequalities in the form of \eqref{eqn:newCH-dyadic} will become crucial. See Section \ref{sec:4} for details.

\end{remark}

\begin{remark}
 Because we only assume almost orthogonality, our atomic functions $ \mathfrak{a}_{P} $ {\rm(}resp.  $\mathfrak{a}_{P_1\times P_2}${\rm)} do not give rise to standard {\rm(}resp. double{\rm)} dyadic martingales. Indeed, $ \mathfrak{a}_{P} $ {\rm(}resp.  $\mathfrak{a}_{P_1\times P_2}${\rm)} are not assumed to have cancellation conditions, as in the classical setting.
 This opens up potential applications to singular integrals in multi-parameter harmonic analysis modeled on domains in several complex variables {\rm(\cite{NS06})},  and to other related questions in stochastic PDEs {\rm (\cite{BOP1,BT1})}.   A final note on the proofs of  statements like \eqref{eqn:newCWW-one-parameter} (resp. \eqref{eqn:double-newCH-dyadic}):  it suffices to consider a suitable decomposition of the function $F_N=\sum F^{(j)}_N$ {\rm(}resp. $F_{N,N}=\sum F^{(j)}_{N,N}${\rm)}  and then prove \eqref{eqn:newCWW-one-parameter} (resp. \eqref{eqn:double-newCH-dyadic}) for each of $F^{(j)}_N$ {\rm(}resp. $F^{(j)}_{N,N}${\rm)} separately, as long as the number of summands is under control.
\end{remark}

\subsection{Applications}\label{sec:1.3}  
We provide some direct applications of our main results above.
 
 \smallskip
 
\noindent{\bf Application I:} Sharp lower bound of the new square functions $S(F_N)$ and $S(F_{N,N})$ beyond the Haar setting introduced in Theorems \ref{thm:newCH-dyadic} and \ref{thm:double-newCH-dyadic}, respectively.

By adapting the idea in \cite{FP} by R. Fefferman and the second author to our setting, we establish the following sharp norm bound for $S(F_N)$. 

\begin{proposition}\label{coro:Lp-CH}
Let all assumptions and notation be the same as in Theorem \ref{thm:newCH-dyadic}.
For any given $Q\in \mathcal D$, let $\big \{\mathfrak{a}_P:\, P\in \mathcal D,\, P\subset Q\big \}$ satisfy \textbf{(P1)}, \textbf{(P2)} and \textbf{(P3)} in Theorem~\ref{thm:newCH-dyadic}. 
We have that for $p\geq 2$,
\begin{equation}\label{eqn:Lp-CH}
   \left \|F_{N}\right\|_{L^p(Q)}\leq C \sqrt{p}\left\|S(F_N)\right\|_{L^p(Q)},
\end{equation}
where the positive constant $C$ is independent of $Q$ and $N$.
\end{proposition}
In particular, this result is a crucial ingredient for the bi-parameter version below.

\begin{proposition}\label{prop:double-Lp-CH-intro}
Let all assumptions and notation be the same as in Theorem \ref{thm:double-newCH-dyadic}. There exists a constant $C>0$ such that for  $p\geq 2$,
\begin{equation}\label{eqn:double-Lp-CH-intro}
	   \left\|F_{N,N}\right\|_{L^p(Q_1\times Q_2,{\rm d}\mu_1{\rm d}\mu_2)}\leq Cp\left\|S(F_{N,N})\right\|_{L^p(Q_1\times Q_2,{\rm d}\mu_1 {\rm d}\mu_2)}.
\end{equation}	
\end{proposition}

\begin{remark}
Note that in both Theorems \ref{thm:newCH-dyadic} and \ref{thm:double-newCH-dyadic} and Propositions \ref{coro:Lp-CH} and \ref{prop:double-Lp-CH-intro} above, the case $N=\infty$ can be verified similarly by using a limiting argument via Fatou's lemma. Hence Propositions \ref{coro:Lp-CH} and \ref{prop:double-Lp-CH-intro} are also available for general $f\in L^p$, $p\geq2$, with $S(f)$ via martingales.

{
Additionally, we highlight that the estimates \eqref{eqn:newCH-dyadic} and \eqref{eqn:newCWW-one-parameter} in Theorem  \ref{thm:newCH-dyadic},  and \eqref{eqn:Lp-CH} in Proposition \ref{coro:Lp-CH} are equivalent; see Section \ref{sec:3} for details.
}
\end{remark}

\begin{remark}
We note that Proposition \ref{coro:Lp-CH} recovers the classical result for the one-parameter square function $\mathcal S_{\mathbb R^{\mathrm d}}(F_N)$ via dyadic martingales by using $S(F_N)(x)\leq \mathcal S_{\mathbb R^{\mathrm d}}(F_N)(x)$ again, while Proposition \ref{prop:double-Lp-CH-intro} provides the missing theory for bi-parameter square function via double dyadic martingales in general $\mathbb R^{\mathrm d_1}\times \mathbb R^{\mathrm d_2}$ {\rm(}note that in $\mathbb R\times \mathbb R$, it was built in \cite{FP} and the bi-disc version was in \cite{Pi86}{\rm)}. Moreover, the multi-parameter version of Proposition \ref{prop:double-Lp-CH-intro} also holds; see {\color{blue}\eqref{eqn:multi-Lp-CH-intro}} in Proposition \ref{prop:multi-newCH-dyadic}.
\end{remark}

%
%
\noindent{\bf Application II:} Sharp lower bound of area functions on the Shilov boundaries.

 Besides $\mathbb R^{\mathrm d_1}\times \mathbb R^{\mathrm d_2}$, one fundamental tensor product space is the Shilov boundary $\widetilde{\mathscr  M} =\mathscr M_1\times  \cdots\times \mathscr M_n$ studied by Nagel and Stein \cite{NS04} which links to the $\bar\partial$-Neumann problem on decoupled  boundaries.
Here
 each $\mathscr M_j$ is an unbounded polynomial domain of finite type $m_j$ defined  as follows. 
\begin{align}\label{manifold M}
  \mathscr M_j:=\big\{ (z_j,w_j)\in\mathbb C^2:\  {\rm Im}(w_j)= \mathcal P_j(z_j) \big\},
\end{align}
where $ \mathcal P_j(z)$ is
a real, subharmonic, non-harmonic polynomial of degree $m_j$. 
Associated with the domain $\mathscr M_j$ is the natural Carnot--Carath\'eodory metric $d_{cc}({\bf x}_j,{\bf y}_j)$ on $\mathscr M_j$ ($\forall\ {\bf x}_j,{\bf y}_j\in \mathscr M_j$) and the measure $\mu_j$ of the nonisotropic ball $B_{cc}({\bf x}_j,r)=\{{\bf y}\in \mathscr M_j: d_{cc}({\bf x}_j,{\bf y}_j)<r\}$, which made $(\mathscr M_j,d_{cc},\mu_j)$ a space of homogeneous type in the sense of Coifman and Weiss. Denote $V({\bf x}_j,r)=\mu_j(B_{cc}({\bf x}_j,r))$.

Let $\mathcal L_j$ be the sub-Laplacian on $\mathscr M_j$, $j=1,\ldots,n$.  For $0<\gamma_j<\infty$ and $\gamma=(\gamma_1,\ldots,\gamma_n)$,  the area function of $f\in L_{loc}^1(\widetilde{\mathscr M})$ is defined as
$$
S_\gamma(f)({\bf x}_1,\ldots,{\bf x}_n) =\left( \iint_{\Gamma_\gamma({\bf x}_1,\ldots,{\bf x}_n)} \Big|{Q^{[1]}_{t_1^2}\cdots Q^{[n]}_{t_n^2}    }(f)({\bf y}_1,\ldots,{\bf y}_n) \Big|^2 { {\rm d}{\mu_1(\bf y_1)} {\rm d} t_1\over V({\bf x}_1,t_1) t_1} \cdots { {\rm d}{\mu_n(\bf y_n)} {\rm d}t_n\over V({\bf x}_n,t_n) t_n} \right)^{1\over2},
$$
where 
{$ Q^{[j]}_{t_j^2} = t_j^2\mathcal L_j e^{-t_j^2\mathcal L_j}$} and $\Gamma_\gamma({\bf x}_1,\ldots,{\bf x}_n) = \Gamma^{[1]}_{\gamma_1}({\bf x}_1) \otimes\cdots\otimes\Gamma^{[n]}_{\gamma_n}({\bf x}_n)$ with 
$$ 
    \Gamma^{[j]}_{\gamma_j}({\bf x}_j) =\left \{ ({\bf y}_j,t_j)\in M_j\times \mathbb R_+:\ d_{cc}( {\bf x}_j, {\bf y}_j )< \gamma_j t_j\right \}. 
$$

 \begin{thmA}  \label{thm:app}
 There exists a positive constant $C$ such that for $p\geq 2$, and for $f\in L^p(\widetilde{\mathscr  M} )$,
\begin{equation}\label{eqn: Lp app}
   \left \|f\right\|_{L^p(\widetilde{\mathscr  M} )}\leq C p^{n\over2}\left\|S_\gamma(f)\right\|_{L^p(\widetilde{\mathscr  M} )}.
\end{equation}
\end{thmA}

\subsection{Outline} This paper is organised as follows. In Section \ref{sec:2}, we recall some preliminaries. Then we prove our main results Theorems \ref{thm:newCH-dyadic} and \ref{thm:double-newCH-dyadic}, along with their corresponding consequences Propositions \ref{coro:Lp-CH} and \ref{prop:double-Lp-CH-intro}  in Sections \ref{sec:3} and  \ref{sec:4}, respectively. We prove Theorem \ref{thm:app}  in Section \ref{sec:5} and provide an appendix in the last section.

\section{Notation and Preliminaries}\label{sec:2}
\setcounter{equation}{0}

We recall the dyadic system and Haar bases on space of homogeneous type in the sense of Coifman and Weiss. 
\subsection{Systems of dyadic cubes}
Throughout this subsection, unless otherwise specified, $\rho$ will be a quasi-metric on the space $X$, i.e., it satisfies the axioms of a metric except for the triangle inequality, which is assumed in the weaker form
\begin{equation}\label{quasi-metric}
    \rho(x,y)\leq A_0\left(\rho(x,z)+\rho(z,y)\right),\  \  \forall\, x,\,y,\, z\in X
\end{equation}
with a constant $A_0\geq 1$. Let $\mu$ be a Borel measure on $X$ which  satisfies the doubling property, namely.  that  there exists a constant
$C>0$ such that
\begin{eqnarray}\label{e2.1}
    \mu(B(x,2r))\leq C \mu(B(x, r)),\quad \forall\,r>0,\,x\in X.
\end{eqnarray}
Then there exist $C, \, n$ such that for all $\lambda\geq 1$ and $x\in X$,
\begin{equation}\label{e2.2}
    \mu(B(x, \lambda r))\leq C\lambda^n \mu(B(x,r)).
\end{equation}

We recall the following construction of dyadic cubes in $(X,\rho)$ due to T. Hyt\"onen and A. Kairema~\cite{HK}, which is
a slight elaboration of the seminal work by M. Christ \cite{Ch}.

\begin{theorem}\label{thm:dyadic}
Suppose that constants $0<c_0\leq C_0<\infty$ and $\delta\in (0,1)$ satisfy
$	12C_0\delta\leq c_0$.
Given a set of point $\{x_\beta^k\}_\beta$, $\beta\in \mathscr{J}_k$, for every $k\in \mathbb{Z}$, with the properties that
$$
    \rho(x_{\beta_1}^k, x_{\beta_2}^k) \  (\beta_1\neq \beta_2)\geq c_0 \delta^k,\   \   \     \  \min_\beta \rho(x,x_\beta^k)<C_0\delta^k    \   \   \   \text{for all  }\ x\in X,
$$	
we can construct a countable family
$\mathcal{D}=\bigcup_{k\in \mathbb{Z}} \mathcal{D}_k,\  \mathcal{D}_k=
    \left\{Q_\beta^k:\  \beta\in \mathscr{J}_k\right\}$
of Borel sets $Q_\beta^k\subseteq X$, with the properties that
\begin{enumerate}[\rm(i)]
\item $X=\bigcup_{\beta\in \mathscr{J}_k} Q_\beta^k$ (disjoint union) for all $k\in \mathbb{Z}$;
\smallskip

\item  if $\ell\geq k$, then either $Q_{\beta_2}^\ell\subseteq	 Q_{\beta_1}^k$ or  $Q_{\beta_1}^k\cap Q_{\beta_2}^\ell=\emptyset$;
\smallskip

\item  for each $(k,\beta_1)$ and each $\ell\leq k$, there exists a unique $\beta_2$ such that $Q_{\beta_1}^k\subseteq Q_{\beta_2}^\ell$;
\smallskip

\item  there exists a positive constant $\mathtt M=\mathtt M(n,\delta)$ such that for each $(k,\beta_1)$,
$$
	    1\leq \# \left\{\beta_2:\  Q_{\beta_2}^{k+1}\subseteq Q_{\beta_1}^k\right\}\leq \mathtt M  \quad  \text{and }\quad 	 Q_{\beta_1}^k= \underset{ \beta_2:\  Q_{\beta_2}^{k+1}\subseteq
         Q_{\beta_1}^k }    {\bigcup} Q_{\beta_2}^{k+1} \ ;
$$ 
\smallskip

\item  for each $(k,\beta)$, $B(x_\beta^k, c_0\delta^k/3)\subseteq Q_\beta^k\subseteq B(x_\beta^k, 2C_0\delta^k)=:B(Q_\beta^k)$ for some $x_\beta^k\in Q_\beta^k$;
\smallskip

\item  if $\ell\geq k$ and $Q_{\beta_2}^\ell\subseteq Q_{\beta_1}^k$, then $B(Q_{\beta_2}^\ell)\subseteq B(Q_{\beta_1}^k)$.
\end{enumerate}
\end{theorem}

 $\mathcal{D}$ is called a system of dyadic cubes,  the set $Q_\beta^k\in \mathcal{D}_k$ is called a dyadic cube of generation $k$ with center point $x_\beta^k$ and sidelength $\ell(Q_{\beta}^k)=\delta^k$.   For any $Q\in \mathcal{D}_{k_0}$, denote
 $Ch^{(\ell)}(Q)=\{P\in \mathcal{D}_{k_0+\ell}:\, P\subseteq Q\}$, $\ell\geq 0$. Then $Ch^{(0)}(Q)=\{Q\}$ and let $Ch(Q):=Ch^{(1)}(Q)$.

 Moreover,  we can construct {\it adjacent dyadic systems}
$$
     \left\{\mathcal{D}^{b}: b=1,2,\ldots, K\right\},
$$
where each $\mathcal{D}^b=\bigcup_{k\in \mathbb{Z}}\mathcal{D}_k^b$ is  a system of dyadic cubes satisfying properties in Theorem~\ref{thm:dyadic}, and $K=K(A_0,\delta)$ is a positive integer. In addition, for every ball $B=B(x,r)\subset X$, there exist $b$ and $Q\in \mathcal{D}^b$ with $B\subset Q$ and ${\rm diam}(Q)\leq Cr$, where $C=C(A_0, \delta)<\infty$. 
For more details, we refer the reader to \cite[Section 4]{HK}.

Throughout the article, for each $b=1,\ldots, K$,   if $^bQ_{\beta_1}^k\in \ \mathcal{D}^b$ has only one child, namely $^bQ_{\beta_2}^{k+1}$,
although $^bQ_{\beta_1}^k$ and $^bQ_{\beta_2}^{k+1}$ are the same set of points,  $^bQ_{\beta_1}^k$ and
$^bQ_{\beta_2}^{k+1}$ are seen as two different cubes. 
%
%

\subsection{Haar bases}  We recall Haar bases associated to $\mathcal D$, constructed in \cite{KLPW}.

Let $\mathtt M_Q:=\# Ch(Q)$ for each $Q\in \mathcal D$. The volume doubling property yields that there exists a positive integer $\mathtt M=\mathtt M(A_0,\delta)<\infty$ such that $1\leq \mathtt M_Q\leq \mathtt M$.

\begin{theorem}{\rm (\cite[Theorem 4.1]{KLPW})}  \label{thm:Haar-rep}
Let  $(X,\rho,\mu)$ be a quasi-metric space satisfying the volume doubling condition.
 For $1<p<\infty$, for each $f\in L^p(X,{\rm d}\mu)$, we have 
$$
   f(x)=m_X(f)+\sum_{Q\in \mathcal D}\sum_{u=1}^{\mathtt M_Q -1} \langle f, h_u^Q\rangle h_u^Q(x),
$$
where the sum converges (unconditionally) both in the $L^p(X,\mu)$-norm and pointwise $\mu$-almost everywhere, and 
$m_X(f)=f_X$ if $\mu(X)<\infty$, $m_X(f)=0$ if $\mu(X)=\infty$.
\end{theorem}

\begin{theorem} {\rm (\cite[Theorem 4.2]{KLPW})} 
\label{thm:Haar-base}
The Haar functions $h_u^Q$, $Q\in \mathcal D$, $u=1,2,\ldots, \mathtt M_Q-1$, have the following properties:
\begin{enumerate}[\rm(i)]
\item $h_u^Q$ is a simple Borel-measurable real function on $X$;

\item $h_u^Q$ is supported on $Q$;

\item $h_u^Q$ is constant on each $R\in Ch(Q)$;

\item $\int h_u^Q d\mu =0$ (cancellation);

\item $\langle h_u^Q, h_{u'}^Q \rangle =0$ for $u\neq u'$,  $u,\, u'\in \left\{1,\ldots, \mathtt M_Q-1\right\}$;

\item the collection 
$
    \left\{\mu(Q)^{-1/2}\mathsf{1}_Q \right\}\cup \left\{h_u^Q:\, u=1,\ldots, \mathtt M_Q-1\right\}
$
is an orthogonal basis for the  space of all functions on $Q$ that are constant on each sub-cube $R\in Ch(Q)$;

\item if $h_u^Q\not\equiv 0$ then
$
    \|h_u^Q\|_{L^p(X,\mu)}\approx \mu(Q_u)^{\frac{1}{p}-\frac{1}{2}} \quad \text{for\   } 1\leq p\leq \infty,
$
where $Q_u\in Ch(Q)$;

\item $\|h_u^Q\|_{L^1(X,\mu)}\cdot \|h_u^Q\|_{L^\infty (X,\mu)}\approx 1$.
\end{enumerate} 
\end{theorem}

\section{One-parameter Hoeffding inequalities:  proof of Theorem \ref{thm:newCH-dyadic}}\label{sec:3}
\setcounter{equation}{0}

This section is devoted to proving  Theorem~\ref{thm:newCH-dyadic}.
Throughout this section, let  $(X,\rho,\mu)$ be a quasi-metric space satisfying the volume doubling condition and assume that $\mu(B)<\infty$ for all balls $B\subset X$.

\smallskip

Recall that for any given $Q\in \mathcal D_{k_0}$ for some $k_0\in \mathbb Z$, the notation $Ch^{(n)}(Q)$ is given by
$$
    Ch^{(n)}(Q)=\left\{P\in \mathcal D _{k_0+n}:\,  P\subseteq Q\right\},\quad n\geq 0.
$$ 
Recall also
\begin{align}\label{FN SFN}
    F_N(x)=\sum_{\substack{P\in  Ch^{(n)}(Q)\\ 0\leq n\leq N}}   \mathfrak{a}_P(x)\quad   \text{and}\quad    
     S(F_N)(x)=\bigg(\sum_{ \substack{P\in  Ch^{(n)}(Q)\\ 0\leq n\leq N}}  \mathfrak{a}_P^2 (x)\bigg)^{1/2},
\end{align}
and similarly, we define $S(\sum \mathfrak a_P)$ as the formal square function for any  given $\sum \mathfrak a_P$.

Based on the localization and constancy of every atom $\mathfrak{a}_P(x)$ therein, one may exploit the almost-orthogonality of atoms to give an iterative algorithm associated to scales of the square function $S$, which plays a crucial role  to deduce our  Chernoff--Hoeffding bound  \eqref{eqn:newCH-dyadic}.

 \smallskip
 
\
Concretely, the organization of this section is as follows:
 
 \begin{itemize}
 \item We begin by proving an alternate version of  Hoeffding's  inequality \eqref{eqn:newCH-dyadic}: there exist positive constants $C$ and $C'$ independent of $N$ and $Q$, such that for any  $\alpha>0$,
\begin{equation}\label{eqn:newCH-dyadic-2}
    \int_Q \exp\left( \,\alpha \big|F_N (x)  \big|\,\right) \, {\rm d}\mu(x)\leq C \exp\bigg(C' \alpha^2 \,\big\|\,  S(F_N)(\cdot)  \big\|_{L^\infty(Q)}^2 \bigg)\,\mu(Q).
 \end{equation}
 
 \item As a consequence of \eqref{eqn:newCH-dyadic-2}, we obtain the one-parameter Chang--Wilson--Wolff result  \eqref{eqn:newCWW-one-parameter} in Theorem~\ref{thm:newCH-dyadic}.
 
 \item Then we show Proposition \ref{coro:Lp-CH} based on \eqref{eqn:newCWW-one-parameter} by adapting the ideas  in \cite{FP}.
  \item Finally, we prove \eqref{eqn:newCH-dyadic} in Theorem~\ref{thm:newCH-dyadic} by applying Proposition \ref{coro:Lp-CH}. 
 
 \end{itemize}
 
 That is, the outline of our proof is as follows.
 $$
      \eqref{eqn:newCH-dyadic-2}  \Longrightarrow  \eqref{eqn:newCWW-one-parameter}  \text{ \ in Theorem~\ref{thm:newCH-dyadic}\ } \Longrightarrow \text{Proposition \ref{coro:Lp-CH} }   \Longrightarrow \eqref{eqn:newCH-dyadic}   \text{ \ in Theorem~\ref{thm:newCH-dyadic}}.
 $$
 
Obviously,  \eqref{eqn:newCH-dyadic}  implies \eqref{eqn:newCH-dyadic-2}. Therefore, all the above estimates are equivalent. Namely, 
$$
      \eqref{eqn:newCH-dyadic-2}  \Longleftrightarrow  \eqref{eqn:newCWW-one-parameter}    \text{ \ in Theorem~\ref{thm:newCH-dyadic}\ } \Longleftrightarrow\text{Proposition \ref{coro:Lp-CH} }   \Longleftrightarrow  \eqref{eqn:newCH-dyadic}  \text{ \ in Theorem~\ref{thm:newCH-dyadic}}.
 $$

 
\medskip

We start our arguments.

\smallskip

 Proof of  \eqref{eqn:newCH-dyadic-2}:

To begin with,  we define 
\begin{align}\label{set G}
    \mathcal G:=\bigcup_{0\leq n\leq N}Ch^{(n)}(Q),
\end{align}
and
\begin{align}\label{set GI GII}
    \mathcal G^{(I)}:=\big\{P\in \mathcal G:\  \|\mathfrak{a}_P\|_{L^\infty} <1\big\}, \qquad \mathcal G^{(II)}:=\big\{P\in \mathcal G:\  \|\mathfrak{a}_P\|_{L^\infty} \geq1\big\}.
\end{align}
Then we split
\begin{equation}\label{eqn:FN-I-II}
    F_N(x)=\sum_{P\in \mathcal G^{(I)}}   \mathfrak{a}_P(x)+\sum_{P\in \mathcal G^{(II)}}  \mathfrak{a}_P(x)=:F_{N}^{(I)}(x)+F_{N}^{(II)}(x).
\end{equation}

 \smallskip

\noindent{\bf Step I.}  We first prove  \eqref{eqn:newCH-dyadic-2} for $F_{N}^{(II)}$ with $\alpha=1$.

\smallskip

We claim that there exists a large constant $C>1$ such that for any $1\leq k\leq N$,
 \begin{equation}\label{eqn:FN-II}
 	   \int_Q \exp\Big ( | F_k^{(II)} (x)|-C  S^2(F_{k}^{(II)})(x)\Big )\, {\rm d}\mu(x)\leq  \int_Q \exp\Big ( |F_{k-1}^{(II)} (x)|-C  S^2(F_{k-1}^{(II)})(x)\Big )\, {\rm d}\mu(x),
\end{equation}
where we define
$$F_{k}^{(II)}(x):= \sum_{\substack{P\in  Ch^{(n)}(Q):\, 0\leq n\leq k\\ \|\mathfrak{a}_P\|_{L^\infty} \geq 1  }}  \mathfrak{a}_P(x),\qquad 0\leq k\leq N;$$
and define
$$
    F_{k}^{(II)}(x):=0
$$
if for all $0\leq n\leq k$ and for all $P\in  Ch^{(n)}(Q)$, we always have $\|\mathfrak{a}_P\|_{L^\infty} < 1$.

We now show \eqref{eqn:FN-II}. Indeed,  for any $R\in Ch^{(k)}(Q)$ with $\|\mathfrak{a}_R\|_{L^\infty} \geq 1$, we have $F_{k}^{(II)}=F_{k-1}^{(II)}+\mathfrak a_R$ and $S^2(F_{k}^{(II)})=S^2(F_{k-1}^{(II)})+\mathfrak a_R^2$ on $R$, moreover, $F_{k-1}^{(II)}$ and $S^2(F_{k-1}^{(II)})$ are both constants on $R$. Note that for any $R\in Ch^{(k)}(Q)$ with $\|\mathfrak{a}_R\|_{L^\infty} < 1$, we simply have $F_{k}^{(II)}=F_{k-1}^{(II)}$ on $R$.

Denote the center of $R$ by $c_R$.  For any $R\in Ch^{(k)}(Q)$ with $\|\mathfrak{a}_R\|_{L^\infty} \geq 1$, we have
  \begin{align*}
     \int_R 	 \exp\Big ( | F_k^{(II)} (x)|-C  S^2(F_{k}^{(II)})(x)\Big )\, {\rm d}\mu(x) 
     \leq \exp\Big (|F_{k-1}^{(II)}(c_R)|-C S^2(F_{k-1}^{(II)})(c_R)\Big )\, 	\mathcal I_R(C)  ,
 \end{align*}
where
 \begin{align}\label{eqn:IR}
 	\mathcal I_R(C):=&\sum_{P\in Ch(R): \,|\mathfrak a_R (c_P)|<1}\int_{P} \exp\left(|\mathfrak a_R(x)|-C|\mathfrak a_R(x)|^2 \right) {\rm d}\mu(x) \nonumber\\
   &\qquad \qquad +\sum_{P\in Ch(R): \,|\mathfrak a_R (c_P)|\geq 1}\int_{P} \exp\left(|\mathfrak a_R(x)|-C|\mathfrak a_R(x)|^2 \right) {\rm d}\mu(x).
 \end{align}
Recall again by \textbf{(P2)},  $\mathfrak{a}_R$ is constant on each subcube $P$ in $Ch(R)$. Hence, in the first term of the right-hand side of \eqref{eqn:IR}, we see that $|\mathfrak{a}_R(x)|=|\mathfrak{a}_R(c_P)|<1$, while in the second term $|\mathfrak{a}_R(x)|=|\mathfrak{a}_R(c_P)|\geq1$.
Then considering the quadratic equation $t-Ct^2$, we see that for $t\in[0,1)$, $t-Ct^2\leq {1\over 4C}$, while for $t\geq1$, $t-Ct^2\leq1-C$. Hence, applying $t=|\mathfrak{a}_R(x)|$ in the right-hand side of \eqref{eqn:IR}, we have
$$
     \mathcal I_R(C)\leq \exp\Big(\frac{1}{4C}\Big )\sum_{P\in Ch(R): \,|\mathfrak a_R (c_P)|<1} \mu(P) + \exp(1-C)\sum_{P\in Ch(R): \,|\mathfrak a_R (c_P)|\geq 1}\mu(P)=:\widetilde{\mathcal I}_R(C).
$$
Since $\mu $ is doubling,  we see that there exists a positive constant $\theta=\theta(A_0,\delta)\in (0,1)$ independent of $R$ such that 
\begin{equation}\label{eqn:double-theta}
	\mu(P)\geq \theta \,\mu(R)\quad \text{ for\  any  \ }   P\in Ch(R).
\end{equation}
 Then it follows from $\exp\Big(\frac{1}{4C}\Big )>1>\exp(1-C)$ that $\widetilde{\mathcal I}_R(C)$
 takes its maximum when there is only one child $P\in Ch(R)$ such that $|\mathfrak a_R(c_P)|\geq 1$. Denote this particular child by $P_{\max}$.
 
 If $\# Ch(R)=1$,  then $\widetilde{\mathcal I}_R(C)<\mu(R)$ since $\exp(1-C)< 1$.
 
Otherwise, we have
 \begin{align*}
     \widetilde{\mathcal I}_R(C) \leq \left[\exp\Big(\frac{1}{4C}\Big )   \left(1-\frac{\mu(P_{\max})}{\mu(R)}\right ) + \exp(1-C) \frac{\mu(P_{\max})}{\mu(R)}\right]\mu(R),
 \end{align*}
 and the larger  $\displaystyle\frac{\mu(P_{\max})}{\mu(R)} $ is, then  the smaller $\widetilde{\mathcal I}_R(C)$ is.  Thus the maximum of $\widetilde{\mathcal I}_R(C)$ is taken when $\displaystyle \frac{\mu(P_{\max})}{\mu(R)} = \theta$.  Hence
 $$
     \mathcal I_R(C)\leq \widetilde{\mathcal I}_R(C)\leq \left[\exp\Big(\frac{1}{4C}\Big )\, (1-\theta ) + \exp(1-C) \,\theta\right] \mu(R).
 $$

 Note that  $\exp\Big(\frac{1}{4C}\Big )\to 1$ and $\exp(1-C)\to 0$ as  $C\to +\infty$, hence there exists a constant $C=C(\theta)>0$ such that  
 $$
     \mathcal I_R(C)\leq \mu(R).
 $$
 Thus
 $$
      \int_R 	 \exp\Big ( | F_k^{(II)} (x)|-C  S^2(F_{k}^{(II)})(x)\Big )\, {\rm d}\mu(x)\leq \int_R 	 \exp\Big ( | F_{k-1}^{(II)} (x)|-C  S^2(F_{k-1}^{(II)})(x)\Big )\, {\rm d}\mu(x).
 $$ 
Summing over $R\in Ch^{(k)}(Q)$ yields \eqref{eqn:FN-II} for each $1\leq k\leq N$. 

In the end, we consider $F_0^{(II)} (x)$ on $Q$. Note that in this case, $F_0^{(II)} (x)$ is either $\mathfrak a_Q (x)$ or zero. Hence, $S^2(F_{0}^{(II)})(x)$ is either $|\mathfrak a_Q (x)|^2$ or zero, respectively.

For both cases, it is direct that for the above fixed $C=C(\theta)$, we have 
\begin{align*}
	      &\int_Q 	 \exp\Big ( | F_0^{(II)} (x)|-C  S^2(F_{0}^{(II)})(x)\Big )\, {\rm d}\mu(x)  \leq  \mu(Q).
\end{align*}
This, in combination with \eqref{eqn:FN-II}, implies that 
\begin{equation}\label{eqn:newCH-dyadic-2-II}
	\int_Q 	 \exp\Big ( | F_N^{(II)} (x)|-C  S^2(F_N^{(II)})(x)\Big )\, {\rm d}\mu(x)\leq \mu(Q),
\end{equation}
that is,  the inequality \eqref{eqn:newCH-dyadic-2}  for $F_{N}^{(II)}$  with $\alpha=1$ holds. Before considering the general $\alpha>0$, we turn to the estimates on $F_N^{(I)}$.

\smallskip
 
 \noindent{\bf Step II.}  We now prove  \eqref{eqn:newCH-dyadic-2} for $F_{N}^{(I)}$, also assuming that $\alpha=1$.

 \smallskip
 

 
 \smallskip
 
The term $F_{N}^{(I)}$ sums over atomic functions $\mathfrak{a}_P$ with small $L^\infty$ norm, which is where the main difficulty for proving \eqref{eqn:newCH-dyadic-2} lies.  
\smallskip

 We begin by noting that, for the construction and selection of the dyadic cubes below, we may not have the consecutive chain of nested cubes, as some cubes $P$ went into the collection for $F_{N}^{(II)}$, handled in {\bf Step I}. 
To ease the burden of notational complexity for the lacunary nested cubes, we may think of 
$a_P\equiv 0$ for those $P$ belongs to $F_{N}^{(II)}$.

\smallskip

Firstly we claim that \eqref{eqn:newCH-dyadic-2} with $\alpha=1$ holds for $F_{N}^{(I)}$ when $\|S(F_N^{(I)})\|_{L^\infty}^2 \leq 10 $. That is,  there exist  constants $\mathsf M=\mathsf M(C_1)>1$ (here $C_1$ is the constant in ${\bf (P3)}$) and $0<c<1$ independent of $N$ and $Q$ such that 
\begin{equation}\label{eqn:k-C-H-base}
    \frac{1}{\mu(Q)}\int_Q \exp\left(  c\left| F_N^{(I)}(x)  \right|\,\right)  {\rm d}\mu(x)\leq \mathsf M \quad {\rm whenever}\   \  \left\|S(F_N^{(I)})\right\|_{L^\infty (Q)}^2\leq 10. 
\end{equation}

\smallskip

\quad {\bf Step II.1.}   Assume \eqref{eqn:k-C-H-base} at the moment  (we leave the proof to  {\bf Step II.2} below). 
We now prove \eqref{eqn:newCH-dyadic-2} with $\alpha=1$ under the condition $\left\|S(F_N^{(I)})\right\|_{L^\infty (Q)}^2>10$.

The key idea is to use \eqref{eqn:k-C-H-base}  as a base case for constructing a suitable iterative argument.

We note again that all the arguments below are valid for the lacunary version (i.e., when some $\mathfrak a_{P}$ in $F_N$ is taken to be zero). 

{Our strategy is to construct the following decomposition 
$$
   \mathcal G^{(I)}=\sum_{k} \mathcal G^{(I)}_k \quad \text{and}\quad  F_N^{(I,k)}:=\sum_{P\in \mathcal G^{(I)}_k}\mathfrak a_{P},
$$
satisfying
\begin{itemize}
	\item $k=1,2,\ldots, \mathsf  K$, where $\mathsf K\approx \|S(F_N^{(I)})\|_{L^\infty(Q)}^2$;
	\item for each $k$ and $x\in Q$, 
	$S^2(F_N^{(I, k)})(x)\leq10$, so we can apply the result \eqref{eqn:k-C-H-base} as the base case;
  \item for any $P\in \mathcal G^{(I)}_k$ and $P'\in \mathcal G^{(I)}_{k+1}$, either $P\supsetneqq P'$ or $P\cap P'=\emptyset$. Hence $F_N^{(I,k)}$ and $S(F_N^{(I,k)})$ are constant on each cube in  $\mathcal G^{(I)}_{k+1}$. 
 \end{itemize}
Based on this construction, the inequality \eqref{eqn:newCH-dyadic-2} follows from an iterative argument.
}

Concretely, we begin by constructing $\mathcal G^{(I)}_1$ via a stopping-time argument.
Let $\left\{P^{(1)}_j\right\}_j$ be  dyadic cubes in $\mathcal G^{(I)}$
 such that 
$$
  \sum_{P\in \mathcal G^{(I)},\, P\supsetneqq P'} |\mathfrak a_{P}|^2+   \max_x |\mathfrak a_{P'}(x)|^2  <9 \ \     \text{on}\  \    P',\   \  \forall\, P'\in\mathcal G^{(I)} \  \   \text{with}\    \    P'\supsetneqq P_j^{(1)}
  \leqno{(C1)}   
$$
and that 
 $$
\sum_{P\in \mathcal G^{(I)},\, P\supsetneqq P^{(1)}_j} |\mathfrak a_{P}|^2+  \max_x |\mathfrak a_{P^{(1)}_j }(x)|^2 \geq 9 \ \     \text{on}\  \     P^{(1)}_j .   \leqno{(C2)}
$$

We first note that this collection  $\left\{P^{(1)}_j\right\}_j $ is non-empty, according to the conditions that $\left\|S(F_N^{(I)})\right\|_{L^\infty (Q)}^2> 10$ and that $\|\mathfrak{a}_P\|_{L^\infty} <1$ for any $P\in \mathcal G^{(I)}$.

Next, we note that by \textbf{(P2)}, the condition (C1) shows  that 
$$
     \sum_{P\in \mathcal G^{(I)},\, P\supsetneqq P'} \left|\mathfrak a_{P}(c_{P'})\right|^2+    \left \|\mathfrak a_{P'} \right\|_{L^\infty(P')}^2  <9,\   \  \forall\, P'\in\mathcal G^{(I)} \  \   \text{with}\    \    P'\supsetneqq P_j^{(1)},
$$
where $c_{P'}$ is the center of $P'$. 

Below, we provide the details for the selection algorithm for these cubes. 

We consider the chain as follows.  
Consider each  $R\in Ch^{(N)}(Q)$ with $\|S(F_N^{(I)})\|_{L^\infty(R)}^2>10$ \footnote{Note that for each bottom $R\in Ch^{(N)}(Q)$ satisfying $\|S(F_N)\|_{L^\infty(R)}^2\leq 10$, one will see that all cubes in the chain of $R$ will belong to  $\mathcal G_1$.}   and list all cubes in $\mathcal G^{(I)}$ containing $R$ as:
$$
    \text{Chain}(R):\quad    Q=P_0\supseteq P_1\supseteq \cdots\supseteq P_N=R,
$$
where $P_i\in Ch^{(i)}(Q)$ for $0\leq i\leq N$.
Then starting from $i=0$ (indeed one may start from $i\geq 9$ due to the fact that $\|\mathfrak{a}_P\|_{L^\infty} <1$ for any $P\in \mathcal G^{(I)}$), we check whether or not $P_i$ satisfies both (C1) and (C2). If yes, we select such $P_i$; if not, then continue to check whether $P_{i+1}$ satisfies both (C1) and (C2).   
Then we find {the first $P_k$} satisfying 
$$
    \sum_{i=0}^{k-1}\left|\mathfrak{a}_{P_i}\right|^2+\max _x\left|\mathfrak{a}_{P_k}(x)\right|^2 \geq  9 \quad \text{on}\  \  P_k.
$$
The existence of $P_k$ follows from $\|S(F_N^{(I)})\|_{L^\infty(R)}^2>10$.

Eventually, this algorithm ensures that  there is \textbf{a unique} $i_0$ with $9< i_0\leq N$ such that $P_{i_0}^{(1)}=:P_R$ in this chain matching (C1) and (C2) is the cube to be selected.

\begin{figure}
\tikzset{every picture/.style={line width=0.75pt}}        

\begin{tikzpicture}[x=0.75pt,y=0.75pt,yscale=-1,xscale=1]
 
\draw  [fill={rgb, 255:red, 248; green, 231; blue, 28 }  ,fill opacity=1 ] (146,22.33) -- (332,22.33) -- (332,208.33) -- (146,208.33) -- cycle ;
 
\draw   (241.67,222) -- (331,222) -- (331,311.33) -- (241.67,311.33) -- cycle ;
 
\draw   (251.33,360.5) -- (319,360.5) -- (319,428.17) -- (251.33,428.17) -- cycle ;
 
\draw   
[fill={rgb, 255:red, 0; green, 255; blue, 255 }  ,fill opacity=1 ]
 (259,439) -- (310,439) -- (310,490) -- (259,490) -- cycle ;   
 
\draw   (264,499) -- (304,499) -- (304,539) -- (264,539) -- cycle ;
 
\draw   (269,584) -- (296,584) -- (296,611) -- (269,611) -- cycle ;
 
\draw   (346,414.5) .. controls (350.67,414.5) and (353,412.17) .. (353,407.5) -- (353,272.67) .. controls (353,266) and (355.33,262.67) .. (360,262.67) .. controls (355.33,262.67) and (353,259.34) .. (353,252.67)(353,255.67) -- (353,117.84) .. controls (353,113.17) and (350.67,110.84) .. (346,110.84) ;
 
\draw    (470,391.84) -- (470,439.84) ;
\draw [shift={(470,441.84)}, rotate = 270] [color={rgb, 255:red, 0; green, 0; blue, 0 }  ][line width=0.75]    (10.93,-3.29) .. controls (6.95,-1.4) and (3.31,-0.3) .. (0,0) .. controls (3.31,0.3) and (6.95,1.4) .. (10.93,3.29)   ;
 
\draw    (470,359.84) -- (470,311.84) ;
\draw [shift={(470,309.84)}, rotate = 90] [color={rgb, 255:red, 0; green, 0; blue, 0 }  ][line width=0.75]    (10.93,-3.29) .. controls (6.95,-1.4) and (3.31,-0.3) .. (0,0) .. controls (3.31,0.3) and (6.95,1.4) .. (10.93,3.29)   ;
 
\draw [line width=1.5]  [dash pattern={on 5.63pt off 4.5pt}]  (283,319) -- (283,334.5) -- (283,352.5) ;
 
\draw [line width=1.5]  [dash pattern={on 5.63pt off 4.5pt}]  (283,546) -- (283,561.5) -- (283,579.5) ;
 
\draw  [fill={rgb, 255:red, 248; green, 231; blue, 28 }  ,fill opacity=1 ] (145.67,222) -- (235,222) -- (235,311.33) -- (145.67,311.33) -- cycle ;
 
\draw  [fill={rgb, 255:red, 248; green, 231; blue, 28 }  ,fill opacity=1 ] (158.33,360.41) -- (226,360.41) -- (226,428.07) -- (158.33,428.07) -- cycle ;
 
\draw  [fill={rgb, 255:red, 248; green, 231; blue, 28 }  ,fill opacity=1 ] (166,438.91) -- (217,438.91) -- (217,489.91) -- (166,489.91) -- cycle ;
 
\draw  [fill={rgb, 255:red, 248; green, 231; blue, 28 }  ,fill opacity=1 ] (171,498.91) -- (211,498.91) -- (211,538.91) -- (171,538.91) -- cycle ;
 
\draw  [fill={rgb, 255:red, 248; green, 231; blue, 28 }  ,fill opacity=1 ] (176,583.91) -- (203,583.91) -- (203,610.91) -- (176,610.91) -- cycle ;
 
\draw [fill={rgb, 255:red, 248; green, 231; blue, 28 }  ,fill opacity=1 ][line width=1.5]  [dash pattern={on 5.63pt off 4.5pt}]  (190,318.91) -- (190,334.41) -- (190,352.41) ;
 
\draw [line width=1.5]  [dash pattern={on 5.63pt off 4.5pt}]  (190,545.91) -- (190,561.41) -- (190,579.41) ;
 
\draw   (137,108.84) .. controls (132.33,108.83) and (130,111.16) .. (129.99,115.83) -- (129.51,344.33) .. controls (129.5,351) and (127.16,354.32) .. (122.49,354.31) .. controls (127.16,354.32) and (129.48,357.66) .. (129.46,364.33)(129.47,361.33) -- (128.99,592.83) .. controls (128.98,597.5) and (131.3,599.83) .. (135.97,599.84) ;

\draw (516,190) node  [font=\small] [align=left] {\begin{minipage}[lt]{143pt}\setlength\topsep{0pt}
large-scale Chain(R)\\[8pt]   $\|S(F^{(I)}_N)\| _{L^{\infty }( R)}^{2}  >10$\\[15pt]
\color{red}
 We put $  P_{0} (=Q),\ P_{1} ,\dotsc ,\ P_{1}{}_{3}$ in $\mathcal{G}^{(I)}_1$, but NOT $P_{14}$ or below \\\\\\
\end{minipage}};
 
\draw (322,272) node   [align=left] {\begin{minipage}[lt]{68pt}\setlength\topsep{0pt}
$  P_{1}$
\end{minipage}};
 
\draw (319,466) node   [align=left] {\begin{minipage}[lt]{68pt}\setlength\topsep{0pt}
$  P_{13}$
\end{minipage}};
 
\draw (320,520) node   [align=left] {\begin{minipage}[lt]{68pt}\setlength\topsep{0pt}
$  P_{14}$
\end{minipage}};
 
\draw (325,630) node  [font=\small] [align=left] {\begin{minipage}[lt]{97pt}\setlength\topsep{0pt}
$ P_{N} =R\in Ch^{( N)}( Q)$
\end{minipage}};
 
\draw (545,265) node  [font=\footnotesize] [align=left] {\begin{minipage}[lt]{261pt}\setlength\topsep{0pt}
(C1) \ $ \forall \ 1\leq k\leq 12$,  \   \   $\displaystyle \sum _{i=0}^{k-1} |a_{P_{i}} |^{2} +\max |a_{P_k}( x) |^{2} < 9\ \text{on}\ P_{k}$\\
\end{minipage}};
 
\draw (500,466) node  [font=\footnotesize] [align=left] {\begin{minipage}[lt]{197.2pt}\setlength\topsep{0pt}
(C2) \ \ $\displaystyle \sum _{i=0}^{12} |a_{P_{i}} |^{2} +\max |a_{P_{13}}( x) |^{2} \geq 9\ \text{on}\ P_{13}$
\end{minipage}};
 
\draw (246,128) node   [align=left] {\begin{minipage}[lt]{68pt}\setlength\topsep{0pt}
$  P'_{0} =  Q =P_{0}$
\end{minipage}};
 
\draw (320,400) node   [align=left] {\begin{minipage}[lt]{68pt}\setlength\topsep{0pt}
$  P_{12}$
\end{minipage}};
 
\draw (484,376) node  [font=\small] [align=left] {\begin{minipage}[lt]{95.2pt}\setlength\topsep{0pt} 
 The choice of $  P_{13} $
\end{minipage}};
 
\draw (570,420) node  [font=\small] [align=left] {\begin{minipage}[lt]{132.6pt}\setlength\topsep{0pt}
\end{minipage}};
 
\draw (226,466) node   [align=left] {\begin{minipage}[lt]{68pt}\setlength\topsep{0pt}
$\displaystyle P'_{13}$
\end{minipage}};
 
\draw (227,520) node   [align=left] {\begin{minipage}[lt]{68pt}\setlength\topsep{0pt}
$\displaystyle P'_{14}$
\end{minipage}};
 
\draw (227,400) node   [align=left] {\begin{minipage}[lt]{68pt}\setlength\topsep{0pt}
$\displaystyle P'_{12}$
\end{minipage}};
 
\draw (228,272) node   [align=left] {\begin{minipage}[lt]{68pt}\setlength\topsep{0pt}
$\displaystyle P'_{1}$
\end{minipage}};
 
\draw (182,630) node  [font=\small] [align=left] {\begin{minipage}[lt]{97pt}
\setlength\topsep{0pt}
$\displaystyle P'_{N} =R^{'} \in Ch^{( N)}( Q)$
\end{minipage}};
 
\draw (60,285) node   
[align=left] {\begin{minipage}[lt]{104pt}\setlength\topsep{0pt}\small 
small-scale $\text{Chain}(R')$\\[8pt]   $\|S(F_N^{(I)})\| _{L^{\infty }( R')}^{2} \leq 10$\\[80pt]
 We put all cubes in  \\[3pt]  $\text{Chain}(R')$ in $  \mathcal{G}^{(I)}_{1}$,  \\[3pt] that is, \\[3pt] $
P'_0 (=Q),\, P'_1 , \,\ldots, \\[3pt]
P'_{13} ,\,  P'_{14} ,\,\ldots,\\[3pt]
P'_{N} =R'
$  
 are all in $\displaystyle \mathcal{G}^{(I)}_{1}$
\end{minipage}};
\end{tikzpicture}
\caption{ The example illustrating the selection of  cubes in the small-scale $\text{Chain}(R')$ and 
large-scale $\text{Chain}(R)$, respectively.}
  \label{fig:1}
\end{figure}
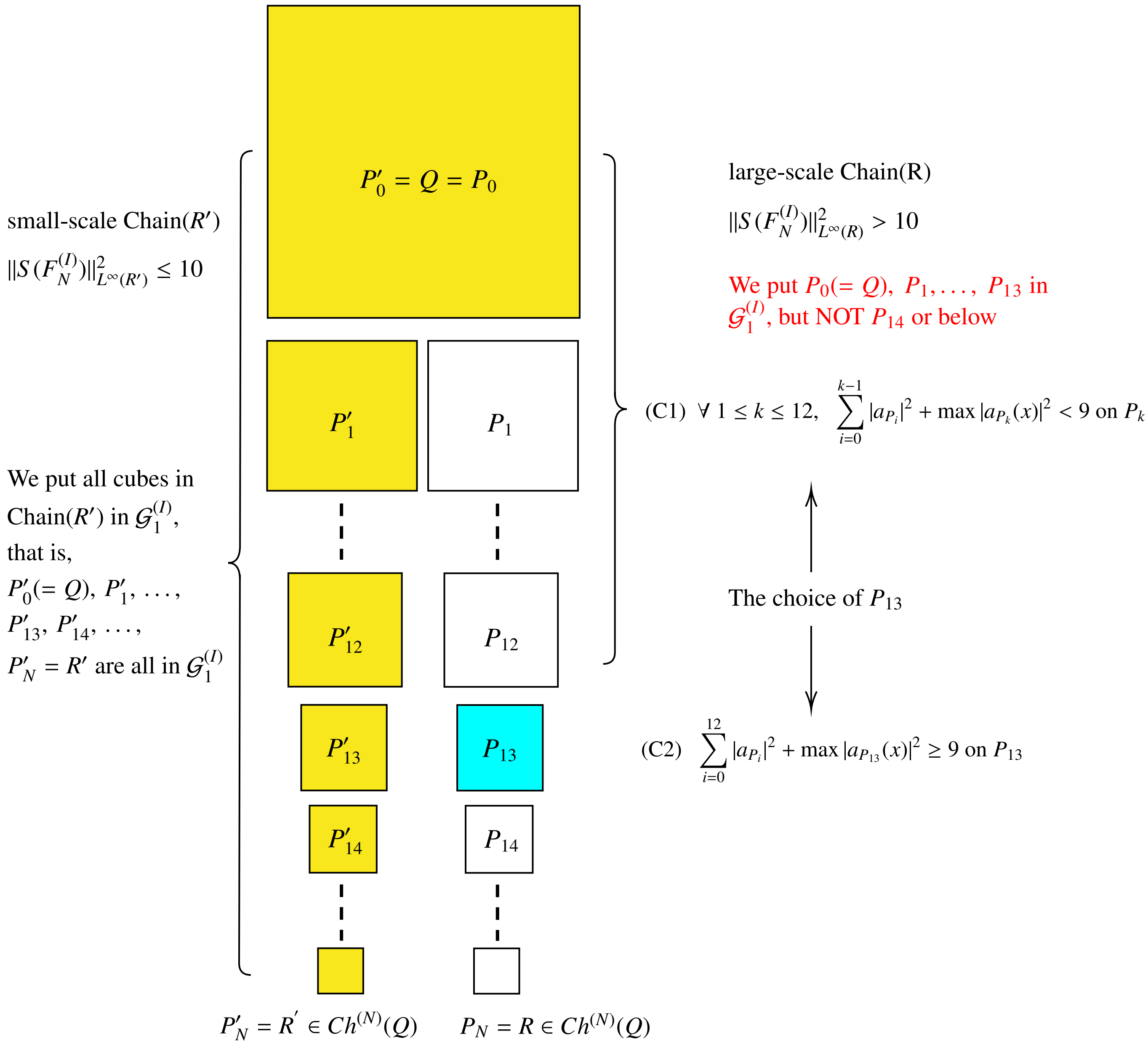

See the second  (on the right-hand side) chain $\text{Chain}(R)$  in Figure \ref{fig:1}
 with $P_{i_0}^{(1)}=P_{13}$ for instance.

Additionally, for any  $R_1,R_2\in Ch^{(N)}(Q)$ with $\|S(F_N^{(I)})\|_{L^\infty}^2>10$ on both $R_1$ and $R_2$, we must have either $P_{R_1}\cap P_{R_2}=\emptyset$ or $P_{R_1}= P_{R_2}$. We emphasize that it's impossible that $P_{R_1}\supsetneqq  P_{R_2}$ or $P_{R_1}\subsetneqq  P_{R_2}$. We argue by contradiction: if $P_{R_1}\supsetneqq  P_{R_2}$,  they both must occur in the chain of  $R_2$, i.e.,
$$
   \text{Chain}(R_2):\quad   Q=P_0\supseteq P_1\supseteq \cdots \supseteq P_{R_1}\supseteq\cdots\supseteq P_{R_2}\supseteq \cdots \supseteq R_2.
$$
But the algorithm requires that $P_{R_1}$ be chosen as the critical cube in this chain, that is,  $P_{R_2}= P_{R_1}$.

Therefore, $\left\{P^{(1)}_j\right\}_j$ is maximal in the sense of 
\begin{itemize}
	\item every $P^{(1)}_j$ is the first  
	cube satisfying  (C2) as a member of the chain in $\mathcal G^{(I)}$ containing $R\in Ch^{(N)}(Q)$.
	\item these $P^{(1)}_j$'s are mutually disjoint.
\end{itemize}
The maximality ensures the non-conflict and existence of $\mathcal G^{(I)}_1$  and $F_N^{(I,1)}$. Moreover,  observe that

(i). if $x\in P^{(1)}_j$, combining the conditions (C1), (C2) and $\|\mathfrak a_{ P^{(1)}_j}\|_{L^\infty}<1$, we have
\begin{equation}\label{eqn:range-chain}
   8< \sum_{P\in \mathcal G^{(I)},\, P\supseteq  P^{(1)}_j} |\mathfrak a_{P}|^2<10 \ \     \text{on}\  \     P^{(1)}_j .   
\end{equation}

(ii).  for $x\in Q\setminus \big(\bigcup_j P^{(1)}_j \big )$, there exists $R\in Ch^{(N)}(Q)$ such that $x\in R$ and $\|S(F_N^{(I)})\|_{L^\infty(R)}^2\leq 10$.

Now we define
$$
    \mathcal G_1^{(I)}:=\left\{P\in \mathcal G^{(I)}:\, P\supseteq P^{(1)}_j \  \text{for \  some\ } j   \ \text{or}\  \   P\cap  \big(\cup_j P^{(1)}_j \big )=\emptyset  
    \right\}\quad \text{and}\quad  F_N^{(I,1)}:=\sum_{P\in \mathcal G_1^{(I)}} \mathfrak a_P.
$$
We note that the condition $P\cap  \big(\cup_j P^{(1)}_j \big )=\emptyset$  above means ``up to measure zero'' with respect to the intersection of boundaries of cubes.
Moreover, if $P$ satisfies $P\cap  \big(\cup_j P^{(1)}_j \big )=\emptyset$, then it shows that $P$ belongs to some $\text{Chain}(R')$ of $R'\in Ch^{(N)}(Q)$ such that 
$ \|S(F_N^{(I)})\|^2_{L^\infty(R')}\leq10  $.
That is, $\text{Chain}(R')$ is on a small scale with respect to the square function.  See the first (on the left-hand side) chain  $\text{Chain}(R')$  in Figure \ref{fig:1}
 with $ \|S(F_N^{(I)})\|^2_{L^\infty(R')}\leq 10  $ for instance.
Therefore,
$$
    S^2(F_N^{(I,1)})(x)\leq10  \  \    \text{for}\  \    x\in Q.
$$

Inductively, using the same algorithm for $k\geq 2$, we can construct the set $\mathcal G_k^{(I)}$ from $\mathcal G^{(I)}\setminus \big(\bigcup_{i\leq k-1} \mathcal G^{(I)}_i\big)$ and the corresponding $F_N^{(I,k)}$ such that $\mathcal G^{(I)}=\bigcup_k \mathcal G^{(I)}_k$ and $F_N^{(I)}=\sum_k F_N^{(I,k)}$. Explicitly, let $\left\{P^{(k)}_j\right\}_j$ be the dyadic cubes in $\mathcal G^{(I)}\setminus \big(\bigcup_{i\leq k-1} \mathcal G^{(I)}_i\big)$ such that

$$
  \sum_{P\in \mathcal G^{(I)}\setminus (\bigcup_{i\leq k-1} \mathcal G^{(I)}_i),\, P\supsetneqq P'} |\mathfrak a_{P}|^2+   \max_x |\mathfrak a_{P'}(x)|^2  <9 \ \     \text{on}\  \    P',\   \  \forall\, P'\in\mathcal G^{(I)}\setminus \Big(\bigcup_{i\leq k-1} \mathcal G^{(I)}_i\Big) \  \   \text{with}\    \    P'\supsetneqq P_j^{(k)}
$$
 while 
 $$
    \sum_{P\in \mathcal G^{(I)}\setminus (\bigcup_{i\leq k-1} \mathcal G^{(I)}_i) ,\, P\supsetneqq P^{(k)}_j} |\mathfrak a_{P}|^2+  \max_x |\mathfrak a_{P^{(k)}_j }(x)|^2 \geq 9 \ \     \text{on}\  \     P^{(k)}_j .  
$$
Then denote
$$
    \mathcal G^{(I)}_k:=\bigg\{P\in \mathcal G^{(I)}\setminus \Big(\bigcup_{i\leq k-1} \mathcal G^{(I)}_i\Big):\, P\supseteq P^{(k)}_j \  \text{for \  some\ } j   \ \text{or}\  \   P\cap  \big(\cup_j P^{(k)}_j \big )=\emptyset  
    \bigg\}\quad \text{and}\quad  F_N^{(I,k)}:=\sum_{P\in \mathcal G^{(I)}_k} \mathfrak a_P.
$$
Obviously, $S^2(F_N^{(I,k)})(x)\leq10$ for any $ x\in Q$. According to the classification above, we have 
$$
    1\leq k\leq \mathsf K:=\left\lfloor \frac{\|S(F_N^{(I)})\|_{L^\infty}^2}{8} \right\rfloor +1,
$$
where $\lfloor \cdot \rfloor$ denotes the floor function. This bound follows from considering the chain of some $R\in Ch^{(N)}(Q)$ such that the maximum of $S^2(F_N^{(I)})$ on $R$ is equal to $\|S^2(F_N^{(I)})\|_{L^\infty(Q)}$.

Clearly,
$$
    F_N^{(I)}=\sum_{k} F_N^{(I,k)}\quad \text{and}\quad S^2(F_N^{(I)})=\sum_k S^2(F_N^{(I,k)}).
$$

For each $k\geq 2$, let $\mathcal G^{(I)}_{k,\max}$ be the set of maximal cubes in $\mathcal G^{(I)}_k$. Then cubes in $\mathcal G^{(I)}_{k,\max}$ are mutually disjoint and the algorithm above ensures that $F_N^{(I,i)}$ is constant on each $W\in \mathcal G^{(I)}_{k,\max}$ for $1\leq i\leq k-1$. Meanwhile, if $P\in \bigcup_{1\leq i\leq k}\mathcal G^{(I)}_i$ while $P\cap (\bigcup_{W\in \mathcal G^{(I)}_{k,\max}} W)=\emptyset$, then $P\in \bigcup_{1\leq i\leq k-1}\mathcal G^{(I)}_i$. Therefore, combining $\|S(F_N^{(I,k)})\|_{L^\infty}^2\leq10$ and \eqref{eqn:k-C-H-base} to see 
\begin{align*}
	&\int_{Q}\exp\bigg( c\Big|\sum_{i=1}^{k} F_N^{(I,i)} (x) \Big|  \bigg)\,  {\rm d}\mu(x)\\
	\leq\,& \sum_{W\in \mathcal G^{(I)}_{k,\max} } \exp\bigg( c\big|\sum_{i=1}^{k-1} F_N^{(I,i)} (c_W)\big|  \bigg)\, \int_{W}\exp\Big ( c\big|F_N^{(I,k)} \big|  \Big )\, {\rm d}\mu(x) +\int_{Q\setminus \bigcup\limits_{W\in \mathcal G^{(I)}_{k,\max}}W  } \exp\bigg( c\Big|\sum_{i=1}^{k-1} F_N^{(I,i)} (x) \Big|  \bigg)\, {\rm d}\mu(x)\\
	\leq \, &  \mathsf M \sum_{W\in \mathcal G^{(I)}_{k,\max} } \exp\bigg( c\big|\sum_{i=1}^{k-1} F_N^{(I,i)} (c_W)\big|  \bigg)\,  \mu(W) +\int_{Q\setminus \bigcup\limits_{W\in \mathcal G^{(I)}_{k,\max}}W  } \exp\bigg( c\Big|\sum_{i=1}^{k-1} F_N^{(I,i)} (x) \Big|  \bigg)\, {\rm d}\mu(x)\\
	\leq \, &  \mathsf M \int_{Q}\exp\bigg( c\Big|\sum_{i=1}^{k-1} F_N^{(I,i)} (x) \Big|  \bigg)\, {\rm d}\mu(x),
\end{align*}
where $0<c<1$ is the constant in \eqref{eqn:k-C-H-base}.
By iteration $\mathsf K$ times, we obtain that whenever $\|\mathfrak a_P\|<1$ for each $P$, we have
\begin{equation}\label{eqn:newCH-dyadic-2-I}
    \int_{Q}\exp\bigg( c\big| F_N^{(I)} (x) \big|  \bigg)\,  {\rm d}\mu(x)\leq \mathsf M^{
    \mathsf K} \mu(Q) \leq \mathsf M \exp\left(\frac{\ln \mathsf M}{8} \|S(F_N^{(I)})\|_{L^\infty}^2\right) \mu(Q).
\end{equation}

\medskip

{\bf Step II.2.}  Now we start to prove \eqref{eqn:k-C-H-base}.  

One may consider the dyadic ${\rm BMO}$ norm (associated to $\mathcal D$) of $F_N^{(I)}$. For any $T\in \mathcal D$ with $T\cap Q\neq \emptyset$, if $Q\subseteq T$, it follows from \textbf{(P3)}  that 
$$
    \frac{1}{\mu(T)}\int_{T} \left| F_N^{(I)}(x)-(F_N^{(I)})_{T}\right| {\rm d}\mu(x)\leq 2|F_N^{(I)}|_{T}\leq 2|F_N^{(I)}|_{Q}\leq 2C_1 \left\|S(F_N^{(I)})\right\|_{L^\infty (Q)};
$$
otherwise, we have  $T\subseteq Q$, then by \textbf{(P2)} and \textbf{(P3)},
\begin{align*}
    \frac{1}{\mu(T)}\int_{T} \left| F_N^{(I)}(x)-(F_N^{(I)})_{T}\right| {\rm d}\mu(x)&=\frac{1}{\mu(T)} \int_{T} \bigg| \sum_{ P\in \mathcal G^{(I)}:\,   P\subseteq T } \mathfrak{a}_P (x)-\bigg(\sum_{  P\in \mathcal G^{(I)}:\, P\subseteq T }  \mathfrak{a}_P \bigg)_{T}\bigg|\, {\rm d}\mu(x)\\
    &\leq 2C_1 \bigg\|S\Big(\sum_{ P\in \mathcal G^{(I)}:\,  \,   P\subseteq T } \mathfrak{a}_P\Big)\bigg\|_{L^\infty (T)}\\
    &\leq 2C_1 \left\|S(F_N^{(I)})\right\|_{L^\infty (Q)}.
\end{align*}
Hence $F_N^{(I)}$ belongs to $  {\rm BMO}_{\mathcal D}$ with $\|F_N^{(I)}\|_{{\rm BMO}_{\mathcal D}}\leq 2C_1 \left\|S(F_N^{(I)})\right\|_{L^\infty (Q)}$,  where $C_1$ is the constant from \textbf{(P3)}.

Furthermore, we claim that applying \textbf{(P1)}, \textbf{(P2)} and \textbf{(P3)} again,  a standard dyadic argument similar to that of the John--Nirenberg type inequality shows that for any $\lambda >0$,
\begin{equation}\label{eqn:J-N}
	\mu\left(\left\{x\in Q:\, \left|F_N^{(I)} (x) - (F_N^{(I)})_Q\right|>\lambda \right\}\right)\leq c_1 \exp\left(-c_2\frac{\lambda}{\|S(F_N^{(I)})\|_{L^\infty(Q)}}\right)  \mu(Q).
\end{equation}
Indeed, by $|F_N^{(I)}|_Q\leq C_1\left  \|S(F_N^{(I)})\right \|_{L^\infty(Q)}$ again,  \eqref{eqn:J-N} follows from its variant
\begin{equation}\label{eqn:J-N-2}
	\mu\Big(\Big\{x\in {Q}:\     |F_N^{(I)}(x)|>\lambda\Big\}\Big)
    \leq   c_1 \exp\bigg(- c_2\f{ \lambda}  {\|S(F_N^{(I)})\|_{L^\infty(Q)}}\bigg)\, \mu({Q})
\end{equation}
for all $\lambda \gtrsim \left\|S(F_N^{(I)})\right\|_{L^\infty(Q)}$.  

Now we start to show \eqref{eqn:J-N-2}.
For each fixed $\lambda\geq 4C_1\|S(F_N^{(I)})\|_{L^\infty(Q)}$, consider the strict superlevel set
$$
    \mathsf S_\lambda:=\Big \{x\in {Q}:   \,  |F_N^{(I)}(x)|>\lambda\Big \}.  
$$
Without loss of generality, we may assume $\mathsf S_\lambda\neq \emptyset$.
Combining $\left|F_N^{(I)} \right|_{Q}\leq \lambda$, using a stopping time argument, one may construct
a set $\left\{Q_{\lambda,i}\right\}_i$ of maximal dyadic cubes $Q_{\lambda,i} \subsetneqq Q$ such that  
\begin{equation}\label{eqn:stopping-time-one}
	\left|F_N^{(I)} \right|_{Q_{\lambda,i}}> \lambda  \quad \text{while}\quad 
	\left|F_N^{(I)} \right|_{\widetilde{Q}_{\lambda,i}}
         \leq \lambda,
\end{equation}
where $\widetilde{Q}_{\lambda,i}$ is the parent of $Q_{\lambda,i}$, that is $Q_{\lambda,i}\in Ch(\widetilde{Q}_{\lambda,i})$. 
Hence up to a subset of measure zero (i.e., the boundaries of cubes), we have 
$$
    \mathsf S_\lambda \subseteq \bigcup_i Q_{\lambda,i} =:\mathsf Q_\lambda.
$$

To continue the argument, we need a subtle upper bound of $ \left|F_N^{(I)}\right|_{Q_{\lambda,i}}$ for any maximal cube $Q_{\lambda, i}$. By \textbf{(P1)} and \textbf{(P2)}, we have
\begin{eqnarray*} 
    \left|F_N^{(I)}\right|_{Q_{\lambda,i}} &\leq & \Bigg| \sum_{ P\in \mathcal G^{(I)}:\,  \,   P\supsetneqq \widetilde{Q}_{\lambda,i} } 
    \mathfrak{a}_P \Bigg|_{Q_{\lambda,i}} + \left|\mathfrak{a}_{\widetilde{Q}_{\lambda,i}} \right|_{Q_{\lambda,i}} +   \Bigg|
    \sum_{ P\in \mathcal G^{(I)}:\,    P\subseteq  Q_{\lambda,i} } 
    \mathfrak{a}_P \Bigg|_{Q_{\lambda,i}}  \nonumber\\
    &=&   \Bigg|  \sum_{ P\in \mathcal G^{(I)}:\,    P\supsetneqq \widetilde{Q}_{\lambda,i} } 
    \mathfrak{a}_P \Bigg|_{\widetilde{Q}_{\lambda,i}}+ \left|\mathfrak{a}_{\widetilde{Q}_{\lambda,i}}\right|_{Q_{\lambda,i}}  +\Bigg|\sum_{ P\in \mathcal G^{(I)}:\,   P\subseteq  Q_{\lambda,i} }  \mathfrak{a}_P\Bigg|_{Q_{\lambda,i}}\nonumber \\
    &=&  \Bigg|F_N^{(I)}-\sum_{ P\in \mathcal G^{(I)}:\,    P\subseteq \widetilde{Q}_{\lambda,i} }  \mathfrak{a}_P \Bigg|_{\widetilde{Q}_{\lambda,i}} + \left|\mathfrak{a}_{\widetilde{Q}_{\lambda,i}} \right|_{Q_{\lambda,i}}  + \Bigg|\sum_{ P\in \mathcal G^{(I)}:\,   P\subseteq  Q_{\lambda,i} }   \mathfrak{a}_P\Bigg|_{Q_{\lambda,i}}.
\end{eqnarray*}
Note that $\mathfrak{a}_{\widetilde{Q}_{\lambda,i}}$ is constant on $Q_{\lambda,i}$ by \textbf{(P2)}, furthermore,
$$
    \left|\mathfrak{a}_{\widetilde{Q}_{\lambda,i}}\right|_{Q_{\lambda,i}}\leq \|\mathfrak{a}_{\widetilde{Q}_{\lambda,i}}\|_{L^\infty(\widetilde{Q}_{\lambda,i})}\leq  \|S(F_N^{(I)})\|_{L^\infty(Q)}.
$$
This, in combination with the fact  $\left|F_N^{(I)}\right|_{\widetilde{Q}_{\lambda,i}}\leq \lambda$ and with \textbf{(P3)}, yields that
\begin{eqnarray}\label{eqn:auxAverg1}
    \left|F_N^{(I)}\right|_{Q_{\lambda,i}} &\leq&  \left|F_N^{(I)} \right|_{\widetilde{Q}_{\lambda,i}}+ \Bigg|\sum_{ P\in \mathcal G^{(I)}:\,   P\subseteq \widetilde{Q}_{\lambda,i} }  \mathfrak{a}_P \Bigg|_{\widetilde{Q}_{\lambda,i}}+(C_1+1 )  \left\|S(F_N^{(I)})\right\|_{L^\infty(Q)} \nonumber\\
     &\leq& \lambda +(2C_1+1) \left \|S(F_N^{(I)})\right\|_{L^\infty(Q)}.
\end{eqnarray}

\smallskip

Now, take a parameter $\lambda'>0$ which will be chosen later. As before, let $\big\{Q_{\lambda+\lambda',j}\big \}$ be the maximal dyadic cubes such that
$|F_N^{(I)}|_{Q_{\lambda+\lambda',j}}>\lambda+\lambda'$. Set $\mathsf Q_{\lambda+\lambda'}:= \bigcup_j Q_{\lambda+\lambda',j} $, and
obviously,  $\mathsf Q_{\lambda+\lambda'}\subseteq \mathsf Q_{\lambda}$.

 From the selection rule of $Q_{\lambda+\lambda',j}$ and $Q_{\lambda,i}$,
if $Q_{\lambda+\lambda',j}\cap Q_{\lambda,i}\neq \emptyset$, then
$Q_{\lambda+\lambda',j}\subseteq Q_{\lambda,i}$. For every fixed $Q_{\lambda,i}$, by using \eqref{eqn:auxAverg1} and \textbf{(P3)},
\begin{align}\label{em1}
    \mu\left(\mathsf Q_{\lambda+\lambda'}\cap Q_{\lambda,i}\right)
    \leq & \f{1}{\lambda+\lambda'} \sum_{j:\  Q_{\lambda+\lambda',j}
    \subset Q_{\lambda,i}}  \int_{Q_{\lambda+\lambda',j}}
    |F^{(I)}_N(x)| \,{\rm d}\mu(x)\nonumber\\
    \leq &  \f{1}{\lambda+\lambda'}  \int_{Q_{\lambda,i}} \left|F_N^{(I)}(x)
    -(F^{(I)}_N)_{Q_{\lambda,i}}\right|  {\rm d}\mu(x) +\f{ \left|
    (F^{(I)}_N)_{Q_{\lambda,i}} \right|}{\lambda+\lambda'}
    \mu(\mathsf Q_{\lambda+\lambda'}\cap Q_{\lambda,i})\nonumber\\
     \leq &  \f{1}{\lambda+\lambda'}  \int_{Q_{\lambda,i}} \bigg|\sum_{P\in \mathcal G^{(I)},\, P\subseteq Q_{\lambda,i}} \mathfrak a_P (x)-\bigg(\sum_{P\in \mathcal G^{(I)},\, P\subseteq Q_{\lambda,i}} \mathfrak a_P\bigg)_{Q_{\lambda,i}}
     \bigg|  {\rm d}\mu(x) +\f{ \left|
    (F_N^{(I)})_{Q_{\lambda,i}} \right|}{\lambda+\lambda'}
    \mu(\mathsf Q_{\lambda+\lambda'}\cap Q_{\lambda,i})\nonumber\\
    \leq & \f{2C_1 \left\|S(F_N^{(I)})\right\|_{L^\infty({Q})}}
    {\lambda+\lambda'} \mu(Q_{\lambda,i}) +\f{\lambda+(2C_1+1)\left\|S(F_N^{(I)})\right\|_{L^\infty({Q})}
    }{\lambda+\lambda'} \mu(\mathsf Q_{\lambda+\lambda'}\cap Q_{\lambda,i}).
\end{align}
Therefore,
$$
    \mu(\mathsf Q_{\lambda+\lambda'}\cap Q_{\lambda,i})  \leq \f{2C_1\left\|S(F_N^{(I)})
    \right\|_{L^\infty({Q})}}{\lambda'-(2C_1+1)\left\|S(F_N^{(I)})\right\|_{L^\infty({Q})} } \mu(Q_{\lambda,i}).
$$
This tells us that for every $\lambda \geq 4C_1\|S(F_N^{(I)})\|_{L^\infty(Q)}$ and $\lambda'>(2C_1+1) \|S(F_N^{(I)})\|_{L^\infty({Q})}$,
$$
    \mu(\mathsf Q_{\lambda+\lambda'}) \leq \f{2C_1\left\|S(F_N^{(I)})\right\|_{L^\infty({Q})}}
    {\lambda'-(2C_1+1) \left\|S(F_N^{(I)})\right\|_{L^\infty({Q})} }\mu(\mathsf Q_{\lambda}).
$$

Now fix $\lambda'=(8C_1+1)\left\|S(F_N^{(I)})\right\|_{L^\infty({Q})}$. It suffices to consider $\lambda \geq  \lambda'$ and there exists a unique $k\geq 1$ such that $k\lambda'\leq  \lambda < (k+1)\lambda '$. The argument above implies that
$$
    \mu(\mathsf S_\lambda)\leq \mu(\mathsf Q_\lambda)\leq \mu(\mathsf Q_{k\lambda'})\leq \frac{1}{3} \mu(\mathsf Q_{(k-1)\lambda'})< \frac{1}{e}\,\mu(\mathsf Q_{(k-1)\lambda'}) .
$$
Iterating $k-1$ times and noting that 
$$
    k-1>\frac{\lambda}{\lambda'}-2=\frac{\lambda} {(8C_1+1)\left\|S(F_N^{(I)})\right\|_{L^\infty({Q})}}-2,
$$
we have 
$$
     \mu(\mathsf S_\lambda)< \exp\left(2-\frac{\lambda} {(8C_1+1)\left\|S(F_N^{(I)})\right\|_{L^\infty({Q})}}\right) \mu(Q),
$$
thus \eqref{eqn:J-N-2} holds with 
\begin{equation}\label{eqn:constants}
	c_1\geq e^2\quad\text{and}\quad c_2\leq \frac{1}{8C_1+1}.
\end{equation}   
Consequently, the following inequality
$$
\mu\Big(\Big\{x\in {Q}:\     |F_N^{(I)}(x)|>\lambda\Big\}\Big)
    \leq   c_1 \exp\bigg(- \frac{c_2}{\sqrt{10}} \lambda\bigg)\, \mu({Q})
$$
holds  for $\lambda\geq 4C_1 \left\|S(F_N^{(I)})\right\|_{L^\infty({Q})}$, whenever $\left\|S(F_N^{(I)})\right\|_{L^\infty({Q})}^2\leq 10$. 

As a consequence, taking
$\displaystyle c=\frac{c_2}{2\sqrt{10}}$, we have for $\|S(F_N^{(I)})\|_{L^\infty}^2 \leq 10$,
\begin{align}\label{eqn:constant-M}
	&\int_Q \exp\left(c|F_N^{(I)}(x)|\right) {\rm d}\mu(x) \nonumber\\
	=&\left\{ \int_0^{4C_1\sqrt{10}} +\int_{4C_1\sqrt{10}} ^\infty   \right\} \mu\Big(\Big\{x\in Q:\, \exp\left(c|F_N^{(I)}(x)|\right)   >\exp\left(c\lambda\right)\Big\}\Big)  \,{\rm d} \exp\left(c\lambda\right)\nonumber\\
	\leq & \exp\left(4cC_1 \sqrt{10}\right) \mu(Q) +\int_{4C_1\sqrt{10}} ^\infty  c_1 \exp\bigg(- \frac{c_2}{\sqrt{10}} \lambda\bigg)   \,\mu(Q) \cdot c \exp(c\lambda)\, {\rm d}\lambda\nonumber\\
	\leq &\,\mathsf M\, \mu(Q),
\end{align}
where $\mathsf M=\mathsf M(C_1)$,
hence the base case \eqref{eqn:k-C-H-base} holds.

The above arguments complete the proof of \eqref{eqn:newCH-dyadic-2} for $F_N^{(I)}$, again for $\alpha=1$. 



\medskip
{\bf Step III.}  The proof of \eqref{eqn:newCH-dyadic-2}  for $F_N$ for every $\alpha>0$ follows by rescaling.

Combining {\bf Step I}  and {\bf Step II}, we obtain that  
$$
	\int_{Q}\exp\bigg( \frac{c}{2}\big| F_N (x) \big|  \bigg)\,  {\rm d}\mu(x)\leq \mathsf M^{1\over2} \exp\left({1\over2}\Big(\frac{\ln \mathsf M}{8} +C(\theta)\Big)\ \|S(F_N)\|_{L^\infty}^2\right) \mu(Q).
$$

Moreover, the same argument  also implies that for any $\alpha>0$,
\begin{align}\label{eqn:weak-newCH-dyadic}
	 \int_{Q}\exp\bigg(  \alpha  \big| F_N (x) \big|  \bigg)\,  {\rm d}\mu(x) &=\int_{Q}\exp\bigg( \frac{c}{2} \, \big|   \frac{2\alpha}{c}  F_N (x) \big|  \bigg)\,  {\rm d}\mu(x) \nonumber\\
	 &\leq \mathsf M^{1\over2} \exp\left({1\over2}\Big(\frac{\ln \mathsf M}{8} +C(\theta)\Big)  \Big \|S\Big ( \frac{2\alpha}{c} \cdot F_N\Big )\Big \|_{L^\infty}^2\right) \mu(Q) \nonumber\\
	 &=\mathsf M^{1\over2} \exp\left( {2\over c^2}\Big(\frac{\ln \mathsf M}{8} +C(\theta)\Big)   \alpha^2  \|S ( F_N ) \|_{L^\infty}^2\right) \mu(Q).
\end{align}
That is, \eqref{eqn:newCH-dyadic-2} holds. 
 \hfill $\Box$


\smallskip

\begin{remark}\label{rem:P3}
As for the almost orthogonality \textbf{(P3)},
the   $L^2$-norm in the left-hand side of \eqref{eqn:P3} can be weakened to the $L^1$-norm. That is, one may weaken the condition \eqref{eqn:P3} as follows:
$$
	 \frac{1}{\mu(R)}  \int_R  \Big|   \sum_{P\in \mathcal G_R} \mathfrak{a}_P(x)\Big| {\rm d}\mu(x)  \leq C_1\,  \bigg\|\,  \bigg(\sum_{ P\in \mathcal G_R}    \mathfrak{a}_P^2 (x)\bigg)^{1/2}\,\bigg\|_{L^\infty(R)} .
$$
This follows from noticing that the iterative argument in {\bf Step II.1} only involves the mean of atomic functions and their sums, and these are the only places in the entire proof where almost orthogonality is used.
\end{remark}

\smallskip

\begin{remark}\label{rem:C1C2}
Note that $\left\{P^{(1)}_j\right\}_j$ satisfying the conditions (C1)  and (C2), can be described equivalently as follows: let $\left\{P^{(1)}_j\right\}_j$ be \textbf{maximal}  dyadic cubes in $\mathcal G$ 
 such that 
$$
    \sum_{P\in \mathcal G,\, P\supsetneqq P^{(1)}_j} |\mathfrak a_{P}|^2 <9  \ \     \text{on}\  \     P^{(1)}_j   \leqno{(C1')}   
$$
while 
$$
    \sum_{P\in \mathcal G,\, P\supsetneqq P^{(1)}_j} |\mathfrak a_{P}|^2+  \max_x |\mathfrak a_{P^{(1)}_j }(x)|^2 \geq 9 \ \     \text{on}\  \     P^{(1)}_j .   \leqno{(C2)}
$$
Although replacing (C1) with (C1') may seem to weaken the condition, the assumed `maximality' of $P^{(1)}_j$ ensures the equivalence. Concretely, (C1) implies (C1') straightforward; the maximality of $P^{(1)}_j$ satisfying (C1') and (C2) deduces that (C1) holds by contradiction.
\end{remark}

\bigskip

{\bf Proof of  \eqref{eqn:newCWW-one-parameter}: }  

Now we show the implication 

\begin{center}
``\eqref{eqn:newCH-dyadic-2} $\Rightarrow$ the one-parameter Chang--Wilson--Wolff type estimate \eqref{eqn:newCWW-one-parameter}''.  
\end{center}

For any $\lambda>0$, it follows from \eqref{eqn:newCH-dyadic-2} with choosing $\displaystyle \alpha=\frac{\lambda}{2C'\|S(F_N)\|_{L^\infty(Q)}^2 }$ that
$$
    \int_Q \exp\bigg( \frac{\lambda}{2C'\|S(F_N)\|_{L^\infty(Q)}^2 }  \Big(|F_N(x)|-\frac{\lambda}{2}\Big)\, \bigg)\, {\rm d}\mu(x)\leq C \mu(Q).
 $$
As a consequence, 
$$
      \mu\left( \left\{x\in Q:\,   |F_N(x)|>\lambda \right\}\right) \leq C  \exp \bigg( - \frac{\lambda^2}{4C'\|S(F_N)\|_{L^\infty(Q)}^2 }   \, \bigg) \, \mu(Q).
$$
Therefore, we obtain \eqref{eqn:newCWW-one-parameter} readily by taking some constant $c_3$ which depends on  $C'$.
 \hfill $\Box$

\bigskip

Furthermore, \eqref{eqn:newCWW-one-parameter}, in combination with some basic results on Orlicz spaces, allows us to modify the method in \cite{FP} to obtain the sharp $L^p$ lower bound for $S(F_N)$ whenever $p\geq 2$.

\bigskip

We first record a local extrapolation lemma.  It will also be used in the proof of
Proposition~\ref{prop:double-Lp-CH}.

\begin{lemma}\label{lem:local-CWW-to-Lp}
Let $Q\in\mathcal D$, let $\mathcal H$ be a finite-dimensional Hilbert space, and let
$\{A_P\}_{P\in\mathcal G}$ be a finite family indexed by dyadic cubes $P\subseteq Q$.  Assume that $A_P$ is supported
on $P$ and is constant on every child of $P$.  Set
$$
 G=\sum_{P\in\mathcal G}A_P,\qquad
 S_{\mathcal A}=\left(\sum_{P\in\mathcal G}\|A_P\|_{\mathcal H}^2\right)^{1/2}.
$$
Suppose that, for every dyadic cube $R\subseteq Q$ and every subfamily
$\mathcal E\subseteq\{P\in\mathcal G:P\subseteq R\}$,
\begin{equation}\label{eqn:local-CWW-hilbert}
 \frac1{\mu(R)}\int_R\exp\left(
 c_0\frac{\|\sum_{P\in\mathcal E}A_P(x)\|_{\mathcal H}^2}
 {\|S_{\mathcal E}\|_{L^\infty(R)}^2}\right)\,{\rm d}\mu(x)\leq C_0,
\end{equation}
where $S_{\mathcal E}=(\sum_{P\in\mathcal E}\|A_P\|_{\mathcal H}^2)^{1/2}$.
The quotient in \eqref{eqn:local-CWW-hilbert} is taken to be zero when
$\|S_{\mathcal E}\|_{L^\infty(R)}=0$.  Then, for $p\geq2$,
\begin{equation}\label{eqn:local-CWW-to-Lp}
       \|G\|_{L^p(Q;\mathcal H)}
       \leq C\sqrt p\,\|S_{\mathcal A}\|_{L^p(Q)},
\end{equation}
where $c_0$ and $C_0$ are independent of $R$, $\mathcal E$, and $\dim\mathcal H$.  The constant $C$ depends only on
$c_0$, $C_0$, and the underlying dyadic structure.
\end{lemma}

\begin{proof}
For $P\in\mathcal G$, denote
$$
 m_P=\frac1{\mu(P)}\int_P S_{\mathcal A}\,{\rm d}\mu.
$$
The quantity $m_P$ records the size of the full square function on the cube indexing $A_P$.  We classify the cubes
according to these averages, rather than by the minimum of $S_{\mathcal A}$ on a cube.  This distinction is needed
because $S_{\mathcal A}$ need not be constant on $P$.
For $j\in\mathbb Z$, let
$$
 \mathcal G_j=\{P\in\mathcal G:2^j<m_P\leq2^{j+1}\},
 \qquad G_j=\sum_{P\in\mathcal G_j}A_P.
$$
If $m_P=0$, then $S_{\mathcal A}=0$ almost everywhere on $P$, and hence $A_P=0$ almost everywhere.  Since the family
$\mathcal G$ is finite, it follows that
$$
 G=\sum_{j\in\mathbb Z}G_j
$$
almost everywhere, with only finitely many nonzero terms.

We first obtain the uniform square-function bound needed to apply the local exponential-square estimate to each layer.
We claim that
\begin{equation}\label{eqn:average-layer-square}
 \|S_{\mathcal G_j}\|_{L^\infty(Q)}
 \leq 2^{j+1}(1+\theta^{-2})^{1/2},
\end{equation}
where $\theta$ is as in \eqref{eqn:double-theta}.  Fix a point $x$, outside the dyadic boundaries, for which
$S_{\mathcal G_j}(x)>0$.  The cubes in $\mathcal G_j$ which contain $x$ form a finite chain.  Let $P_x$ be the smallest
one and write
$$
 B_x(y)^2=\sum_{\substack{P\in\mathcal G_j\\P\supsetneqq P_x}}\|A_P(y)\|_{\mathcal H}^2,
 \qquad y\in P_x.
$$
For every strict ancestor $P$ of $P_x$, the function $A_P$ is constant on the child of $P$ containing $P_x$.
Consequently, $B_x$ is constant on $P_x$.  Since $B_x\leq S_{\mathcal A}$ pointwise on $P_x$, averaging over $P_x$
gives
$$
 B_x(x)\leq m_{P_x}\leq2^{j+1}.
$$
If $P_x'$ is the child of $P_x$ containing $x$, then $A_{P_x}$ is constant on $P_x'$ and
$$
 m_{P_x}\geq\frac1{\mu(P_x)}\int_{P_x'}S_{\mathcal A}\,{\rm d}\mu
 \geq\frac{\mu(P_x')}{\mu(P_x)}\|A_{P_x}(x)\|_{\mathcal H}
 \geq\theta\|A_{P_x}(x)\|_{\mathcal H}.
$$
Every cube in $\mathcal G_j$ containing $x$ is either $P_x$ or a strict ancestor of $P_x$.  It follows that
$$
 S_{\mathcal G_j}(x)^2=B_x(x)^2+\|A_{P_x}(x)\|_{\mathcal H}^2
 \leq2^{2j+2}(1+\theta^{-2}),
$$
which proves \eqref{eqn:average-layer-square}.

Notice that this argument does not require $S_{\mathcal A}$ to be constant on $P_x$.  The average defining $m_{P_x}$,
the constancy of the proper-ancestor terms on $P_x$, and the lower measure bound for the child containing $x$ are
exactly the three ingredients used in \eqref{eqn:average-layer-square}.

Let $K_\theta=2(1+\theta^{-2})^{1/2}$, and let $\{P_{j,i}\}_i$ be the maximal cubes in $\mathcal G_j$.  These cubes are
pairwise disjoint.  Denote by $\mathcal G_{j,i}$ the cubes in $\mathcal G_j$ contained in $P_{j,i}$, and set
$$
 G_{j,i}=\sum_{P\in\mathcal G_{j,i}}A_P,\qquad E_j=\bigcup_iP_{j,i}.
$$
The square function of the restricted family $\mathcal G_{j,i}$ is bounded by that of $\mathcal G_j$.  Thus
\eqref{eqn:local-CWW-hilbert} and \eqref{eqn:average-layer-square} give
$$
 \frac1{\mu(P_{j,i})}\int_{P_{j,i}}
 \exp\left(\frac{c_0}{K_\theta^2}
 \frac{\|G_{j,i}(x)\|_{\mathcal H}^2}{2^{2j}}\right)\,{\rm d}\mu(x)
 \leq C_0.
$$
We recall the Orlicz estimate in the form needed below.  Let
$$
 \Phi(t)=\exp(t^2)-1,\qquad \Psi(t)=t\big(\log(e+t)\big)^{1/2},\qquad t\geq0,
$$
and, for a cube $P$, denote the corresponding normalized Luxemburg norms by
$$
 \|h\|_{\Phi,P}=\inf\left\{\lambda>0:\frac1{\mu(P)}
 \int_P\Phi\left(\frac{|h|}{\lambda}\right)\,{\rm d}\mu\leq1\right\}
$$
and $\|h\|_{\Psi,P}$, respectively.  A Young function complementary to $\Phi$ is controlled by $\Psi$ and is
equivalent to $\Psi$ for large arguments.  Hence the generalized H\"older inequality gives
$$
 \frac1{\mu(P)}\int_P|h(x)k(x)|\,{\rm d}\mu(x)
 \leq C\|h\|_{\Phi,P}\|k\|_{\Psi,P}.
$$
The preceding exponential-square estimate implies
$$
 \big\|\|G_{j,i}\|_{\mathcal H}\big\|_{\Phi,P_{j,i}}\leq C2^j.
$$
Moreover, the $L(\log L)^{1/2}$ norm of a nonnegative function $v$ satisfies
$$
 \|v\|_{\Psi,P}\leq\frac C{\mu(P)}
 \left(v(P)\int_Pv\log\left(e+\frac{v}{v_P}\right)\,{\rm d}\mu\right)^{1/2}.
$$
These are the $\exp(L^2)$--$L(\log L)^{1/2}$ estimates used in
\cite[Theorem~11.2 and Corollary~11.1]{Wi08}.  Applying them on $P_{j,i}$, for an $A_1$ weight $v$ with
$v(P_{j,i})>0$, gives
\begin{align*}
 \int_{P_{j,i}}\|G_{j,i}\|_{\mathcal H}v\,{\rm d}\mu
 &\leq C2^j v(P_{j,i})^{1/2}
 \left(\int_{P_{j,i}}v
 \log\left(e+\frac{v}{v_{P_{j,i}}}\right)\,{\rm d}\mu\right)^{1/2}
 \leq C2^j[v]_{A_1}^{1/2}v(P_{j,i}).
\end{align*}
Here the second inequality follows, as in \cite[(2.15)]{Wi08}, from
$$
 \int_Pv\log\left(e+\frac{v}{v_P}\right)\,{\rm d}\mu
 \leq C\int_P\mathcal M(v\mathsf1_P)\,{\rm d}\mu
 \leq C[v]_{A_1}v(P).
$$
The cubes for which $v(P_{j,i})=0$ make no contribution and may be omitted.  Since the maximal cubes are disjoint and
$G_j=G_{j,i}$ on $P_{j,i}$, we obtain
$$
 \int_Q\|G_j\|_{\mathcal H}v\,{\rm d}\mu
 \leq C2^j[v]_{A_1}^{1/2}v(E_j).
$$

Let $\mathcal M_{\mathcal D}$ denote the dyadic maximal operator.  Since
$m_{P_{j,i}}>2^j$,
$$
 E_j\subseteq
 \{x\in Q:\mathcal M_{\mathcal D}(S_{\mathcal A}\mathsf1_Q)(x)>2^j\}.
$$
Using the triangle inequality for $G=\sum_jG_j$ and then summing the last weighted estimate over $j$, we have
\begin{align}\label{eqn:weighted-all-layers}
 \int_Q\|G\|_{\mathcal H}v\,{\rm d}\mu
 &\leq C[v]_{A_1}^{1/2}\sum_{j\in\mathbb Z}2^jv(E_j)
 \leq C[v]_{A_1}^{1/2}
 \int_Q\mathcal M_{\mathcal D}(S_{\mathcal A}\mathsf1_Q)v\,{\rm d}\mu.
\end{align}
Here we used the pointwise estimate
$$
 \sum_{j\in\mathbb Z}2^j
 \mathsf1_{\{\mathcal M_{\mathcal D}(S_{\mathcal A}\mathsf1_Q)>2^j\}}
 \leq2\mathcal M_{\mathcal D}(S_{\mathcal A}\mathsf1_Q).
$$

It remains to choose the weight.  Let $2\leq p<\infty$, and take $\varphi\geq0$ with
$\|\varphi\|_{L^{p'}(Q)}=1$.  Extend $\varphi$ by zero outside $Q$ and define the Rubio de Francia majorant
$$
 v=\sum_{k=0}^\infty
 \frac{\mathcal M^k\varphi}{(2\|\mathcal M\|_{L^{p'}(X)})^k}.
$$
The term with $k=0$ shows that $v\geq\varphi$.  The other two properties follow directly from
the definition.  Indeed,
\begin{align*}
 \|v\|_{L^{p'}(X)}
 &\leq\sum_{k=0}^\infty
 \frac{\|\mathcal M^k\varphi\|_{L^{p'}(X)}}
 {(2\|\mathcal M\|_{L^{p'}(X)})^k}
 \leq\sum_{k=0}^\infty2^{-k}=2,
\end{align*}
and, by the sublinearity of $\mathcal M$,
\begin{align*}
 \mathcal Mv
 \leq\sum_{k=0}^\infty
 \frac{\mathcal M^{k+1}\varphi}
 {(2\|\mathcal M\|_{L^{p'}(X)})^k}
 =2\|\mathcal M\|_{L^{p'}(X)}
 \sum_{k=1}^\infty\frac{\mathcal M^k\varphi}
 {(2\|\mathcal M\|_{L^{p'}(X)})^k}
 \leq2\|\mathcal M\|_{L^{p'}(X)}v.
\end{align*}
In particular, $v\in A_1$ and the maximal operator estimate in \cite{Buck} yields
$$
 [v]_{A_1}\leq2\|\mathcal M\|_{L^{p'}(X)}\leq Cp.
$$
The dyadic maximal operator is bounded on $L^p(X)$ uniformly for $p\geq2$.  Hence \eqref{eqn:weighted-all-layers},
H\"older's inequality, and the preceding bounds give
\begin{align*}
 \int_Q\|G\|_{\mathcal H}\varphi\,{\rm d}\mu
 &\leq\int_Q\|G\|_{\mathcal H}v\,{\rm d}\mu
 \leq C[v]_{A_1}^{1/2}
 \|\mathcal M_{\mathcal D}(S_{\mathcal A}\mathsf1_Q)\|_{L^p(X)}
 \|v\|_{L^{p'}(X)}
 \leq C\sqrt p\,\|S_{\mathcal A}\|_{L^p(Q)}.
\end{align*}
Duality proves \eqref{eqn:local-CWW-to-Lp}.
\end{proof}

\bigskip

{\bf Proof of Proposition \ref{coro:Lp-CH}:}
Let $\mathcal G=\bigcup_{0\leq n\leq N}Ch^{(n)}(Q)$.  We verify the local hypothesis of
Lemma~\ref{lem:local-CWW-to-Lp}.  Fix a dyadic cube $R\subseteq Q$ and a subfamily
$$
 \mathcal E\subseteq\{P\in\mathcal G:P\subseteq R\}.
$$
The assertion is immediate if $\mathcal E$ is empty.  Otherwise, choose $L$ so that every cube in $\mathcal E$ belongs
to $\bigcup_{0\leq k\leq L}Ch^{(k)}(R)$.  For a cube $P$ in these generations, define
$$
 b_P=\begin{cases}
 \mathfrak a_P,&P\in\mathcal E,\\
 0,&P\notin\mathcal E.
 \end{cases}
$$
The functions $b_P$ satisfy \textbf{(P1)} and \textbf{(P2)}.  To verify \textbf{(P3)}, let $W\subseteq R$ be
dyadic, let $\mathcal F$ be any subfamily of the cubes indexing the $b_P$'s which are contained in $W$, and put
$$
 \mathcal F_0=\{P\in\mathcal F:b_P\neq0\}.
$$
Then
$$
 \sum_{P\in\mathcal F}b_P=\sum_{P\in\mathcal F_0}\mathfrak a_P,
 \qquad
 \sum_{P\in\mathcal F}|b_P|^2=\sum_{P\in\mathcal F_0}|\mathfrak a_P|^2.
$$
Consequently, \textbf{(P3)} for the original atoms, applied to $\mathcal F_0$, gives
$$
 \left(\frac1{\mu(W)}\int_W\left|\sum_{P\in\mathcal F}b_P\right|^2\,{\rm d}\mu\right)^{1/2}
 \leq C_1\left\|\left(\sum_{P\in\mathcal F}|b_P|^2\right)^{1/2}\right\|_{L^\infty(W)}.
$$
Thus the $b_P$'s satisfy \textbf{(P3)} with the same constant $C_1$.  Set
$$
 G_{\mathcal E}=\sum_{P\in\mathcal E}\mathfrak a_P,
 \qquad
 S_{\mathcal E}=\left(\sum_{P\in\mathcal E}|\mathfrak a_P|^2\right)^{1/2}.
$$
If $\|S_{\mathcal E}\|_{L^\infty(R)}=0$, then \textbf{(P3)} gives $G_{\mathcal E}=0$ almost everywhere on $R$,
and the required assertion is immediate.  Otherwise, applying \eqref{eqn:newCWW-one-parameter} with top cube $R$
and terminal generation $L$, we obtain
$$
 \frac1{\mu(R)}\int_R\exp\left(
 c_3\frac{|G_{\mathcal E}(x)|^2}
 {\|S_{\mathcal E}\|_{L^\infty(R)}^2}
 \right)\,{\rm d}\mu(x)\leq c_4.
$$
The constants are independent of $R$, $\mathcal E$, and $L$.  Hence \eqref{eqn:local-CWW-hilbert} holds for
$\mathcal H=\mathbb C$ and $A_P=\mathfrak a_P$.  Since the corresponding functions $G$ and $S_{\mathcal A}$ are
$F_N$ and $S(F_N)$, respectively, Lemma~\ref{lem:local-CWW-to-Lp} gives \eqref{eqn:Lp-CH}.
\hfill $\Box$

\smallskip

At the end of this section, we show \eqref{eqn:newCH-dyadic} by using Proposition \ref{coro:Lp-CH}.

\medskip

{\bf Proof of  \eqref{eqn:newCH-dyadic}:}

It follows easily via the Taylor series expansion of 
the exponential in $\int_{Q} \exp\big(  \left|\alpha  F_N (x) \right|  \big)  \,{\rm d}\mu(x) $, Stirling's approximation, and the precise dependence 
on the $L^k$ norm together with the important observation that the constant $C$ in that $L^p$ estimate \eqref{eqn:Lp-CH}  for $p\geq2$ is independent of $p$. 
 \hfill $\Box$

\medskip
This finishes the proofs of Theorem~\ref{thm:newCH-dyadic} and Proposition \ref{coro:Lp-CH}, and establishes their equivalence.

\smallskip



\section{Bi-parameter Hoeffding inequalities:  proof of Theorem \ref{thm:double-newCH-dyadic}}\label{sec:4}
\setcounter{equation}{0}

For  $i=1,2$, let  $(X_i,\rho_i,\mu_i)$ be a quasi-metric space with a doubling measure $\mu_i$ satisfying  $\mu_i(B)<\infty$ for all balls $B\subseteq X_i$,
and let $\mathcal D^{[i]}$ be a system of dyadic cubes on $X_i$.

We prove Theorem \ref{thm:double-newCH-dyadic} in this section by employing a reduction machinery based on the one-parameter version in Section \ref{sec:3}.
 
For any given $Q_1\times Q_2\in \mathcal D^{[1]}\times \mathcal D^{[2]}$, let $\left\{\mathfrak a_{P_1\times P_2}:\, P_1\times P_2\in \mathcal D^{[1]}\times \mathcal D^{[2]}, \,  P_1\subseteq Q_1,\, P_2\subseteq Q_2\right\}$ be ``atomic functions'' satisfying \textsf{(P1)}, \textsf{(P2)} and \textsf{(P3)}  in Theorem \ref{thm:double-newCH-dyadic}. Recall 
$$
    F_{N,N}(  x_1,  x_2)=\sum_{\substack{P_1\in  Ch^{(n)}(Q_1)\\ 0\leq n\leq N}} \sum_{\substack{P_2\in  Ch^{(n)}(Q_2)\\ 0\leq n\leq N}} \mathfrak{a}_{P_1\times P_2} ( x_1,  x_2)
$$
and 
$$
   S(F_{N,N})(  x_1, x_2) =\Bigg(\sum_{\substack{P_1\in  Ch^{(n)}(Q_1)\\ 0\leq n\leq N}} \sum_{\substack{P_2\in  Ch^{(n)}(Q_2)\\ 0\leq n\leq N}}  { \mathfrak{a}^2_{P_1\times P_2}    ( x_1,x_2)  }\Bigg)^{1/2}.
$$

For any fixed $x_2\in Q_2$, we have 
$$
    F_{N,N}(x_1, x_2):=\sum_{\substack{P_1\in  Ch^{(n)}(Q_1)\\ 0\leq n\leq N}}    \mathfrak a_{P_1}^{(N)}(x_1, x_2),\quad \text{where}  
    \quad 
    \mathfrak a_{P_1}  ^{(N)}=\sum_{\substack{P_2\in  Ch^{(n)}(Q_2)\\ 0\leq n\leq N}}    \mathfrak a_{P_1\times P_2}(x_1, x_2).
$$
Since  $\left \{\mathfrak a_{P_1\times P_2}\right \}$ satisfies \textsf{(P1)}, \textsf{(P2)} and \eqref{eqn:double-P3-revise-I} in
\textsf{(P3)},   for any fixed $x_2\in Q_2$,
$\left\{\mathfrak a_{P_1}^{(N)}\right\}$ as a sequence of functions in $x_1$, satisfies the properties \textbf{(P1)}, \textbf{(P2)} and \textbf{(P3)} listed in Theorem \ref{thm:newCH-dyadic}. Hence, for any fixed $x_2\in Q_2$ and $\alpha>0$, Theorem \ref{thm:newCH-dyadic} implies that 
\begin{equation}\label{eqn:fix-newCH-dyadic}
       \int_{Q_1} \exp\left(\alpha \left|F_{N, N}(x_1, x_2) \right| \right) {\rm d}\mu_1(x_1)\leq c_1   \int_{Q_1} \exp\left( c_2 \alpha^2   \big(S_1(F_{N,N})\big)^2(x_1, x_2) \right){\rm d}\mu_1(x_1),	
\end{equation}
where 
\begin{align}\label{eqn:S1}
	 S_1(F_{N,N})(x_1,x_2):=
	 \Bigg(\sum_{\substack{P_1\in  Ch^{(n)}(Q_1)\\ 0\leq n\leq N}} {\Big|\sum_{\substack{P_2\in  Ch^{(n)}(Q_2)\\ 0\leq n\leq N}}  \mathfrak{a}_{P_1\times P_2}(  x_1, x_2) \Big | ^2  }\,\Bigg)^{1/2}.
\end{align}
 
Furthermore, we have the sharp $L^p$ lower bounds for  $S_1 (F_{N,N})$ whenever $p\geq 2$, as a direct consequence of Proposition  \ref{coro:Lp-CH}. This is also the first crucial step of showing Theorem \ref{thm:double-newCH-dyadic}.

\smallskip

\begin{corollary}\label{coro:part-I-Lp-CH}
With the same assumptions and notation as in Theorem \ref{thm:double-newCH-dyadic}, there exists a constant $C>0$ such that for any fixed $x_2\in Q_2$ and $p\geq 2$,
\begin{equation}\label{eqn:part-I-Lp-CH}
	   \left\|F_{N,N}(\cdot, x_2)\right\|_{L^p(Q_1,{\rm d}\mu_1)}\leq C\sqrt{p} \left\|S_1(F_{N,N})(\cdot, x_2)\right\|_{L^p(Q_1,{\rm d}\mu_1)}.
\end{equation}
\end{corollary}

\smallskip

Based on \eqref{eqn:part-I-Lp-CH}, 
one suitable strategy of proving  Theorem \ref{thm:double-newCH-dyadic} from the relationship  (see the argument below in this section) between Hoeffding's inequalities and sharp $L^p$ lower bound for the square function $S(F_{N,N})$  is to first  show the following estimate.

\begin{proposition}\label{prop:double-Lp-CH}
Let all assumptions and notation be the same as in Theorem \ref{thm:double-newCH-dyadic}. There exists a constant $C>0$ such that for any fixed $x_1\in Q_1$ and $p\geq 2$,
\begin{equation}\label{eqn:part-II-Lp-CH}
	   \left\|S_1(F_{N,N})(x_1, \cdot)\right\|_{L^p(Q_2,{\rm d}\mu_2)}\leq C\sqrt{p} \left\|S(F_{N,N})(x_1,\cdot)\right\|_{L^p(Q_2,{\rm d}\mu_2)}.
\end{equation}
\end{proposition}

As a consequence of  Corollary \ref{coro:part-I-Lp-CH} and Proposition \ref{prop:double-Lp-CH}, we see that 
Proposition \ref{prop:double-Lp-CH-intro} holds.

Next we will show that Proposition \ref{prop:double-Lp-CH} is closely connected with the following variation  of Hoeffding's inequality. 

\begin{proposition}\label{prop:one-shift-bi-I}
	Under the same assumptions and notation as in Theorem \ref{thm:double-newCH-dyadic}, there exist constants $c>1$ and $C>0$ such that 
	for any $Q_1$, $Q_2$ and $N$, and for any fixed $x_1\in Q_1$ and $\alpha>0$,
\begin{equation}\label{eqn:one-shift-bi-I}
   	\int_{Q_2} \exp\left(  \alpha S_1(F_{N,N})(x_1,x_2)  \right){\rm d}\mu_2(x_2)\leq C\int_{Q_2} \exp\left(  c \alpha^2 S^2(F_{N,N})(x_1,x_2)\right){\rm d}\mu_2(x_2).
\end{equation}
\end{proposition}
 \smallskip

Remarkably, as mentioned in the first section, if we define $S_1(F_{N,N})$ via the version \eqref{eqn: 1.9}, then the (almost) orthogonality of atomic functions is destroyed by the $L^\infty$ norm there, which prevents us from obtaining  Proposition \ref{prop:one-shift-bi-I}. This not only shows the necessity of our definition of $S_1(F_{N,N})$ given in \eqref{eqn:S1} without the $L^\infty$ norm, but also, in turn,  highlights the fundamental  role of our Theorem \ref{thm:newCH-dyadic} with the new square function $S(F_N)$, which implies  Corollary \ref{coro:part-I-Lp-CH}.

When taking atomic functions as double dyadic martingale differences $d_{n,m}$ on $[0,1]^2$, Proposition \ref{prop:one-shift-bi-I} can be obtained by applying a reduction argument; see  \cite[Lemma 2.2]{Pi86} by the second author. 

{However, one may note that the proof of \cite[Lemma 2.2]{Pi86} depends not only on the cancellation of $d_{n,m}$, but also on the special property that $d_{n,m}$ as a one-variable function of $x_2$ (with $x_1$ fixed),  merely changes sign on either half of the dyadic intervals of length $2^{-m+1}$. This special property  is invalid for the setting of higher dimensional space $\mathbb R^d$ or general spaces of homogeneous type.} Indeed,  let $\{\Delta_{P_1\times P_2}(f)\}$ be double martingale differences of some function $f$  on a product space of homogeneous type,  the argument in \cite[Lemma 2.2]{Pi86}  
 can deduce that 
\begin{align}\label{eqn:one-shift-bi-I-weak}
		\int_{Q_2} \exp\left(  \sqrt{1+S^2_1(\mathcal F_{N,N})     (x_1,x_2)} -c \widetilde{S}^2(\mathcal F_{N,N})(x_1,x_2)\right){\rm d}\mu_2(x_2)\leq C\mu_2(Q_2),
\end{align}
where 
$$
    \mathcal F_{N,N}(x_1,x_2):=\sum_{\substack{P_1\in  Ch^{(n)}(Q_1)\\ 0\leq n\leq N}} \sum_{\substack{P_2\in  Ch^{(n)}(Q_2)\\ 0\leq n\leq N}} \Delta_{P_1\times P_2}(f)(x_1,x_2)
$$
and 
$$
     \widetilde{S}^2(\mathcal F_{N,N})(x_1,x_2):=\sum_{\substack{P_1\in  Ch^{(n)}(Q_1)\\ 0\leq n\leq N}} 
     \sum_{\substack{P_2\in  Ch^{(n)}(Q_2)\\ 0\leq n\leq N}} \left \|\Delta_{P_1\times P_2}(f)\right\|_{L^\infty(P_1\times P_2)}^2\mathsf 1_{P_1\times P_2}(x_1,x_2).
$$
It's clear that $\widetilde{S}(F_{N,N})\geq S(F_{N,N})$ pointwise, and the crucial ingredient for \eqref{eqn:one-shift-bi-I-weak} is  Hoeffding's lemma based on the cancellation of $\Delta_{P_1\times P_2}(f)$, this creates an extreme limitation for many potential applications.

{Additionally, even if we assume the cancellation property,  when studying the 3-parameter Hoeffding type inequality,  the analog of Proposition \ref{prop:one-shift-bi-I} will play the role of Theorem \ref{thm:newCH-dyadic} for bi-parameter theory, and consequently \eqref{eqn:one-shift-bi-I-weak},  where the $L^\infty$ structure occurs, will not be suitable. }

The argument in the biparameter setting considers, as in the earlier one parameter context,  two cases: the ``small scale" and "large scale" regimes. Here we highlight that, using our trick for estimating $\mathcal I_R(C)$ in \eqref{eqn:IR},  under the large scale condition, an improvement on \eqref{eqn:one-shift-bi-I-weak} can be obtained in which the square function $\widetilde S$ is replaced by the square function $S$ defined via double martingale differences, or defined via our atomic functions $\mathfrak a_{P_1\times P_2}$. The large scale regime is defined by $\|\mathfrak a_{P_1\times P_2}\|_{L^\infty}\geq1$, and the argument for this case is given in Step I below.
As in Theorem \ref{thm:newCH-dyadic}, the small-scale case is more complicated,
and we address it in the subsequent Step II  by modifying the classification algorithm and the John--Nirenberg type argument of Theorem \ref{thm:newCH-dyadic}. 

Therefore, our Proposition \ref{prop:one-shift-bi-I} for the new square function without $L^\infty$-norm  is not only a suitable substitution of \eqref{eqn:one-shift-bi-I-weak} from $\widetilde S$ to $S$ for the double martingale differences, but also sharp enough for various applications on martingales or in more general settings where one only assumes almost orthogonality. {Both its expression and argument are also flexible and necessary for advancing the study of multi-parameter Hoeffding type inequalities.} These, in turn,  lead to potential applications to singular integrals in multi-parameter harmonic analysis and other related questions.

\bigskip

Similar to the organization of Section \ref{sec:3}, we begin by claiming a rough version of \eqref{eqn:one-shift-bi-I}: for any fixed $x_1\in Q_1$ and $\alpha>0$,
\begin{equation}\label{eqn:aux-S1-S-weak-0}
    \int_{Q_2}\exp\bigg(\alpha S_1(F_{N,N}) (x_1,x_2)\bigg)\, {\rm d}\mu_2(x_2)\leq C  \exp\left(c\alpha^2  \|S(F_{N,N})(x_1,\cdot)\|_{L^\infty(Q_2)}^2\right) \mu_2(Q_2).
\end{equation}

Once \eqref{eqn:aux-S1-S-weak-0} is established, Theorem \ref{thm:double-newCH-dyadic} follows readily.

\bigskip

{\bf Proof of  Theorem \ref{thm:double-newCH-dyadic}: } 

As a consequence of \eqref{eqn:aux-S1-S-weak-0}, similar to the proof of  \eqref{eqn:newCWW-one-parameter}, we can obtain the corresponding exponential square integrability: there exist positive constants $c_3'$ and $c_4'$ such that for any fixed $x_1\in Q_1$,
\begin{equation}\label{eqn:newCWW-S1-S-0}
       \int_{Q_2} \exp\left(    c_3'    \frac { |S_1(F_{N,N})(x_1,x_2) |^2 }   {\|S(F_{N,N})(x_1,\cdot )\|_ {L^\infty(Q_2) }^2  }
       \right)  {\rm d}\mu_2(x_2)\leq c_4' \, \mu_2(Q_2).
\end{equation}
Thus
\begin{equation}\label{eqn:newCWW-S1-S}
       \int_{Q_2} \exp\left(    c_3'    \frac { |S_1(F_{N,N})(x_1,x_2) |^2 }   {\|S(F_{N,N})\|_ {L^\infty(Q_1\times Q_2) }^2  }
       \right)  {\rm d}\mu_2(x_2)\leq c_4' \, \mu_2(Q_2)
\end{equation}
for any $x_1\in Q_1$.

Additionally, it follows from \eqref{eqn:newCH-dyadic} that for any fixed $x_2\in Q_2$,
\begin{equation}\label{eqn:newCH-dyadic-coro1}
       \int_{Q_1} \exp\left(\alpha \left|F_{N,N}(x_1,x_2)\right|  
       \right)  {\rm d}\mu_1(x_1)\leq c_1' \int_{Q_1} \exp\left( c_2'\alpha^2   S_1^2  (F_{N,N})  (x_1,x_2)
       \right)  {\rm d}\mu_1(x_1)	,
\end{equation}
where $c_1'$ and $c_2'$ depends on the parameter $C_2$ in the almost orthogonal property \textsf{(P3)}. 

Therefore, taking $\displaystyle
\alpha=\sqrt{\frac{c_3'}{c_2'}}\,
\|S(F_{N,N})\|_{L^\infty(Q_1\times Q_2)}^{-1}$, one may combine
\eqref{eqn:newCH-dyadic-coro1} and \eqref{eqn:newCWW-S1-S} to see
\begin{align*}
&\int_{Q_1\times Q_2} \exp\Bigg(   \sqrt{\frac{c_3'}  {c_2' }     }    \frac{   | F_{N,N}( x_1, x_2)|   }{\|S(F_{N,N})\|_{L^\infty(Q_1\times Q_2)}}\Bigg)  \,{\rm d}\mu_1(  x_1) \, {\rm d}\mu_2( x_2)\\
\leq & c_1' \int_{Q_1\times Q_2} \exp\Bigg(    c_3'    \frac{   |S_1( F_{N,N})( x_1, x_2)|^2   }{\|S(F_{N,N})\|_{L^\infty(Q_1\times Q_2)}^2   }\Bigg)  \,{\rm d}\mu_1(  x_1) \, {\rm d}\mu_2( x_2)\\
 \leq & c_4' \mu_1(Q_1)\,\mu_2(Q_2).
\end{align*}
That is, \eqref{eqn:double-newCH-dyadic} holds and we complete the proof of Theorem \ref{thm:double-newCH-dyadic}.
 \hfill $\Box$

\bigskip

Another straightforward consequence of  \eqref{eqn:newCWW-S1-S-0} is Proposition \ref{prop:double-Lp-CH}, which in turn, allows us to obtain Proposition \ref{prop:one-shift-bi-I}.  From the above, our outline includes 
\begin{subequations}
\begin{equation}\label{F-a}
   \eqref{eqn:aux-S1-S-weak-0}  \Longrightarrow  \eqref{eqn:newCWW-S1-S-0} \overset{\eqref{eqn:newCH-dyadic}}{ \Longrightarrow  }  \text{Theorem \ref{thm:double-newCH-dyadic}}   
\end{equation}
and 
\begin{equation}\label{F-b}
    \eqref{eqn:newCWW-S1-S-0}    \Longrightarrow \text{Proposition \ref{prop:double-Lp-CH}}  \Longrightarrow  \text{Proposition \ref{prop:one-shift-bi-I}}   \   (\Longrightarrow    \eqref{eqn:aux-S1-S-weak-0} ). 
\end{equation} 
\end{subequations}

\smallskip

To complete the proofs of these estimates and their implications, we begin by showing  \eqref{eqn:aux-S1-S-weak-0}.

\bigskip

{\bf Proof of  \eqref{eqn:aux-S1-S-weak-0}: } 



For any fixed $x_1\in Q_1$ and $0\leq k\leq N$, set
$$
   A_k(x_1,x_2) :=S_1^2(F_{N,k})(x_1,x_2)=\sum_{\substack{P_1\in  Ch^{(n)}(Q_1)\\ 0\leq n\leq N}} \Big|\sum_{\substack{P_2\in  Ch^{(n)}(Q_2)\\ 0\leq n\leq k}}  \mathfrak{a}_{P_1\times P_2} (  x_1, x_2)\Big|^2
 $$
and 
$$
    D_k(x_1,x_2):=\sum_{\substack{P_1\in  Ch^{(n)}(Q_1)\\ 0\leq n\leq N}} \Big|\sum_{ P_2\in  Ch^{(k)}(Q_2) }  \mathfrak{a}_{P_1\times P_2}(  x_1, x_2) \Big |^2 .
$$

Recall the fact that for each $0\leq n\leq k$ and $(x_1,x_2)\in Q_1\times Q_2$, there exists a unique $P_2\in Ch^{(n)}(Q_2)$ occurring in $A_k(x_1,x_2)$. Similarly, there exists a unique $P_2\in Ch^{(k)}(Q_2)$ occurring in $D_k(x_1,x_2)$. Hence
$$
    D_k(x_1,x_2)=\sum_{\substack{P_1\in  Ch^{(n)}(Q_1)\\ 0\leq n\leq N}} \sum_{ P_2\in  Ch^{(k)}(Q_2) }   \mathfrak{a}^2_{P_1\times P_2}(  x_1, x_2)  =\sum_{P_2\in Ch^{(k)}(Q_2)}  \left( \mathcal A_{P_2}\right)^2 (x_1,x_2),
$$
where 
$$ \mathcal A_{P_2}(x_1,x_2):=\bigg( \sum_{\substack{P_1\in  Ch^{(n)}(Q_1)\\ 0\leq n\leq N}}\mathfrak{a}^2_{P_1\times P_2}(  x_1, x_2)\bigg)^{1\over2} .$$

Hence
$$
    S^2(F_{N,k})(x_1,x_2)=\sum_{n=0}^k D_n(x_1,x_2) , \quad 0\leq k\leq N.
$$

Denote also
$$
     r_k(x_1,x_2):=2\sum_{\substack{P_1\in  Ch^{(n)}(Q_1)\\ 0\leq n\leq N}}  \Bigg[\sum_{\substack{P_2\in  Ch^{(n)}(Q_2)\\ 0\leq n\leq k-1}}  \mathfrak{a}_{P_1\times P_2} (  x_1, x_2)\cdot\sum_{ P_2\in  Ch^{(k)}(Q_2) }  \mathfrak{a}_{P_1\times P_2}(  x_1, x_2)\Bigg],\qquad k\geq 1.
$$
Then we see that 
\begin{equation}\label{eqn:Ak}
	  A_k(x_1,x_2)=A_{k-1}(x_1,x_2)+ D_k(x_1,x_2)+ r_k(x_1,x_2).
\end{equation}
Moreover,  by Cauchy--Schwarz inequality,
$$
    r_k^2(x_1,x_2)\leq 4 A_{k-1}(x_1,x_2) D_k(x_1,x_2).
$$

Recall that we split atomic functions in $F_N$ into two parts in terms of their $L^\infty$ norms. Similarly, let
$$
   \mathcal G_2^{(I)}:=\bigg\{P_2\in\bigcup_{0\leq n\leq N}Ch^{(n)}(Q_2):\   \left\|\mathcal A_{P_2}(x_1,\cdot)\right\|_{L^\infty(P_2)}< 1\bigg\}
$$
and
$$
   \mathcal G_2^{(II)}:=\bigg\{P_2\in\bigcup_{0\leq n\leq N}Ch^{(n)}(Q_2):\   \left\|\mathcal A_{P_2}(x_1,\cdot)\right\|_{L^\infty(P_2)}\geq 1 \bigg\}.
$$

Denote 
$$
      F_{N,N}^{(I)}:= \sum_{\substack{P_1\in  Ch^{(n)}(Q_1)\\ 0\leq n\leq N}} \sum_{P_2\in \mathcal G_2^{(I)}} \mathfrak{a}_{P_1\times P_2} \quad \text{and}\quad      
       F_{N,N}^{(II)}:= \sum_{\substack{P_1\in  Ch^{(n)}(Q_1)\\ 0\leq n\leq N}} \sum_{P_2\in \mathcal G_2^{(II)}} \mathfrak{a}_{P_1\times P_2} .
$$
Then by a similar argument \eqref{eqn:weak-newCH-dyadic} in {\bf  Step III} in the proof of \eqref{eqn:newCH-dyadic-2}, it suffices to show that \eqref{eqn:aux-S1-S-weak-0} holds for some  constant $\alpha\in (0,1]$ with  $F_{N,N}$ replaced by $F_{N,N}^{(I)}$ and by $F_{N,N}^{(II)}$, respectively. 

\smallskip

\noindent {\it Step I.} We begin by considering \eqref{eqn:aux-S1-S-weak-0}  with $\alpha=1$ for $F_{N,N}^{(II)}$.

Without loss of generality, one may assume that $\mathcal G_2^{(II)}=\bigcup_{0\leq n\leq N}Ch^{(n)}(Q_2) $, that is, $F_{N,N}^{(II)}=F_{N,N}$. 
The argument here involves an iteration on the generations of dyadic cubes in the second summation that occurs in $S_1(F_{N,N})$, which is motivated by \cite[Lemma 2.2]{Pi86}.  However, unlike the original work \cite{Pi86} which relies on the vanishing property of martingale differences, the feasibility of the iteration below is similar to the argument of \eqref{eqn:newCH-dyadic-2-II} about $F_N^{(II)}$ in  {\bf  Step I} in the proof of \eqref{eqn:newCH-dyadic-2}. 

To prove \eqref{eqn:aux-S1-S-weak-0} with $\alpha=1$, we turn to show a slightly stronger estimate that  for any $1\leq k\leq N$,
\begin{align}\label{eqn:one-shift-bi-iteration}
	  &\int_{Q_2} \exp\left( \sqrt{1+S^2_1(F_{N,k})(x_1,x_2)} -c   S^2(F_{N,k})(x_1,x_2)\right){\rm d}\mu_2(x_2)\nonumber\\
	  \leq &\int_{Q_2} \exp\left(  \sqrt{1+S^2_1(F_{N,k-1})(x_1,x_2)} -c  S^2(F_{N,k-1})(x_1,x_2)\right){\rm d}\mu_2(x_2).
\end{align}
Once \eqref{eqn:one-shift-bi-iteration} is verified,  \eqref{eqn:aux-S1-S-weak-0} with $\alpha=1$ follows from the following estimate:
\begin{align*}
	&\int_{Q_2} \exp\left(  \sqrt{1+S^2_1(F_{N,0})(x_1,x_2)} -c  S^2(F_{N,0})(x_1,x_2)\right){\rm d}\mu_2(x_2)\\
	=& \int_{Q_2} \exp\Bigg(  \sqrt{1+  \sum_{\substack{P_1\in  Ch^{(n)}(Q_1)\\ 0\leq n\leq N}}   \mathfrak a^2_{P_1\times Q_2}      (x_1,x_2) } -c  \sum_{\substack{P_1\in  Ch^{(n)}(Q_1)\\ 0\leq n\leq N}}
	 \mathfrak a^2_{P_1\times Q_2} (x_1,x_2)\Bigg)  \, {\rm d}\mu_2(x_2)\\
	\leq & \int_{Q_2} \exp\Bigg(   1+ \frac{1}{2} \sum_{\substack{P_1\in  Ch^{(n)}(Q_1)\\ 0\leq n\leq N}}  \mathfrak a^2_{P_1\times Q_2}  (x_1,x_2)      -c  \sum_{\substack{P_1\in  Ch^{(n)}(Q_1)\\ 0\leq n\leq N}}
	  \mathfrak a^2_{P_1\times Q_2}  (x_1,x_2) \Bigg)  \,{\rm d}\mu_2(x_2)\\
	\leq&  e \,\mu_2(Q_2),
\end{align*}
where in the first inequality we use the fact  $\sqrt{1+\Box}\leq 1+\Box/2$ for any $\Box>0$, and in the last inequality we use the condition that $c>1$. 

It remains to show \eqref{eqn:one-shift-bi-iteration}.
For arbitrary $R_2\in Ch^{(k)}(Q_2)$, note that 
\begin{itemize}
\item $A_{k-1}=S_1^2(F_{N,k-1})$ is constant on $R_2$.
\smallskip
\item $S(F_{N, k-1})$ and $S(F_{N, k})$ are both constants on $R_2$, and for $x_2\in  R_2$,
\begin{align*}
	S^2(F_{N, k}) (x_1,\cdot)&=S^2(F_{N,k-1}) (x_1,x_2)+\sum_{\substack{P_1\in  Ch^{(n)}(Q_1)\\ 0\leq n\leq N}}    \mathfrak{a}^2_{P_1\times R_2}   (x_1,x_2) \\
	&= S^2(F_{N,k-1})(x_1,x_2)+   D_k(x_1,x_2).
\end{align*}
\end{itemize}
Based on the above arguments,  applying \eqref{eqn:Ak} and using the inequality  $\sqrt{1+\Box}\leq 1+\Box/2$ again, we obtain that 
\begin{align}\label{eqn:aux-iteration-S1-case1}
    &\int_{R_2} \exp\left( \sqrt{1+S^2_1(F_{N,k})(x_1,x_2)} -c   S^2(F_{N,k})(x_1,x_2)\right)  {\rm d}\mu_2(x_2)\nonumber\\
    =&\int_{R_2} \exp\left(\sqrt{1+A_{k-1}(x_1, x_2)} \cdot \sqrt{1+\frac{D_k(x_1,x_2)+r_k(x_1,x_2) }{1+ A_{k-1}(x_1, x_2) }} -c   S^2(F_{N,k})(x_1,x_2)\right)  {\rm d}\mu_2(x_2)\nonumber\\ 
    \leq& \int_{R_2}\exp\bigg(\sqrt{1+A_{k-1}(x_1,c_{R_2})}  +\frac{D_k(x_1,x_2)+r_k(x_1,x_2)}{2\sqrt{1+A_{k-1}(x_1,c_{R_2})}}\bigg. \nonumber\\
    &\qquad \qquad \bigg.  - cS^2(F_{N, k-1})(x_1, c_{R_2})-c D_k(x_1,x_2) \bigg)\,  {\rm d}\mu_2(x_2) \nonumber\\
    \leq & \exp\left(\sqrt{1+A_{k-1}(x_1,c_{R_2})}  - cS^2(F_{N, k-1})(x_1, c_{R_2}) \right) \nonumber\\
    & \quad \times \int_{R_2} \exp\left(  \frac{1}{2\sqrt{1+A_{k-1}(x_1,c_{R_2})} } r_k(x_1,x_2)-\Big (c-\frac{1}{2}\Big)\,  D_k(x_1,x_2)  \right){\rm d}\mu_2(x_2),
 \end{align}
 where $c_{R_2}$ denotes the center of $R_2\in Ch^{(k)}(Q_2)$.

Note that
$$
    \left(\frac{1}{2\sqrt{1+A_{k-1}(x_1,c_{R_2})} } r_k(x_1,x_2)\right)^2 \leq \frac{A_{k-1}(x_1,c_{R_2})  D_k(x_1,x_2)  }   {1+A_{k-1}(x_1,c_{R_2})}\leq D_k(x_1,x_2),
$$
and  $D_k(x_1,\cdot)$ are constant on each child of $R_2$. 
Thus,
\begin{align*}
    &\int_{R_2} \exp\left(  \frac{1}{2\sqrt{1+A_{k-1}(x_1,c_{R_2})} } r_k(x_1,x_2)-\Big (c-\frac{1}{2}\Big)\,  D_k(x_1,x_2)  \right){\rm d}\mu_2(x_2)\\
    \leq &\sum_{P_2\in Ch(R_2): \,|  D_k (x_1,c_{P_2})|<1}\int_{P_2} \exp\left(  \sqrt{D_k(x_1,x_2) }-\Big (c-\frac{1}{2}\Big)\,  D_k(x_1,x_2)  \right) {\rm d}\mu_2(x_2) \nonumber\\
    & \quad +\sum_{P_2\in Ch(R_2): \,| D_k (x_1,c_{P_2})|\geq 1}\int_{P_2} \exp\left(  \sqrt{D_k(x_1,x_2) }-\Big (c-\frac{1}{2}\Big)\,  D_k(x_1,x_2)  \right) {\rm d}\mu_2(x_2)\\
     =:&\widetilde{\mathcal I}_{R_2} \big (c-\frac{1}{2}\big ) .
\end{align*}
 Similar to the estimate for $\mathcal I_R(C)$ at  \textbf{Step I} in the proof of \eqref{eqn:newCH-dyadic-2},  due to  $\|D_k(x_1,\cdot)\|_{L^\infty (R_2)}\geq 1$ assumed in this case, one can take $c=c(\theta)$ sufficiently large, where $\theta$ is the parameter in \eqref{eqn:double-theta},
such that 
$$
    \widetilde{\mathcal I}_{R_2}\big (c-\frac{1}{2}\big ) \leq \mu_2(R_2).
$$
Plugging this into \eqref{eqn:aux-iteration-S1-case1}, we get
\begin{align*}
	&\int_{R_2} \exp\left( \sqrt{1+S^2_1(F_{N,k})(x_1,x_2)} -c   S^2(F_{N,k})(x_1,x_2)\right)  {\rm d}\mu_2(x_2)\\
	\leq& \exp\left(\sqrt{1+A_{k-1}(x_1,c_{R_2})}  - cS^2(F_{N, k-1})(x_1, c_{R_2}) \right) \mu_2(R_2).
\end{align*}
 Therefore,   \eqref{eqn:one-shift-bi-iteration} holds, as desired.

\medskip

\noindent{\it Step II}. It remains to show \eqref{eqn:aux-S1-S-weak-0} with some constant $\alpha\in (0,1)$ for $F_{N,N}^{(I)}$. Similarly, one may assume that $\mathcal G_2^{(I)}=\bigcup_{0\leq n\leq N}Ch^{(n)}(Q_2) $, that is, $F_{N,N}^{(I)}=F_{N,N}$.

Note that for any fixed $x_1\in Q_1$,
\begin{equation}\label{eqn:small-size-bi}
    \left( \mathcal A_{P_2}\right)^2 (x_1,x_2)=\sum_{\substack{P_1\in  Ch^{(n)}(Q_1)\\ 0\leq n\leq N}}    \mathfrak{a}^2_{P_1\times P_2}(  x_1, x_2) <1,\quad \forall\, x_2\in P_2,
\end{equation} 
and
$$
    S^2 (F_{N,N})  (x_1,x_2)=\sum_{\substack{P_2\in  Ch^{(n)}(Q_2)\\ 0\leq n\leq N}}  \left( \mathcal A_{P_2}\right)^2 (x_1,x_2).
$$
Hence the ``atomic function'' $\mathcal A_{P_2}(x_1,\cdot)$   plays the role of $\mathfrak a_P$ with $\|\mathfrak a_P\|_{L^\infty}<1$ in \eqref{eqn:newCH-dyadic-2}. This motivates us to modify the approach (concretely, these starting from {\bf Steps II} therein) in the proof of  \eqref{eqn:newCH-dyadic-2} to attack 
\eqref{eqn:aux-S1-S-weak-0}.  
We now show the details in the following steps.
 
 \smallskip
 
Step II-i.  Similar to {\bf  Step II.2} in the proof of  \eqref{eqn:newCH-dyadic-2}, one can verify that 
there exist positive  constants $c'$ and $\mathsf M'$, such that 
\begin{equation}\label{eqn:bi-base-case}
    \int_{Q_2}\exp\bigg(c'  S_1(F_{N,N}) (x_1,x_2)\bigg)\, {\rm d}\mu_2(x_2)\leq \mathsf M' \mu_2(Q_2)\quad \text{whenever}\   \   \left \|S(F_{N,N})(x_1,\cdot)\right\|_{L^\infty(Q_2)}^2 \leq 10.
\end{equation}

Since $S_1(F_{N,N})$ is not defined linearly as $F_N$ in Theorem~\ref{thm:newCH-dyadic}, and to clarify the specific application of almost orthogonality \textsf{(P3)}, we will provide the detailed proof for the readers' convenience.

Concretely, it suffices to show
\begin{equation}\label{eqn:J-N-bi-base}
	\mu_2\Big(\Big\{x_2\in Q_2:\     |S_1(F_{N,N})(x_1,x_2)|>\lambda\Big\}\Big)
    \leq   c'_1 \exp\bigg(- c'_2\f{ \lambda}  {\|S(F_{N,N})(x_1,\cdot)\|_{L^\infty(Q_2)}}\bigg)\, \mu_2(Q_2)
\end{equation}
for all $\lambda \gtrsim \left\|S(F_{N,N})(x_1,\cdot)\right\|_{L^\infty(Q_2)}$.

For any fixed $x_1\in Q_1$, we apply \eqref{eqn:double-P3-revise-II} in \textsf{(P3)} to see
\begin{align}\label{eqn:P3-S1}
	\left|S_1(F_{N,N})(x_1,\cdot )\right|_{Q_2}\leq \,&\bigg(\frac{1}{\mu_2(Q_2)} \int_{Q_2}   \sum_{\substack{P_1\in  Ch^{(n)}(Q_1)\\ 0\leq n\leq N}} \Big|\sum_{\substack{P_2\in  Ch^{(n)}(Q_2)\\ 0\leq n\leq N}}  \mathfrak{a}_{P_1\times P_2}(  x_1, x_2) \Big | ^2  {\rm d}\mu_2(x_2) \bigg)^{1/2} \nonumber\\
	\leq \,&      \bigg(\frac{1}{\mu_2(Q_2)}   \sum_{\substack{P_1\in  Ch^{(n)}(Q_1)\\ 0\leq n\leq N}}    C_2^2  \int_{Q_2}   \sum_{\substack{P_2\in  Ch^{(n)}(Q_2)\\ 0\leq n\leq N}}  \big|\mathfrak{a}_{P_1\times P_2}(  x_1, x_2) \big | ^2  {\rm d}\mu_2(x_2) \bigg)^{1/2}\nonumber\\
	=\,& C_2  \left\|S(F_{N,N})(x_1,\cdot)\right\|_{L^\infty(Q_2)}.
\end{align}

For any $\lambda\geq C_2\left\|S(F_{N,N})(x_1,\cdot)\right\|_{L^\infty(Q_2)}$, consider the level set
$$
    \mathcal S_\lambda:=\Big \{x_2\in Q_2:\, S_1(F_{N,N}) (x_1,x_2)>\lambda\Big \}
$$
and assume $\mathcal S_\lambda \not=\emptyset$. Due to \eqref{eqn:P3-S1}, there exists a set $\{Q_{\lambda, i}\}_i$ of maximal dyadic cubes $Q_{\lambda, i}\subsetneqq Q_2$ such that 
\begin{equation}\label{eqn:stopping-time-bi}
	\left(S_1(F_{N,N})(x_1,\cdot) \right)_{Q_{\lambda,i}}> \lambda  \quad \text{while}\quad \left(S_1(F_{N,N})(x_1,\cdot) \right)_{\widetilde{Q}_{\lambda,i}}
         \leq \lambda,
\end{equation}
where $\widetilde{Q}_{\lambda,i}$ is the parent of $Q_{\lambda,i}$.
Then
$$
    \mathsf S_\lambda \subseteq \bigcup_i Q_{\lambda,i} =:\mathsf Q_\lambda.
$$

Next we consider the upper bound of $ \left(S_1(F_{N,N})(x_1,\cdot)\right)_{Q_{\lambda,i}}$ for any maximal cube $Q_{\lambda, i}$. 
Note that
\begin{align*} 
    \left(S_1(F_{N,N})(x_1,\cdot)\right)_{Q_{\lambda,i}} \leq \,& \Bigg( \sqrt{\sum_{\substack{P_1\in  Ch^{(n)}(Q_1)\\ 0\leq n\leq N}}    \bigg|\sum_{\substack{P_2\in Ch^{(n)}(Q_2):\,   P_2\supsetneqq \widetilde{Q}_{\lambda,i}\\ 0\leq n\leq N}} 
    \mathfrak{a}_{P_1\times P_2}(x_1,\cdot) \bigg|^2 }\, \Bigg)_{Q_{\lambda,i}} +
    \Bigg( \sqrt{\sum_{\substack{P_1\in  Ch^{(n)}(Q_1)\\ 0\leq n\leq N}}    \big| 
    \mathfrak{a}_{P_1\times \widetilde{Q}_{\lambda,i}}(x_1,\cdot) \big|^2 }\,\Bigg)_{Q_{\lambda,i}}\\
    &+\Bigg( \sqrt{\sum_{\substack{P_1\in  Ch^{(n)}(Q_1)\\ 0\leq n\leq N}}    \bigg|\sum_{\substack{P_2\in Ch^{(n)}(Q_2):\,   P_2\subseteq  Q_{\lambda,i}\\ 0\leq n\leq N}} 
    \mathfrak{a}_{P_1\times P_2}(x_1,\cdot) \bigg|^2 }\, \Bigg)_{Q_{\lambda,i}}\\
    =:\, & I+II+III.
\end{align*}
By  \textsf{(P2)},  
 the atomic functions in the term $I$ are all constant on $\widetilde{Q}_{\lambda,i}$, thus
\begin{align*}
	I =\, &   \Bigg( \sqrt{\sum_{\substack{P_1\in  Ch^{(n)}(Q_1)\\ 0\leq n\leq N}}    \bigg|\sum_{\substack{P_2\in Ch^{(n)}(Q_2):\,   P_2\supsetneqq \widetilde{Q}_{\lambda,i}\\ 0\leq n\leq N}} 
    \mathfrak{a}_{P_1\times P_2}(x_1,\cdot) \bigg|^2 }\, \Bigg)_{\widetilde Q_{\lambda,i}}\\
    =\, &   \Bigg( \sqrt{\sum_{\substack{P_1\in  Ch^{(n)}(Q_1)\\ 0\leq n\leq N}}    \bigg|\sum_{\substack{P_2\in Ch^{(n)}(Q_2)\\ 0\leq n\leq N}} 
    \mathfrak{a}_{P_1\times P_2}(x_1,\cdot) -   \sum_{\substack{P_2\in Ch^{(n)}(Q_2):\,   P_2\subseteq\widetilde{Q}_{\lambda,i}\\ 0\leq n\leq N}} 
    \mathfrak{a}_{P_1\times P_2}(x_1,\cdot) \bigg|^2 }\, \Bigg)_{\widetilde Q_{\lambda,i}}\\
    \leq \,& \big (S_1(F_{N,N}) (x_1,\cdot)\big )_{\widetilde Q_{\lambda,i}} +     \Bigg( \sqrt{\sum_{\substack{P_1\in  Ch^{(n)}(Q_1)\\ 0\leq n\leq N}}    \bigg|\sum_{\substack{P_2\in Ch^{(n)}(Q_2):\,   P_2\subseteq \widetilde{Q}_{\lambda,i}\\ 0\leq n\leq N}} 
    \mathfrak{a}_{P_1\times P_2}(x_1,\cdot) \bigg|^2 }\, \Bigg)_{\widetilde Q_{\lambda,i}}.
\end{align*}
Note that $\big (S_1(F_{N,N})(x_1,\cdot)\big )_{\widetilde Q_{\lambda,i}}\leq \lambda$ by \eqref{eqn:stopping-time-bi}. Similar to \eqref{eqn:P3-S1},  applying \eqref{eqn:double-P3-revise-II} in \textsf{(P3)} again, we have 
$$
     \Bigg( \sqrt{\sum_{\substack{P_1\in  Ch^{(n)}(Q_1)\\ 0\leq n\leq N}}    \bigg|\sum_{\substack{P_2\in Ch^{(n)}(Q_2):\,   P_2\subseteq \widetilde{Q}_{\lambda,i}\\ 0\leq n\leq N}} 
    \mathfrak{a}_{P_1\times P_2}(x_1,\cdot) \bigg|^2 }\, \Bigg)_{\widetilde Q_{\lambda,i}}\leq C_2  \left\|S(F_{N,N})(x_1,\cdot)\right\|_{L^\infty(Q_2)}.
$$
These imply that 
$$
  I\leq  \lambda + C_2  \left\|S(F_{N,N})(x_1,\cdot)\right\|_{L^\infty(Q_2)}.
$$

By definition,
$$
  II= \sqrt{\sum_{\substack{P_1\in  Ch^{(n)}(Q_1)\\ 0\leq n\leq N}}    \big| 
    \mathfrak{a}_{P_1\times \widetilde{Q}_{\lambda,i}}(x_1,c_{Q_{\lambda,i}}) \big|^2 } \leq \left\|S(F_{N,N})(x_1,\cdot)\right\|_{L^\infty(Q_2)},
$$
and moreover, by \eqref{eqn:double-P3-revise-II} in \textsf{(P3)} again,
$$
     III\leq C_2  \left\|S(F_{N,N})(x_1,\cdot)\right\|_{L^\infty(Q_2)}.
$$

Hence, combining the estimates of $I$, $II$ and $III$, we have
\begin{equation}\label{eqn:auxAverg1-bi}
    \left(S_1(F_{N,N})(x_1,\cdot)\right)_{Q_{\lambda,i}} \leq \lambda +(2C_2+1) \left\|S(F_{N,N})(x_1,\cdot)\right\|_{L^\infty(Q_2)}.
\end{equation}
       
For any $\lambda'>0$, let $\big\{Q_{\lambda+\lambda',j}\big \}$ be the maximal dyadic cubes such that
$\left( S_1(F_{N,N})(x_1,\cdot) \right)_{Q_{\lambda+\lambda',j}}>\lambda+\lambda'$. Also denote $\mathsf Q_{\lambda+\lambda'}:= \bigcup_j Q_{\lambda+\lambda',j} $,
obviously, $\mathsf Q_{\lambda+\lambda'}\subseteq \mathsf Q_{\lambda}$.

For every fixed $Q_{\lambda,i}$,
\begin{align*} 
   & \mu_2\left(\mathsf Q_{\lambda+\lambda'}\cap Q_{\lambda,i}\right)\\
    \leq & \f{1}{\lambda+\lambda'} \sum_{j:\  Q_{\lambda+\lambda',j}
    \subset Q_{\lambda,i}}  \int_{Q_{\lambda+\lambda',j}}
    S_1(F_{N,N})(x_1,x_2) \,{\rm d}\mu_2(x_2) \\
    \leq &  \f{1}{\lambda+\lambda'}  \int_{Q_{\lambda,i}} \left| S_1(F_{N,N})(x_1,x_2)
    -\big ( S_1(F_{N,N})(x_1,\cdot)\big )_{Q_{\lambda,i}}\right|  {\rm d}\mu_2(x_2) +\f{ \big ( S_1(F_{N,N})(x_1,\cdot)\big )_{Q_{\lambda,i}}}{\lambda+\lambda'}
    \mu_2(\mathsf Q_{\lambda+\lambda'}\cap Q_{\lambda,i}).
\end{align*}

For any $x_2\in Q_{\lambda, i}$, by \textsf{(P2)} again,
\begin{align*}
	&\left| S_1(F_{N,N})(x_1,x_2)
    -\big ( S_1(F_{N,N})(x_1,\cdot)\big )_{Q_{\lambda,i}}\right|\\
    \leq & \frac{1}{\mu_2(Q_{\lambda, i})}\int_{Q_{\lambda, i} }   \sqrt{   \sum_{\substack{P_1\in  Ch^{(n)}(Q_1)\\ 0\leq n\leq N}}   \Big|\sum_{\substack{P_2\in  Ch^{(n)}(Q_2)\\ 0\leq n\leq N}}  \mathfrak{a}_{P_1\times P_2}(  x_1, x_2)-\sum_{\substack{P_2\in  Ch^{(n)}(Q_2)\\ 0\leq n\leq N}}  \mathfrak{a}_{P_1\times P_2}(  x_1, y_2) \Big | ^2    }  \,{\rm d}\mu_2(y_2)\\
    \leq & \frac{1}{\mu_2(Q_{\lambda, i})}\int_{Q_{\lambda, i} }   \sqrt{   \sum_{\substack{P_1\in  Ch^{(n)}(Q_1)\\ 0\leq n\leq N}}   \Big|\sum_{\substack{P_2\in  Ch^{(n)}(Q_2):\, P_2\subseteq Q_{\lambda, i}\\ 0\leq n\leq N}}  \mathfrak{a}_{P_1\times P_2}(  x_1, x_2)-\sum_{\substack{P_2\in  Ch^{(n)}(Q_2):\, P_2\subseteq Q_{\lambda, i}\\ 0\leq n\leq N}}  \mathfrak{a}_{P_1\times P_2}(  x_1, y_2) \Big | ^2    }  \,{\rm d}\mu_2(y_2)\\
    \leq &    \sqrt{   \sum_{\substack{P_1\in  Ch^{(n)}(Q_1)\\ 0\leq n\leq N}}   \Big|\sum_{\substack{P_2\in  Ch^{(n)}(Q_2):\, P_2\subseteq Q_{\lambda, i}\\ 0\leq n\leq N}}  \mathfrak{a}_{P_1\times P_2}(  x_1, x_2) \Big | ^2    }      +\Bigg(\sqrt{   \sum_{\substack{P_1\in  Ch^{(n)}(Q_1)\\ 0\leq n\leq N}}   \Big|\sum_{\substack{P_2\in  Ch^{(n)}(Q_2):\, P_2\subseteq Q_{\lambda, i}\\ 0\leq n\leq N}}  \mathfrak{a}_{P_1\times P_2}(  x_1, \cdot) \Big | ^2    }  \Bigg)_{Q_{\lambda ,i}}.
\end{align*}
Hence,
\begin{align*}
   & \mu_2\left(\mathsf Q_{\lambda+\lambda'}\cap Q_{\lambda,i}\right)\\
    \leq & \frac{2}{\lambda+\lambda'}\Bigg(\sqrt{   \sum_{\substack{P_1\in  Ch^{(n)}(Q_1)\\ 0\leq n\leq N}}   \Big|\sum_{\substack{P_2\in  Ch^{(n)}(Q_2):\, P_2\subseteq Q_{\lambda, i}\\ 0\leq n\leq N}}  \mathfrak{a}_{P_1\times P_2}(  x_1, \cdot) \Big | ^2    }  \Bigg)_{Q_{\lambda ,i}} \mu_2(Q_{\lambda,i})+  \f{ \big ( S_1(F_{N,N})(x_1,\cdot)\big )_{Q_{\lambda,i}}}{\lambda+\lambda'}
    \mu_2(\mathsf Q_{\lambda+\lambda'}\cap Q_{\lambda,i}).
\end{align*}
For the first term in the last inequality above, we first apply H\"older's inequality and then 
 \eqref{eqn:double-P3-revise-II} in \textsf{(P3)}. Then combining with \eqref{eqn:auxAverg1-bi}, we see that 
\begin{align*}
   & \mu_2\left(\mathsf Q_{\lambda+\lambda'}\cap Q_{\lambda,i}\right)\\
    \leq\, & \f{2C_2 \left\|S(F_{N,N})(x_1,\cdot)\right\|_{L^\infty({Q_2})}}
    {\lambda+\lambda'} \mu_2(Q_{\lambda,i}) +\f{\lambda+(2C_2+1)\left\|S(F_{N,N})(x_1,\cdot)\right\|_{L^\infty({Q_2})}
    }{\lambda+\lambda'} \mu_2(\mathsf Q_{\lambda+\lambda'}\cap Q_{\lambda,i}).
\end{align*}
That is, the analog of \eqref{em1} holds for our purpose here. The remaining argument for \eqref{eqn:J-N-bi-base} is similar to that in {\bf  Step II.2} of  \eqref{eqn:newCH-dyadic-2}, and we skip the details.

\medskip

 Step II-ii. 
 Recall that the assumption \eqref{eqn:small-size-bi} made in this case plays a role analogous to the small-scale assumption  $\|\mathfrak a_P\|_{L^\infty}<1$  in \eqref{eqn:newCH-dyadic-2}. Similarly,
 using \eqref{eqn:bi-base-case} as the base case, one can apply an analogy of the classification algorithm given in {\bf  Step II.1} of  \eqref{eqn:newCH-dyadic-2} to obtain
$$
     \mathcal G_2:=\bigcup_{0\leq n\leq N} Ch^{(n)}(Q_2)=\sum_k \mathcal G_{2,k}.
$$
We begin by constructing $\mathcal G_{2,1}$. 
Let $\left\{P^{(1)}_j\right\}_j$ be  dyadic cubes in $\mathcal G_2$
 such that 
$$
    \sum_{P\in \mathcal G_2,\, P\supsetneqq P'} \big|\mathcal A_{P}(x_1,\cdot)\big|^2+   \max_{x_2} \big|\mathcal A_{P'}(x_1,x_2)\big|^2  <9 \ \     
    \text{on}\  \    P',\   \  \forall\, P'\in \mathcal G_2\  \  \text{with}\  \  P'\supsetneqq P_j^{(1)},
  \leqno{(\mathcal C1)}   
$$
 while 
 $$
    \sum_{P\in \mathcal G,\, P\supsetneqq P^{(1)}_j} \big|\mathcal A_{P}(x_1,\cdot)\big|^2+  \max_{x_2} \big|\mathcal A_{P^{(1)}_j }(x_1,x_2)\big|^2 \geq 9 \ \     \text{on}\  \     P^{(1)}_j .   \leqno{(\mathcal C2)}
$$
Now denote
$$
    \mathcal G_{2,1}=\left\{P\in \mathcal G_2:\, P\supseteq P^{(1)}_j \  \text{for \  some\ } j   \ \text{or}\  \   P\cap  \big(\cup_j P^{(1)}_j \big )=\emptyset  
    \right\}\quad \text{and}\quad  F_{N,N}^{(1)}:=\sum_{\substack{P_1\in  Ch^{(n)}(Q_1)\\ 0\leq n\leq N}} \sum_{ P_2\in \mathcal G_{2,1}} \mathfrak{a}_{P_1\times P_2}.
$$
Similar to the one-parameter version in the proof of \eqref{eqn:newCH-dyadic-2}, one may verify  the non-conflict and existence of $\mathcal G_{2,1}$.
Then
$$
   S_1(F_{N,N}^{(1)}) =\sqrt{\sum_{\substack{P_1\in  Ch^{(n)}(Q_1)\\ 0\leq n\leq N}} \bigg|\sum_{ P_2\in \mathcal G_{2,1}} \mathfrak{a}_{P_1\times P_2}  \bigg|^2  }\quad 
   \text{and}  \quad   S^2(F_{N,N}^{(1)}) =\sum_{\substack{P_1\in  Ch^{(n)}(Q_1)\\ 0\leq n\leq N}} \sum_{ P_2\in \mathcal G_{2,1}} \big|\mathfrak{a}_{P_1\times P_2} \big|^2 .
$$
Due to the assumption $\|\mathcal A_{P_2}(x_1,\cdot)\|_{L^\infty(P_2)}<1$,
it's clear that 
$$
    S^2(F_{N,N}^{(1)}) (x_1,x_2)\leq10 \  \   \text{for}\  \   x_2\in Q_2.
$$
Hence the base case \eqref{eqn:bi-base-case} holds when replacing $F_{N,N}$ therein by $F_{N,N}^{(1)}$.

For $k\geq 2$, we construct $\mathcal G_{2,k}$ from $\mathcal G_{2}\setminus (\bigcup_{i=1}^{k-1}\mathcal G_{2,i})$ in a similar manner. That is,  let $\left\{P^{(k)}_j\right\}_j$ be  dyadic cubes in $  \mathcal G_{2}\setminus (\bigcup_{i=1}^{k-1}\mathcal G_{2,i})$ 
 such that 
$$
    \sum_{P\in \mathcal G_{2}\setminus (\bigcup_{i=1}^{k-1}\mathcal G_{2,i}),\, P\supsetneqq P'} \big|\mathcal A_{P}(x_1,\cdot)\big|^2+   \max_{x_2} \big|\mathcal A_{P'}(x_1,x_2)\big|^2  <9 \ \     
    \text{on}\  \    P',\   \  \ \ \forall\, P'\in \mathcal G_{2}\bigg\backslash \Big (\bigcup_{i=1}^{k-1}\mathcal G_{2,i}\Big ) \   \text{with}\    P'\supsetneqq P_j^{(k)}, 
$$
while 
$$
    \sum_{P\in \mathcal G_{2}\setminus (\bigcup_{i=1}^{k-1}\mathcal G_{2,i}),\, P\supsetneqq P_j^{(k)} } \big|\mathcal A_{P}(x_1,\cdot)\big|^2+  \max_{x_2} \big|\mathcal A_{P_j^{(k)}  }(x_1,x_2)\big|^2 \geq 9 \ \     \text{on}\  \     P_j^{(k)} .   
$$
Now denote for $k\geq 2$,
$$
    \mathcal G_{2,k}=\left\{P\in \mathcal G_{2}\setminus \Big (\bigcup_{i=1}^{k-1}\mathcal G_{2,i}\Big ):\, P\supseteq P^{(k)}_j \  \text{for \  some\ } j   \ \text{or}\  \   P\cap  \big(\cup_j P^{(k)}_j \big )=\emptyset  
    \right\}
$$
and
$$
    F_{N,N}^{(k)}:=\sum_{\substack{P_1\in  Ch^{(n)}(Q_1)\\ 0\leq n\leq N}} \sum_{ P_2\in \mathcal G_{2,k}} \mathfrak{a}_{P_1\times P_2}.
$$

Note that 
$$
    S_1(F_{N,N})\leq \sum_k S_1(F_{N,N}^{(k)}) \quad \text{and}\quad  S^2(F_{N,N})=\sum_k S^2(F_{N,N}^{(k)})=\sum_k \sum_{P_2\in \mathcal G_{2,k}}  \big(\mathcal A_{P_2}\big)^2.
$$
Again, the construction of $\mathcal G_{2,k}$ yields that 
$\displaystyle  \left\|S(F_{N,N}^{(k)}) (x_1,\cdot)\right\|_{L^\infty (Q_2)} ^2\leq10$ and 
$$
     1\leq k\leq \mathsf K_2:=\left\lfloor \frac{\|S(F_{N,N})(x_1,\cdot)\|_{L^\infty(Q_2)}^2}{8}\right\rfloor+1.
$$

For each $k\geq 2$,  since $\sum_{i=1}^{k-1}S_1(F_{N,N}^{(i)})$ are constant on each maximal cube in $\mathcal G_{2,k}$,  by the iteration  in {\bf  Step II} of 
\eqref{eqn:newCH-dyadic-2}, it follows from \eqref{eqn:bi-base-case} to see
$$
        \int_{Q_2}\exp\bigg(c'  \sum_{i=1}^{k}S_1(F_{N,N}^{(i)}) (x_1,x_2)\bigg)\, {\rm d}\mu_2(x_2)\leq \mathsf M'   \int_{Q_2}\exp\bigg(c'  \sum_{i=1}^{k-1}S_1(F_{N,N}^{(i)}) (x_1,x_2)\bigg)\, {\rm d}\mu_2(x_2).
$$
As a consequence, via iterating on $k$, we have
$$
    \int_{Q_2}\exp\bigg(c'   S_1(F_{N,N}) (x_1,x_2)\bigg)\, {\rm d}\mu_2(x_2)\leq  \int_{Q_2}\exp\bigg(c'  \sum_{i=1}^{\mathsf K_2}S_1(F_{N,N}^{(i)}) (x_1,x_2)\bigg)\, {\rm d}\mu_2(x_2)\leq  (\mathsf M')^{\mathsf K_2}   \mu_2(Q_2).
$$
Hence, there exist positive  constants $C$ and $C'$ such that 
\begin{equation}\label{eqn:aux-S1-S-weak}
    \int_{Q_2}\exp\bigg(c' S_1(F_{N,N}) (x_1,x_2)\bigg)\, {\rm d}\mu_2(x_2)\leq C  \exp\left(C'  \|S(F_{N,N})(x_1,\cdot)\|_{L^\infty(Q_2)}^2\right) \mu_2(Q_2).
\end{equation}

From the above, we complete the proof of \eqref{eqn:aux-S1-S-weak-0}.
\hfill $\Box$

\bigskip

Up to now, we have proved all estimates and implications in \eqref{F-a}. We turn to considering its counterpart \eqref{F-b}. We start by  proving Proposition \ref{prop:double-Lp-CH} based on the Chang--Wilson--Wolff type estimate \eqref{eqn:newCWW-S1-S-0}.

\medskip

{\bf Proof of Proposition \ref{prop:double-Lp-CH}:}
Fix $x_1\in Q_1$.  We apply Lemma~\ref{lem:local-CWW-to-Lp} in the second variable with the Hilbert space
$$
 \mathcal H=\ell^2(\mathcal P_1),
 \qquad \mathcal P_1=\bigcup_{0\leq n\leq N}Ch^{(n)}(Q_1).
$$
Let
$$
 \mathcal P_2=\bigcup_{0\leq n\leq N}Ch^{(n)}(Q_2).
$$
For each $P_2\in\mathcal P_2$, define
$$
 A_{P_2}(x_2)
 =\{\mathfrak a_{P_1\times P_2}(x_1,x_2)\}_{P_1\in\mathcal P_1}\in\mathcal H.
$$
An entry is understood to be zero if the corresponding rectangle does not occur in $F_{N,N}$.  By \textsf{(P1)} and
\textsf{(P2)}, $A_{P_2}$ is supported on $P_2$ and is constant on every child of $P_2$.  Moreover,
$$
 \left\|\sum_{P_2}A_{P_2}(x_2)\right\|_{\mathcal H}
 =\left(\sum_{P_1\in\mathcal P_1}
 \left|\sum_{P_2\in\mathcal P_2}\mathfrak a_{P_1\times P_2}(x_1,x_2)\right|^2\right)^{1/2}
 =S_1(F_{N,N})(x_1,x_2),
$$
while
$$
 \left(\sum_{P_2}\|A_{P_2}(x_2)\|_{\mathcal H}^2\right)^{1/2}
 =\left(\sum_{P_2\in\mathcal P_2}\sum_{P_1\in\mathcal P_1}
 |\mathfrak a_{P_1\times P_2}(x_1,x_2)|^2\right)^{1/2}
 =S(F_{N,N})(x_1,x_2).
$$

We next verify the local exponential-square hypothesis.  Fix a dyadic cube $R_2\subseteq Q_2$ and a subfamily
$$
 \mathcal E\subseteq\{P_2\in\mathcal P_2:P_2\subseteq R_2\}.
$$
The assertion is immediate if $\mathcal E$ is empty.  Otherwise, choose $L$ so that
$$
 \mathcal E\subseteq\bigcup_{0\leq k\leq L}Ch^{(k)}(R_2).
$$
Put
$$
 \mathcal P_2(R_2,L)=\bigcup_{0\leq k\leq L}Ch^{(k)}(R_2)
$$
and define the zero-padded family
$$
 \widetilde A_{P_2}=
 \begin{cases}
 A_{P_2},&P_2\in\mathcal E,\\
 0,&P_2\in\mathcal P_2(R_2,L)\setminus\mathcal E.
 \end{cases}
$$
This is a finite family based on the top cube $R_2$.  Its nonzero entries are exactly those indexed by
$\mathcal E$, while the omitted generations and cubes contribute zero.
For $P_1\in\mathcal P_1$, restrict the rectangle family to those $P_1\times P_2$ with $P_2\in\mathcal E$.  If
$x_1\notin P_1$, all the corresponding atoms vanish.  If $x_1\in P_1$, then for every dyadic cube $W_2\subseteq R_2$
and every subfamily $\mathcal F\subseteq\{P_2\in\mathcal E:P_2\subseteq W_2\}$,
\begin{align*}
 &\left(\int_{W_2}\left|\sum_{P_2\in\mathcal F}
 \mathfrak a_{P_1\times P_2}(x_1,x_2)\right|^2\,{\rm d}\mu_2(x_2)\right)^{1/2}
 \leq C_2\left(\int_{W_2}\sum_{P_2\in\mathcal F}
 |\mathfrak a_{P_1\times P_2}(x_1,x_2)|^2\,{\rm d}\mu_2(x_2)\right)^{1/2}.
\end{align*}
This is precisely \eqref{eqn:double-P3-revise-II} for the restricted family.  Squaring this estimate and summing over
$P_1\in\mathcal P_1$ gives
$$
 \left(\int_{W_2}\left\|\sum_{P_2\in\mathcal F}A_{P_2}(x_2)\right\|_{\mathcal H}^2
 \,{\rm d}\mu_2(x_2)\right)^{1/2}
 \leq C_2\left(\int_{W_2}\sum_{P_2\in\mathcal F}\|A_{P_2}(x_2)\|_{\mathcal H}^2
 \,{\rm d}\mu_2(x_2)\right)^{1/2}.
$$
For an arbitrary subfamily of $\mathcal P_2(R_2,L)$, deleting the indices outside $\mathcal E$ changes neither side
of this estimate.  Hence the same estimate holds for the zero-padded family.
In particular, no factor depending on the dimension of $\mathcal H$ occurs.  The support and child-constancy conditions
for $\{\widetilde A_{P_2}\}_{P_2\in\mathcal P_2(R_2,L)}$ are also inherited from \textsf{(P1)} and \textsf{(P2)}.
Set
$$
 G_{\mathcal E}=\sum_{P_2\in\mathcal E}A_{P_2},
 \qquad
 S_{\mathcal E}=\left(\sum_{P_2\in\mathcal E}\|A_{P_2}\|_{\mathcal H}^2\right)^{1/2}.
$$
If $\|S_{\mathcal E}\|_{L^\infty(R_2)}=0$, then the preceding Hilbert-valued \textsf{(P3)} estimate gives
$G_{\mathcal E}=0$ almost everywhere on $R_2$, and the assertion is immediate.  Otherwise, we may repeat the proof of
\eqref{eqn:newCWW-S1-S-0} on $R_2$.  Since that proof uses the preceding estimate with the same constant $C_2$, it gives
$$
 \frac1{\mu_2(R_2)}\int_{R_2}\exp\left(
 c_3'\frac{\|G_{\mathcal E}(x_2)\|_{\mathcal H}^2}
 {\|S_{\mathcal E}\|_{L^\infty(R_2)}^2}
 \right)\,{\rm d}\mu_2(x_2)\leq c_4'.
$$
The constants are independent of $R_2$, $\mathcal E$, $x_1$, and the cardinality of $\mathcal P_1$.  Thus
\eqref{eqn:local-CWW-hilbert} holds for the family $\{A_{P_2}\}_{P_2\in\mathcal P_2}$.  Applying
Lemma~\ref{lem:local-CWW-to-Lp}, together with the two identities above, gives
$$
 \|S_1(F_{N,N})(x_1,\cdot)\|_{L^p(Q_2)}
 \leq C\sqrt p\,
 \|S(F_{N,N})(x_1,\cdot)\|_{L^p(Q_2)},\qquad p\geq2.
$$
This is \eqref{eqn:part-II-Lp-CH}.
\hfill $\Box$

\medskip

Consequently, we will obtain Proposition \ref{prop:one-shift-bi-I} straightforwardly.

\medskip

{\bf Proof of Proposition \ref{prop:one-shift-bi-I}:} 

Due to the constant for sharp $L^p$ square function estimates in Proposition \ref{prop:double-Lp-CH} is independent of $x_1\in Q_1$ and $p\geq 2$, thus the argument for Proposition \ref{prop:one-shift-bi-I} is similar to that of \eqref{eqn:newCH-dyadic}  in Theorem \ref{thm:newCH-dyadic}, and  $S_1(F_{N,N})$, $S(F_{N,N})$ in  Proposition \ref{prop:one-shift-bi-I} are the corresponding roles of $F_N$,  $S(F_N)$ in \eqref{eqn:newCH-dyadic}, respectively. We omit the details.
\hfill $\Box$

\smallskip

\begin{remark}\label{rem:global-omega}
Similarly, the arguments in Sections \ref{sec:3} and \ref{sec:4} above  can be extended to the whole space. Concretely, let  $\{\mathfrak a_{P_1\times P_2}: \, P_1\times P_2\in \mathcal D^{[1]}\times \mathcal D^{[2]} \} $ be atomic functions satisfying  the assumption of Theorem \ref{thm:double-newCH-dyadic} when fixing on  any given $Q_1\times Q_2$. Define for $M,N\in \mathbb N$,
$$
     \mathbf F_{M,N}(x_1,x_2):=\sum_{P_1\times P_2\in \bigcup_{i=-M}^N    \mathcal D_i^{[1]}  \times  \bigcup_{j=-M}^N     \mathcal D_j^{[2]}} \mathfrak a_{P_1\times P_2}(x_1,x_2).
$$
Then for $p\geq 2$,
$$
    \|\mathbf F_{M,N}\|_{L^p(X_1\times X_2)}\leq C \sqrt{p}\ \|S_1(\mathbf  F_{M,N})\|_{L^p(X_1\times X_2)}
$$
and 
$$
      \|S_1(\mathbf F_{M,N})\|_{L^p(X_1\times X_2)}\leq C \sqrt{p}\ \|S(\mathbf F_{M,N})\|_{L^p(X_1\times X_2)}.
$$
The proof is rather standard based on our Theorems \ref{thm:newCH-dyadic} and   \ref{thm:double-newCH-dyadic}, and we skip the details.

 In particular, if $S(\mathbf F_{M,N})\in L^\infty(X_1\times X_2)$, then we also have $S_1(\mathbf F_{M,N})\in L^\infty(X_1\times X_2)$ by noting $S_1(\mathbf F_{M,N})\leq \sqrt{(M+N+1)  \big  ((M+N+1)  \|S(\mathbf F_N)\|_{L^\infty}\big  )^2   }$. Therefore, for any given open set $\Omega\subseteq X_1\times X_2$ with $\mu(\Omega)<\infty$, set
$$
     \mathbf F_{\Omega, M,N} (x_1,x_2):=\sum_{\substack{ P_1\times P_2\in \bigcup_{i=-M}^N
     \mathcal D_i^{[1]}      \times  \bigcup_{j=-M}^N     \mathcal D_j^{[2]}  \\  P_1\times P_2\subseteq \Omega }} \mathfrak a_{P_1\times P_2}(x_1,x_2).
$$
We  then have for $p\geq 2$,
$$
    \|\mathbf F_{\Omega, M,N}\|_{L^p(\Omega)}=\|\mathbf F_{\Omega,M,N}\|_{L^p(X_1\times X_2)}\leq C \sqrt{p}\ \|S_1(\mathbf F_{\Omega, M,N})\|_{L^p(X_1\times X_2)}\leq C \sqrt{p}\ \|S_1(\mathbf F_{M,N})\|_{L^p(\Omega)}
$$
and 
$$
      \|S_1(\mathbf F_{\Omega, M,N})\|_{L^p(\Omega)}\leq C \sqrt{p}\ \|S(\mathbf F_{M,N})\|_{L^p(\Omega)}.
$$
This, in combination with the implications (replacing $Q_1\times Q_2$ therein by $\Omega$) in \eqref{F-a} and \eqref{F-b}, deduces that 
\begin{align}\label{eqn:one-shift-bi-I-Omega}
   	&\int_{\Omega} \exp\Big (  \alpha S_1 \left (\mathbf F_{\Omega, M,N} \right )(x_1,x_2)  \Big )\,  {\rm d}\mu_1(  x_1) \,  {\rm d}\mu_2(x_2)\nonumber\\
   	\leq  &  C\int_{\Omega} \exp\Big (  c \alpha^2 S^2\left ( \mathbf F_{\Omega, M,N}\right )(x_1,x_2)\Big )\,  {\rm d}\mu_1(  x_1) \,{\rm d}\mu_2(x_2)
\end{align}
for any $\alpha>0$, 
and
\begin{equation}\label{eqn:double-newCH-dyadic-Omega}
 \int_{\Omega} \exp\Bigg(c' \frac{   \left| \mathbf F_{\Omega, M,N} ( x_1, x_2)\right|   }{\left\|S(\mathbf F_{M,N})\right \|_{L^\infty(X_1\times X_2)}}\Bigg)  \,{\rm d}\mu_1(  x_1) \, {\rm d}\mu_2( x_2) \leq C' \mu(\Omega).
\end{equation}
\end{remark}

\medskip

From the above, we have completed the proofs of estimates and their implications listed in \eqref{F-a} and \eqref{F-b}, concerning the bi-parameter theory of Hoeffding's inequalities, Chang--Wilson--Wolff type result and sharp square function estimates.
In particular, double dyadic martingale differences satisfy the assumptions \textsf{(P1)}, \textsf{(P2)} and \textsf{(P3)} assumed in this article, hence our results answered the question raised in Section \ref{sec:1} affirmatively:  
we  can establish \eqref{eqn:double-ESI intro} in the higher dimensional space $\mathbb R^{\mathrm d_1}\times \mathbb R^{\mathrm d_2}$ and in more general settings such as products of spaces of homogeneous type.

\smallskip
\begin{remark}
As previously discussed, to extend the one-parameter results to the bi-parameter analogs,  the intermediate square function $S_1$ should be of the form presented in this article, where its definition avoids the use of the $L^\infty$ norm of atomic functions. Likewise, the bi-parameter type square function $S(F_{N,N})$ should also exclude the $L^\infty$ structure when applied to investigate the 3-parameter theory.  
\end{remark}

\begin{remark}
Specifically, if we only focus on the bi-parameter setting,  we may also consider the square function defined by $\widetilde S^2(\sum_{P_1}\sum_{P_2}\Delta_{P_1\times P_2}(f)):=\sum_{P_1}\sum_{P_2} \|\Delta_{P_1\times P_2}(f)\|_{L^\infty(P_1\times P_2)}^2$ associated to dyadic martingale differences $\Delta_{P_1\times P_2}(f)$ for some function $f$. There will be more tools to use, due to $\widetilde S \geq S$ pointwise and the cancellation property of $\Delta_{P_1\times P_2}(f)$. One may refer to the second part of the appendix in the last section.
\end{remark}

\medskip

At the end of this section, we further discuss the multi-parameter version of the above results, which shows that the question raised in Section \ref{sec:1} can further be extended to multi-parameter settings, and the corresponding Hoeffding's inequalities, Chang--Wilson--Wolff type result  and sharp square function estimates also hold.

 \begin{proposition}   \label{prop:multi-newCH-dyadic} 
 Let $n\geq 2$.
 For   $j=1,\ldots, n$, let  $(X_j,\rho_j,\mu_j)$ be a quasi-metric space satisfying the volume doubling condition and let $\mathcal D^{[j]}$ be a system of dyadic cubes on $X_j$. For any given $Q_1\times\cdots\times Q_n\in \mathcal D^{[1]}\times\cdots \times \mathcal D^{[n]}$, let $\left \{\mathfrak a_{P_1\times \cdots\times P_n}:\, P_1\times\cdots \times  P_n\in \mathcal D^{[1]}\times\cdots \times  \mathcal D^{[n]}, \,  P_j\subseteq Q_j,\, j=1,\ldots,n\right\}$ be ``atomic functions'' satisfying
 \begin{itemize}
	\item[\textrm   {(P1)} ] $\mathfrak{a}_{P_1\times \cdots\times P_n}$ is supported on $P_1\times\cdots\times  P_n$.
	\smallskip
	\item[\textrm {(P2)} ]  $\mathfrak{a}_{P_1\times \cdots\times P_n}$ is constant on each sub-rectangle in $ Ch(P_1\times\cdots \times  P_n):=Ch(P_1)\times \cdots \times Ch(P_n)$.	
	\smallskip
	\item[\textrm {(P3)} ]  \textbf{{\rm(}almost orthogonality{\rm)}} there exists  $C>0$ such that for any dyadic rectangle $R_1\times \cdots  \times  R_n\subseteq Q_1\times \cdots \times Q_n$ and any subcollection 
	$\mathcal G_{R_1\times\cdots \times  R_n}\subseteq \{P_1\times \cdots \times P_n\in \mathcal D^{[1]}\times \cdots \times \mathcal D^{[n]} :\, P_1\times \cdots \times P_n \subseteq R_1\times \cdots\times  R_n\}$, 
\begin{align}\label{eqn:multi-P3-revise-I}
	& \bigg(\int_{R_1} \bigg |   \sum_{P_1}  \sum_{P_2,\ldots, P_n:\, P_1\times \cdots \times P_n\in \mathcal G_{R_1\times \cdots \times R_n}} \mathfrak a_{P_1\times\cdots\times  P_n} (x_1,\ldots, x_n)\bigg |^2  {\rm d}\mu_1(x_1)\bigg)^{1/2}\\
	& \leq\, C  \, \bigg(\int_{R_1}   \sum_{P_1}  \bigg | \sum_{P_2,\ldots, P_n:\, P_1\times \cdots\times P_n\in \mathcal G_{R_1\times\cdots \times  R_n}} \mathfrak a_{P_1\times\cdots \times  P_n} (x_1,\ldots, x_n)\bigg |^2  {\rm d}\mu_1(x_1)\bigg)^{1/2}\nonumber
\end{align}
for each   fixed $x_j\in R_j$ with $2\leq  j\leq n$, and 
\begin{align}\label{eqn:multi-P3-revise-II}
	&\bigg(\int_{R_j} \bigg |    \sum_{P_j,\ldots,P_n:\, P_1\times \cdots \times P_n\in \mathcal G_{R_1\times \cdots\times R_n}} \mathfrak a_{P_1\times \cdots\times P_n} (x_1,\ldots, x_n) \bigg |^2  {\rm d}\mu_j(x_j)\bigg)^{1/2}\\
	&\leq\,  C \,  \bigg(\int_{R_j}     \sum_{P_j} \bigg |  \sum_{P_{j+1},\ldots, P_n:\, P_1\times P_2\in \mathcal G_{R_1\times \cdots \times R_n}} \mathfrak a_{P_1\times \cdots\times P_n} (x_1,\ldots, x_n)\bigg|^2  {\rm d}\mu_j(x_j)\bigg)^{1/2}, \qquad \forall j=2,\ldots, n \nonumber
\end{align}
for any fixed $P_1,\cdots, P_{j-1}$ and  $x_i$ with $i\neq j$.
\end{itemize}
Denote
$$
 F_{N,\ldots,N}(  x_1, \ldots, x_n)=\sum_{j=1}^n\sum_{\substack{P_j\in  Ch^{(n_j)}(Q_j)\\ 0\leq n_j\leq N}}  \mathfrak{a}_{P_1\times\cdots\times  P_n} ( x_1,\ldots,  x_n)
$$
and 
$$
S(F_{N,\ldots,N})( x_1, \ldots,  x_n) =\Bigg(\sum_{j=1}^n\sum_{\substack{P_j\in  Ch^{(n_j)}(Q_j)\\ 0\leq n_j\leq N}}    { \mathfrak{a}^2_{P_1\times \cdots\times P_n}  (x_1,\ldots, x_n)     }\Bigg)^{1/2}. 
$$
Then the following estimates hold.
 \begin{itemize}

 \item[ (a) ] \textbf{{\rm (}sharp $L^p$ lower bound of $S(F_{N,\ldots, N})${\rm).}} \  \    
 There exists a constant $C=C(n)>0$ depending on $n$ such that for  $p\geq 2$,
\begin{equation}\label{eqn:multi-Lp-CH-intro}
	   \left\|F_{N,\ldots, N}\right\|_{L^p(Q_1\times\cdots\times  Q_n,\,  {\rm d}\mu_1\cdots {\rm d}\mu_n)}\leq C  p^{n/2}\left\|S(F_{N,\ldots,N})\right\|_{L^p(Q_1\times\cdots\times  Q_n,   \,  {\rm d}\mu_1\cdots {\rm d}\mu_n)  }.
\end{equation}

\smallskip

\item[(b) ]  \textbf{{\rm (}multi-parameter Hoeffding inequalities{\rm).}} \  \  
There exist constants $c>1$ and $C>0$ such that 
	\begin{align}\label{eqn:one-shift-bi-multi-1}
   &	\int_{Q_1\times \cdots\times Q_{n-1}} \exp\left(  \left| \alpha  F_{N,\ldots, N}(x_1,\ldots, x_n) \right|^{2/n} \right) {\rm d}\mu_1(x_1)\cdots {\rm d}\mu_{n-1}(x_{n-1}) \nonumber\\
   \leq C &   \int_{Q_1\times \cdots\times Q_{n-1}} \exp\left(  c \alpha^2 S_{n-1}^2(F_{N,\ldots,N})(x_1,\ldots,x_n)\right)  {\rm d}\mu_1(x_1)\cdots {\rm d}\mu_{n-1}(x_{n-1})
   \end{align}
for any fixed $x_n\in Q_n$  and $\alpha>0$,
and 
\begin{equation}\label{eqn:one-shift-bi-multi-2}
   	\int_{Q_n} \exp\left(  \alpha  S_{n-1}(F_{N,\ldots, N})(x_1,\ldots, x_n)  \right){\rm d}\mu_n(x_n)\leq C\int_{Q_n} \exp\left(  c \alpha^2 S_n ^2(F_{N,\ldots,N})(x_1,\ldots,x_n)\right){\rm d}\mu_n(x_n)
\end{equation}
for any fixed $x_j\in Q_j$ with $1\leq j\leq n-1$ and $\alpha>0$, where
$$
   S_{n-1}(F_{N,\ldots, N})(x_1,\ldots, x_n)=\Bigg( \sum_{i=1}^{n-1}\sum_{\substack{P_i\in  Ch^{(m_i)}(Q_i)\\ 0\leq m_i\leq N}}   \Bigg|  \sum_{\substack{P_n\in  Ch^{(m_n)}(Q_n)\\ 0\leq m_n\leq N}}    { \mathfrak{a}_{P_1\times \cdots\times P_n}  (x_1,\ldots, x_n)  \Bigg|^2   }\Bigg)^{1/2}.
$$
\smallskip

\item[ (c) ] \textbf{{\rm (}multi-parameter Chang--Wilson--Wolff type estimate{\rm).}} \  \  
There exist positive constants $\widetilde c$ and $\widetilde C$ independent of $Q_j$ ($1\leq j\leq n$)  such that for any $N\in \mathbb N_+$,
\begin{equation}\label{eqn:multi-newCH-dyadic}
 \int_{Q_1\times\cdots \times  Q_n} \exp\Bigg(\widetilde c \frac{   | F_{N,\ldots, N}( x_1,\ldots, x_n)|^{2/n}   }{\|S(F_{N,\ldots ,N})\|_{L^\infty(Q_1\times \cdots\times Q_n)}^{2/n}}\Bigg)  \,{\rm d}\mu_1(  x_1) \cdots {\rm d}\mu_n( x_n) \leq \widetilde C \prod_{j=1}^n \mu_j(Q_j) .
\end{equation}

\end{itemize}

 \end{proposition}

\smallskip

\begin{proof}
One may define
$$
   S_j(F_{N,\ldots, N})(x_1,\ldots, x_n)=\Bigg( \sum_{i=1}^{j}\sum_{\substack{P_i\in  Ch^{(m_i)}(Q_i)\\ 0\leq m_i\leq N}}   \Bigg| \sum_{k=j+1}^{n}\sum_{\substack{P_k\in  Ch^{(m_k)}(Q_k)\\ 0\leq m_k\leq N}}    { \mathfrak{a}_{P_1\times \cdots\times P_n}  (x_1,\ldots, x_n)  \Bigg|^2   }\Bigg)^{1/2},\  \  \text{for}\  \   j=1,\ldots, n.
$$
In particular, $S_n(F_{N,\ldots, N})=S(F_{N,\ldots, N})$. 

\smallskip

(a).  Due to \textrm   {(P1)}, \textrm   {(P2)} and \textrm   {(P3)} assumed, the arguments for  Propositions \ref{coro:Lp-CH} and   \ref{prop:double-Lp-CH-intro} are also applicable, hence  there exists a constant $C>0$ such that 
\begin{equation}\label{eqn:multi-squareLp-1-old}
	    \left \|F_{N,\ldots, N}(\cdot, x_2,\ldots, x_n)\right \|_{L^p(Q_1 ,\, {\rm d}\mu_1)}
	    \leq   C\sqrt{p} \left \|S_1(F_{N,\ldots, N})(\cdot, x_2,\ldots, x_n)\right \|_{L^p( Q_1,\,{\rm d}\mu_1)}
\end{equation}
hold for any fixed $x_k\in Q_k$ with $k\neq 1$ and $p\geq 2$,
 and 
\begin{align}\label{eqn:multi-squareLp-j-old}
   & \left \|S_j(F_{N,\ldots, N})(x_1,\ldots, x_{j-1},\, \cdot,\,  x_{j+1},\ldots, x_n)\right \|_{L^p(Q_{j+1}, \,{\rm d}\mu_{j+1 })}\nonumber\\
   \leq C &  \sqrt{p} \left \|S_{j+1}(F_{N,\ldots, N})(x_1,\ldots, x_{j-1},\, \cdot,\,  x_{j+1},\ldots, x_n)\right \|_{L^p(Q_{j+1}  ,\, {\rm d}\mu_{j+1} )},\quad 1\leq j\leq n-1
\end{align}
hold for any fixed $x_k\in Q_k$ with $k\neq j+1$ and $p\geq 2$.

Therefore,
\begin{align*}
	  \left\|F_{N,\ldots, N}\right\|_{L^p(Q_1\times\cdots\times  Q_n,\,  {\rm d}\mu_1\cdots {\rm d}\mu_n)} &  \leq C\sqrt{p}  \left\|S_1(F_{N,\ldots, N})\right\|_{L^p(Q_1\times\cdots\times  Q_n,\,  {\rm d}\mu_1\cdots {\rm d}\mu_n)}\\
	   & \leq C^2 p \left\|S_2(F_{N,\ldots, N})\right\|_{L^p(Q_1\times\cdots\times  Q_n,\,  {\rm d}\mu_1\cdots {\rm d}\mu_n)}\\
	   &\leq \cdots\cdots\\
	   &\leq C^n p^{n/2} \left\|S(F_{N,\ldots,N})\right\|_{L^p(Q_1\times\cdots\times  Q_n,   \,  {\rm d}\mu_1\cdots {\rm d}\mu_n)  }.
\end{align*}
 That is, {\color{blue}\eqref{eqn:multi-Lp-CH-intro}} holds.

\medskip

(b). By combining \eqref{eqn:multi-squareLp-1-old} and \eqref{eqn:multi-squareLp-j-old}  again, we have that for $p\geq2$ and for any fixed $x_n\in Q_n$,
\begin{align*}
	  \left\|F_{N,\ldots, N}(\cdots, x_n)\right\|_{L^p(Q_1\times\cdots\times  Q_{n-1},\,  {\rm d}\mu_1\cdots {\rm d}\mu_{n-1})} 
	 \leq    C^{n-1}  p^{(n-1)/2} \left\|S_{n-1}(F_{N,\ldots,N})(\cdots, x_n)\right\|_{L^p(Q_1\times\cdots\times  Q_{n-1},   \,  {\rm d}\mu_1\cdots {\rm d}\mu_{n-1})}. 
\end{align*}
This, implies that for any $\alpha>0$ and for any  $p=nk$ with $k\in \mathbb N_+$, 
\begin{align*}
	&\left\|(\alpha F_{N,\ldots, N}  (\cdots, x_n))^{2/n}\right\|_{L^{nk}(Q_1\times\cdots\times  Q_{n-1},\,  {\rm d}\mu_1\cdots {\rm d}\mu_{n-1})} ^{nk}\\
	 =& \left\|\alpha F_{N,\ldots, N} (\cdots, x_n)\right\|_{L^{2k}(Q_1\times\cdots\times  Q_{n-1},\,  {\rm d}\mu_1\cdots {\rm d}\mu_{n-1})} ^{2k} \\
	\leq & C^{2k (n-1)}  (2k)^{k (n-1)}  \alpha^{2k  }  \left\| S_{n-1}(F_{N,\ldots, N}) (\cdots, x_n)\right\|_{L^{2k}(Q_1\times\cdots\times  Q_{n-1},\,  {\rm d}\mu_1\cdots {\rm d}\mu_{n-1})} ^{2k}. 
\end{align*}

Let $\mu(Q_1\times \cdots \times Q_{n-1})=\mu_1(Q_1)\times \cdots \times \mu_{n-1}(Q_{n-1})$.
Similar to the proof of  \eqref{eqn:newCH-dyadic} and using 
$\displaystyle
    \frac{e^k}{k^k}  \leq c \frac{1}{(k-1)!} 
$
 for $k\geq 2$ again, we have for any fixed $x_n\in Q_n$ that 
\begin{align*}
	  &	\int_{Q_1\times \cdots\times Q_{n-1}} \exp\left(  \left| \alpha  F_{N,\ldots, N}(x_1,\ldots, x_n) \right|^{2/n} \right) {\rm d}\mu_1(x_1)\cdots {\rm d}\mu_{n-1}(x_{n-1})\\
	\leq & e^2 \mu(Q_1\times \cdots \times Q_{n-1}) \\
	&\qquad+\sum_{k=2}^\infty e^{k+1} \mu\left (\left \{(x_1,\ldots, x_{n-1})\in Q_1\times \cdots\times Q_{n-1}:\,   k\leq \left| \alpha  F_{N,\ldots, N}(x_1,\ldots, x_n) \right|^{2/n} <k+1\right \}\right )\\
		\leq & e^2 \mu(Q_1\times \cdots \times Q_{n-1}) + \sum_{k=2}^\infty \frac{e^{k+1}}{k^{nk}}    \left\|(\alpha F_{N,\ldots, N}  (\cdots, x_n))^{2/n}\right\|_{L^{nk}(Q_1\times\cdots\times  Q_{n-1},\,  {\rm d}\mu_1\cdots {\rm d}\mu_{n-1})} ^{nk}\\
	\leq & e^2 \mu(Q_1\times \cdots \times Q_{n-1}) + \sum_{k=2}^\infty \frac{e^{k+1}}{k^{nk}}   
	C^{2k (n-1)}  (2k)^{k (n-1)}  \alpha^{2k  }  \left\| S_{n-1}(F_{N,\ldots, N}) (\cdots, x_n)\right\|_{L^{2k}(Q_1\times\cdots\times  Q_n,\,  {\rm d}\mu_1\cdots {\rm d}\mu_{n-1})} ^{2k} 
	 \\
	\leq & e^2 \mu(Q_1\times \cdots \times Q_{n-1})  +e\cdot  \int_Q \sum_{k=2}^\infty \frac{e^k}{k^k}   \left (2^{n-1}C^{2(n-1)} \alpha^2  S_{n-1}^2(F_{N,\ldots, N})(x_1,\ldots, x_n)  \right)^k\,
	{\rm d}\mu_1(x_1) \cdots {\rm d}\mu_{n-1}(x_{n-1})  \\
	\leq & e^2 \mu(Q_1\times \cdots \times Q_{n-1})  +ec  \int_Q   2^{n-1}C^{2(n-1)} \alpha^2  S_{n-1}^2(F_{N,\ldots, N})(x_1,\ldots, x_n)    \\
	&\qquad \qquad  \qquad \qquad \qquad \qquad  \times    \exp\left(2^{n-1}C^{2(n-1)} \alpha^2  S_{n-1}^2(F_{N,\ldots, N})(x_1,\ldots, x_n)   \right)  \,
	{\rm d}\mu_1(x_1) \cdots {\rm d}\mu_{n-1}(x_{n-1})\\
	\leq& (e^2+ec)  \int_Q \exp\left(2^n C^{2(n-1)} \alpha^2  S_{n-1}^2(F_{N,\ldots, N})(x_1,\ldots, x_n)   \right)  \,
	{\rm d}\mu_1(x_1) \cdots {\rm d}\mu_{n-1}(x_{n-1}).
\end{align*}
We obtain the Hoeffding inequality \eqref{eqn:one-shift-bi-multi-1}.

The argument of another multi-parameter Hoeffding inequality \eqref{eqn:one-shift-bi-multi-2} is similar to that of  Proposition \ref{prop:one-shift-bi-I}, and we skip it.

\medskip

(c).   Similar to the proof of  \eqref{eqn:newCWW-one-parameter}, we can apply \eqref{eqn:multi-squareLp-j-old} with $j=n-1$ to  obtain an analogue of \eqref{eqn:newCWW-S1-S-0} that there exist positive constants $c'$ and $C'$ such that  for any fixed $x_j\in Q_j$ with $1\leq j\leq n-1$,
\begin{equation}\label{eqn:multi-newCH-dyadic-2}
 \int_{  Q_n} \exp\Bigg(c'\frac{   | S_{n-1}(F_{N,\ldots, N})( x_1,\ldots, x_n)|^2   }{\|S(F_{N,\ldots ,N})\|_{L^\infty(Q_1\times \cdots\times Q_n)}^2}\Bigg)  \,{\rm d}\mu_n(  x_n)  \leq   C'   \mu_n(Q_n) .
\end{equation}
This, combined with \eqref{eqn:one-shift-bi-multi-1} and taking $\alpha=\sqrt{\frac{c'}{c}}\|S(F_{N,\ldots, N})\|_{L^\infty(Q_1\times \cdots\times Q_n)}^{-1}$ therein, implies that
\begin{align*}
	& \int_{Q_1\times\cdots \times  Q_n} \exp\Bigg( \Big(\frac{c'}{c}\Big  )^{1/n }  \frac{   |   F_{N,\ldots, N}( x_1,\ldots, x_n)|^{2/n}   }{\|S(F_{N,\ldots ,N})\|_{L^\infty(Q_1\times \cdots\times Q_n)}^{2/n}}\Bigg)  \,{\rm d}\mu_1(  x_1) \cdots {\rm d}\mu_n( x_n) \\
	&\leq C \int_{Q_1\times\cdots \times  Q_n} \exp\Bigg(c'\frac{   | S_{n-1}(F_{N,\ldots, N})( x_1,\ldots, x_n)|^2   }{\|S(F_{N,\ldots ,N})\|_{L^\infty(Q_1\times \cdots\times Q_n)}^2}\Bigg)   \,{\rm d}\mu_1(  x_1) \cdots {\rm d}\mu_n( x_n) \\
	&\leq C C'  \prod_{j=1}^n \mu_j(Q_j) .
\end{align*}
Therefore \eqref{eqn:multi-newCH-dyadic} holds and we complete the proof of Proposition \ref{prop:multi-newCH-dyadic}.
\end{proof}
 
\smallskip

\begin{remark}{
The one-parameter  Chang--Wilson--Wolff type result \eqref{eqn:newCWW-one-parameter} can be deduced by a rough Hoeffding inequality \eqref{eqn:newCH-dyadic-2}. However, 
to show the bi-parameter and multi-parameter versions, i.e., \eqref{eqn:double-newCH-dyadic} and \eqref{eqn:multi-newCH-dyadic}, respectively, we have to combine a more subtle Hoeffding inequality associated to $S_{n-1}$ (see \eqref{eqn:newCH-dyadic-coro1} for $n=2$ and \eqref{eqn:one-shift-bi-multi-1} for $n\geq 3$) with an intermediate Chang--Wilson--Wolff type result for the function $S_{n-1}$ and ``its square function'' $S_n$ (see \eqref{eqn:newCWW-S1-S-0} for $n=2$ and \eqref{eqn:multi-newCH-dyadic-2} for $n\geq 3$). The role of our Hoeffding's inequalities \eqref{eqn:newCH-dyadic-coro1} and \eqref{eqn:one-shift-bi-multi-1} without using the $L^\infty$ structure of atomic functions is fundamental and crucial in multi-parameter settings.
}
\end{remark}

\medskip

\section{Application: bounded multi-parameter square functions on the Shilov boundary}\label{sec:5}
\setcounter{equation}{0}

Let  $\widetilde {\mathscr M} =\mathscr M_1\times  \cdots\times \mathscr M_n$ be the Shilov boundary studied by Nagel and Stein \cite{NS04}, with each $\mathscr M_j$
the polynomial model domain as recalled in Section \ref{sec:1.3}. We first recall the fundamental properties of the heat kernel of the sub-Laplacian on
each $\mathscr M_j$. For the sake of simplicity, we drop the subscript $j$ in Section \ref{sec:5.1} below.

\subsection{Heat kernels for the sub-Laplacian on $\mathscr M$}\label{sec:5.1}
We note that (\cite{NS04}) $\mathscr M$ can be identified with 
$\mathbb C\times\mathbb R=\{ (z,\mathfrak t): z\in\mathbb C, \mathfrak t\in\mathbb R\}$. The basic (0,1) Levi vector field is then 
$\bar Z= {\partial\over\partial \bar z}-i {\partial \mathcal P\over\partial \bar z}{\partial\over\partial \mathfrak t}$, where $i^2=-1$, and we write $\bar Z = X_1+iX_2$. The real vector fields $\{X_1, X_2\}$ and their commutators of order up to $m$  span the tangent space at each point.

To construct the Littlewood--Paley square function, Nagel and Stein in \cite{NS04} 
considered the sub-Laplacian $\mathcal {L}$ on $\mathscr M$ in self-adjoint
form, given by
\begin{eqnarray*}
\mathcal {L}=-\sum_{j=1}^2 X_j^2.
\end{eqnarray*}
The kernel  $H(s,{\bf x}, {\bf y})$ of the heat semigroup $H_s=e^{-s\mathcal{L}}$ has some nice properties (see Proposition 2.3.1
in \cite{NS04} and Theorem 2.3.1 in \cite{NS01a}). We restate some of them. Below, $\partial_X^L$ denotes $X^I = Y_1^{i_1} Y_2^{i_2}\cdots  Y_\ell^{i_\ell}$ (with $ Y_1,\dots, Y_\ell\in\{X_{1},X_{2}\}$) acting on the variable ${\bf x}$ such that $d(I)=i_1+\cdots+i_\ell=L$, and 
$\partial_Y^K$ denotes $X^I = Y_1^{i_1} Y_2^{i_2}\cdots  Y_\ell^{i_\ell}$ (with $ Y_1,\dots, Y_\ell\in\{X_{1},X_{2}\}$) acting on the variable ${\bf y}$ such that $d(I)=i_1+\cdots+i_\ell=K$.
\begin{itemize}
\item[(1)] $H(s,{\bf x},{\bf y})\in C^{\infty}\big([0,\infty)\times \mathscr M\times \mathscr M
\backslash \{s=0\ {\rm and}\ {\bf x}= {\bf y} \}\big)$, $H(s,{\bf x}, {\bf y})=H(s,{\bf y}, {\bf x})$ and 
$H(s,{\bf x}, {\bf y})\geq 0.$

\item[(2)] For every integer $N\geq 0$, 
\begin{align*}
|\partial_s^j\partial_X^L\partial_Y^K H(s,{\bf x}, {\bf y})|
& \lesssim  \frac{\displaystyle 1 }{\displaystyle
(d_{cc}({\bf x}, {\bf y} )+\sqrt{s})^{2j+K+L} }\frac{\displaystyle 1 }{\displaystyle
V({\bf x}, d_{cc}({\bf x},{\bf y}) )+V ({\bf x},\sqrt s) +V ({\bf y}, \sqrt s ) } \bigg(\frac{\displaystyle
\sqrt{s} }{\displaystyle d_{cc}({\bf x}, {\bf y})+\sqrt{s} } \bigg)^{N\over 2},
\end{align*}
where $V({\bf x}, d_{cc}({\bf x},{\bf y}) ) = \mu(B_{cc}({\bf x},d_{cc}({\bf x},{\bf y})))$.
 {The definition of the Carnot--Caratheodory metric $d_{cc}$ and the related nonisotropic ball $B_{cc}$ (as well as the measure $\mu$) were given in Section \ref{sec:1.3}).}

\item[(3)] For all $(s, {\bf x})\in (0,\infty)\times \mathscr M$, $\int_{\mathscr M} H(s,{\bf x}, {\bf y}) \,{\rm d}\mu({\bf y})=1.$

\item[(4)]  For $1\leq p\leq \infty$, $\|H_s(f)\|_{L^p(\mathscr M)}\leq
\|f\|_{L^p(\mathscr M)}$.
\end{itemize}

To introduce the reproducing identity and the Littlewood--Paley
square function, Nagel and Stein define a bounded operator
$Q_s=2s{\displaystyle{\partial H_s}\over\displaystyle \partial s}$,
$s>0$, on $L^2(M)$. Denote by $q_s({\bf x},{\bf y})$ the kernel of $Q_s$. Then
from the estimates of $H(s,{\bf x}, {\bf y})$, we have that
\begin{itemize}
\item[(a)] $q_s({\bf x}, {\bf y})\in C^{\infty}\big(\mathscr  M\times \mathscr M \backslash \{ {\bf x}= {\bf y} \}\big).$

\item[(b)] For every integer $N\geq 0$,
\begin{align*}
    \left|\partial_X^L\partial_Y^K q_s({\bf x}, {\bf y})\right|  \lesssim &
    \frac{\displaystyle 1 }{\displaystyle (d_{cc}({\bf x},{\bf y})+\sqrt{s})^{K+L} }
    \frac{\displaystyle 1 }{\displaystyle V({\bf x}, d_{cc}({\bf x},{\bf y}))+V({\bf x},\sqrt{s})
    +V({\bf y}, \sqrt{s}) } \bigg(\frac{\displaystyle \sqrt{s}}{\displaystyle d_{cc}({\bf x}, {\bf y})+\sqrt{s} } \bigg)^{N\over 2}.
\end{align*}

\item[(c)] $\int_{\mathscr M} q_s({\bf x}, {\bf y})\, {\rm d}\mu({\bf y})=\int_{\mathscr M} q_s({\bf x}, {\bf y})\, {\rm d} \mu({\bf x})=0.$
\end{itemize}

The reproducing identity was established via the operators
$\{Q_s\}_{s>0}$, which play an important role in the Littlewood--Paley
theory and the boundedness of singular integral operators. 

\begin{proposition}[\cite{NS04}]\label{prop-of-reproducing formula by nagel and stein}
Let $Q_s^2=Q_s\cdot Q_s$. For $f\in L^2(\mathscr M)$,
\begin{eqnarray}\label{reproducing formula by nagel and stein}
    \int_0^\infty Q_s^2(f){{\rm d}s\over s} =f,
\end{eqnarray}
where the integral on the left is defined as
$\lim\limits_{\epsilon\rightarrow 0}\int_{\epsilon}^{1/\epsilon}Q_s^2(f){{\rm d}s\over s}$, with the limit
taken in the $L^2(\mathscr M)$ norm.
\end{proposition}

\subsection{Shilov boundary $\widetilde{\mathscr M}=\mathscr M_1\times\cdots\times \mathscr M_n$}\label{sec:5.2}

Suppose $j=1,2\ldots,n$.
Let $\mathcal L_j$ be the sub-Laplacian on $\mathscr M_j$, $H^{[j]}_{t_j}=e^{-t_j\mathcal{L}_j}$ ($t_j\geq0$) and $Q^{[j]}_{t_j}=2t_j{\displaystyle{\partial H^{[j]}_{t_j}}\over\displaystyle \partial t_j}$.  For $0<\gamma_j<\infty$ and $\gamma=(\gamma_1,\ldots,\gamma_n)$,  the area function of $f\in L_{loc}^1(\widetilde{\mathscr M})$ is defined as
$$
    S_\gamma(f)({\bf x}_1,\ldots,{\bf x}_n) =\left( \iint_{\Gamma_\gamma({\bf x}_1,\ldots,{\bf x}_n)} \Big| Q^{[1]}_{t_1^2}\cdots Q^{[n]}_{t_n^2}(f)({\bf y}_1,\ldots,{\bf y}_n) \Big|^2 {{\rm d} \mu_1({\bf y_1}) {\rm d} t_1\over V({\bf x}_1,t_1) t_1} \cdots { {\rm d}  \mu_n({\bf y_n}) {\rm d} t_n\over V({\bf x}_n,t_n) t_n} \right)^{1\over2},
$$
where  $\Gamma_\gamma({\bf x}_1,\ldots,{\bf x}_n) = \Gamma^{[1]}_{\gamma_1}({\bf x}_1) \otimes\cdots\otimes\Gamma^{[n]}_{\gamma_n}({\bf x}_n)$, with 
$ \Gamma^{[j]}_{\gamma_j}({\bf x}_j) = \left\{ ({\bf y}_j,t_j)\in\mathscr  M_j\times \mathbb R_+:\ d_{cc}( {\bf x}_j, {\bf y}_j )< \gamma_j t_j\right\}$.

To ease the burden of notational complexity, we now consider $\widetilde{\mathscr M}=\mathscr M_1\times \mathscr M_2$ below.
The argument for general $\widetilde{\mathscr M}=\mathscr M_1\times\cdots\times \mathscr M_n$ with $n\geq 3$ is similar.

For a fixed $\rho>0$, let $\phi\in C_0^\infty (\mathbb R)$ be real-valued and even, ${\rm supp\, } \phi \subset (-\rho, \rho)$, and set $\psi(x)=x^2\widehat{\phi}(x)$, $x\in \mathbb R$. Then the kernel $K_{\psi(t_j \sqrt{\mathcal L_j})}$  of $\psi(t_j \sqrt{\mathcal L_j})$, $j=1,\,2$, satisfies
\begin{itemize}
\item[(i)] $\displaystyle 
   {\rm supp\,} K_{\psi(t_j \sqrt{\mathcal L_j})} \subset \left\{ ({\bf x}_j, {\bf y}_j)\in \mathscr M_j\times \mathscr M_j:\,   d_{cc}({\bf x}_j, {\bf y}_j)\leq \rho t_j\right\}$.

\item[(ii)]  for all ${\bf x}_j, \, {\bf x}_j^{\prime},\, {\bf y}_j \in \mathscr M_j$ and all $t_j>0$, 
\begin{equation}\label{eqn:psi-scale}  
    \Big | K_{\psi(t_j \sqrt{\mathcal L_j})} ({\bf x}_j, {\bf y}_j)\Big|\leq \frac{C}{V({\bf x}_j, t_j)},   
\end{equation}
and there exists $\theta\in (0,1]$,
\begin{equation}\label{eqn:Holder}
    \left|K_{\psi(t_j \sqrt{\mathcal L_j})}({\bf x}_j , {\bf y}_j)-K_{\psi_(t_j \sqrt{\mathcal L_j})} ({\bf x}_j^{\prime}, {\bf y}_j)\right| \leq C\min \left\{\left(\frac{d_{cc}({\bf x}_j, {\bf x}_j^{\prime})}{t_j}\right)^{\theta}, 1\right\} \frac{1}{V({\bf y}_j, t_j)}.
\end{equation}
 
\item [(iii)] the kernel $K_{\psi (s_j \sqrt{\mathcal L_j}) \psi (t_j \sqrt{\mathcal L_j})}({\bf x}_j, {\bf y}_j)$ of $\psi (s_j \sqrt{\mathcal L_j}) \psi (t_j \sqrt{\mathcal L_j})$ satisfies
\begin{equation}\label{eqn:al-orth}
    \left|K_{\psi (s_j \sqrt{\mathcal L_j}) \psi (t_j \sqrt{\mathcal L_j})}({\bf x}_j, {\bf y}_j)\right| \leq \min \left(\frac{s_j}{t_j}, \frac{t_j}{s_j}\right) \frac{C}{V({\bf x}_j, \max (s_j, t_j))}, \quad \forall s_j, t_j>0, \  \   \forall {\bf x}_j, {\bf  y}_j \in \mathscr M_j,
\end{equation}
\end{itemize}

\smallskip

For $j=1,2$, let $\mathscr D^{[j]}$ be a system of dyadic cubes in the homogeneous space $(\mathscr M_j,d_{cc},\mu_j)$. Recall that the sidelength $\ell(I_j)=\delta^k$ for $I\in \mathscr D_k^{[j]}$,  we define
$$
    I_+=I\times \big[\rho^{-1}\delta^{k+1}, \rho^{-1}\delta^k\big],
$$
where $\rho>0$ is the parameter such that ${\rm supp\, } \phi \subset (-\rho, \rho)$ as aforementioned.
For any dyadic rectangle $R=I\times J$ with $I\in \mathscr D^{[1]}$ and $J\in \mathscr D^{[2]}$, denote
  $ R_+=I_+\times J_+$.

By the $H^\infty$-functional calculus, for any $f\in L^2(\widetilde {\mathscr M})$, we have the  Calder\'on reproducing formula 
\begin{equation}\label{eqn:Calderon-reproducing}
    f({\bf x}_1,{\bf x}_2)=c_\psi \int_0^\infty \int_0^\infty \psi\big(t_1\sqrt{\mathcal L_1}\big) \psi\big(t_2\sqrt{\mathcal L_2}\big) Q_{t_1^2}^{[1]} Q_{t_2^2}^{[2]} (f) ({\bf x}_1, {\bf x}_2) \frac{{\rm d} t_1}{t_1}\frac{ {\rm d} t_2}{t_2},
\end{equation} 
where the limit converges in $L^2(\widetilde {\mathscr M})$.  This gives the atomic decomposition
\begin{align*}
    &f({\bf x}_1, {\bf x}_2)\\
    &=\sum_{R=I\times J\in \mathscr D^{[1]}\times \mathscr D^{[2]}}  c_{\psi} \iint_{R_+} K_{\psi\big(t_1\sqrt{\mathcal L_1}\big) }({\bf x}_1, {\bf y}_1) K_{\psi\big(t_2\sqrt{\mathcal L_2}\big)}({\bf x}_2, {\bf y}_2) Q_{t_1^2}^{[1]} Q_{t_2^2}^{[2]} (f) ({\bf y}_1, {\bf y}_2) \frac{   {\rm d} \mu_1( {\bf y}_1) {\rm d} t_1}{t_1}\frac{ {\rm d}\mu_2( {\bf y}_2)  {\rm d}t_2}{t_2}\\
    &=:\sum_{R=I\times J\in \mathscr D^{[1]}\times \mathscr D^{[2]}}    a_R(f)({\bf x}_1, {\bf x}_2).
\end{align*}
Set
\begin{equation}\label{eqn:lambdaR}
  \lambda_R(f):=\left( \iint_{R_+}  \left| Q_{t_1^2}^{[1]} Q_{t_2^2}^{[2]} (f) ({\bf y}_1, {\bf y}_2) \right|^2 \frac{{\rm d} \mu_1({\bf y}_1)  {\rm   d}t_1}{t_1}\frac{{\rm d}\mu_2( {\bf y}_2) {\rm d}t_2}{t_2}\right)^{1/2}.
\end{equation}
From properties of the kernels of $\psi\big(t_j\sqrt{\mathcal L_j}\big)$, there exists an absolute constant $C>0$ such that 
\begin{itemize}
\item[(a)]  $\supp {a}_R(f)\subseteq 3R:= 3I\times 3J$. 
\smallskip

\item  [(b)] $\|a_R(f)\|_{L^\infty}\leq C \lambda_R(f)\mu(R)^{-1/2}$, where $\mu(R)=\mu_1(I)\times \mu_2(J)$.

\smallskip

\item[(c)] for any ${\bf x}_1, \, {\bf x}'_1\in \mathscr M_1 $ and ${\bf x}_2,\, {\bf x}'_2\in \mathscr M_2$, 
$$
   \left| a_R(f)({\bf x}_1, {\bf x}_2)-a_R(f)({\bf x}'_1, {\bf x}_2)\right|\leq C  \min \left\{ \bigg(\frac{d_1({\bf x}_1, {\bf x}'_1)}{\ell(I)}\bigg)^{\theta}, 1\right\}   \lambda_R(f) \,\mu(R)^{-1/2}
$$
and
$$
   \left| a_R(f)({\bf x}_1, {\bf x}_2)-a_R(f)({\bf x}_1, {\bf x}'_2)\right|\leq C  \min\left\{  \bigg(\frac{d_2({\bf x}_2, {\bf x}'_2)}{\ell(J)}\bigg)^{\theta} , 1\right\} \lambda_R(f)   \,\mu(R)^{-1/2}.
$$  
\end{itemize}

For any open  set $\Omega\subset  \widetilde{\mathscr M}$, define
$$
    F_{\Omega}({\bf x}_1,{\bf x}_2)=\sum_{\substack{R\subseteq \Omega\\ R:\, \text{dyadic}}}  a_R(f)({\bf x}_1,{\bf x}_2).
$$

Our main result in this section is the following two-parameter version of the exponential square integrability by Chang--Wilson--Wolff \cite{CWW}, as well as the extension  of the second author's work \cite[Theorem 2]{Pi86} on bounded double square functions. 

\begin{theorem}\label{thm:d-SWW}
Suppose $f\in L^1_{loc}(\widetilde {\mathscr M})$ and $S_\gamma(f)\in L^\infty(\widetilde{\mathscr M})$. Then there exist constants $c_1$ depending of $\gamma$, and  $c_2$ independent of $\gamma$ such that for all open set $\Omega\in\widetilde{\mathscr  M}$ with finite measure 
\begin{align}\label{eqn:d-SWW}
    \int_{\Omega} \exp\bigg(c_1  \frac {|F_{\Omega}({\bf x}_1,{\bf x}_2)|}{ \|S_\gamma(f)\|_{ L^\infty(\widetilde{\mathscr  M})} }\bigg) \,{\rm d} \mu_1({\bf x}_1)  {\rm d}\mu_2({\bf x}_2) \leq c_2  \mu(\Omega).
\end{align} 
\end{theorem}

\begin{proof}
The argument is based on  Theorem \ref{thm:double-newCH-dyadic} associated to double martingale differences, and properties of the atomic functions $a_R(f)$.

Let $\left\{\mathcal D^{[j],b}:\, b=1,2,\ldots, K_j\right\}$ be adjacent dyadic systems on $\mathscr M_j$ for $j=1,2$. 
Then it follows from the doubling volume property that there exists a constant $C>0$ such that for any dyadic rectangle 
$R=I\times J\in \mathscr D^{[1]}\times \mathscr D^{[2]}$, there exists $E(R)\in   \mathcal D^{[1],b_1}\times \mathcal D^{[2],b_2}$ satisfying
\begin{equation}\label{eqn:injective}
	   R\subseteq \supp {a}_R(f) \subseteq E(R)\quad \text{and}\quad \mu(E(R))\leq C \mu(R).
\end{equation}
 
Moreover, note that 
$$
  \widetilde{\Omega} := \bigcup_{\substack{R\subseteq \Omega\\ R:\, \text{dyadic}}}  E(R)\supseteq \Omega\quad \text{and}\quad \mu(\widetilde{\Omega} )\lesssim \mu(\Omega).
$$
Hence one may prove \eqref{eqn:d-SWW} by replacing the domain  $\Omega$ of integration therein by $\widetilde{ \Omega}$. Alternatively, one may assume in the sequel that 
$$
    F_{\Omega}({\bf x}_1,{\bf x}_2)=\sum_{\substack{E(R)\subseteq \Omega\\ R:\, \text{dyadic}}}  a_R(f)({\bf x}_1,{\bf x}_2)
$$
and continue to show \eqref{eqn:d-SWW} with this notation.

Then we split $F_{\Omega}$ as
$$
      F_\Omega({\bf x}_1, {\bf x}_2)=\sum_{\bf b} F_\Omega^{\bf b} ({\bf x}_1, {\bf x}_2),
$$
where ${\bf b}=(b_1, b_2)\in \{1,\ldots, K_1\}\times \{1,\ldots, K_2\}$, and 
$$
    F_\Omega^{\bf b} ({\bf x}_1, {\bf x}_2)=\sum_{\substack{  E(R)\in \mathcal D^{[1],b_1} \times  \mathcal D^{[2], b_2} \\ E(R)\subseteq \Omega}}  a_R(f)({\bf x}_1, {\bf x}_2).
$$

It suffices to prove that for each ${\bf  b}=(b_1, b_2)$,
$$
    \int_{\Omega} \exp\bigg(c_1\frac {|F_{\Omega}^{\bf b}    ({\bf x}_1,{\bf x}_2)|}{ \|S_\gamma(f)\|_{ L^\infty(\widetilde{\mathscr  M})}}\bigg) \,{\rm d} \mu_1({\bf x}_1) {\rm d} \mu_2({\bf x}_2)  \leq c_2\mu(\Omega).
$$
Moreover, for $M\geq 1$, let 
$$
    F_{\Omega,M}^{\bf b}:= \sum_{\substack{ E(R)\in \bigcup_{i=-M}^M  \mathcal D_i^{[1],b_1} \times  \bigcup_{j=-M}^M \mathcal D_j^{[2], b_2}  \\  E(R)\subseteq \Omega  }}  a_R(f)({\bf x}_1, {\bf x}_2),
$$
and $ F_{\Omega,M}^{\bf b}$ is a bounded function  on $\Omega$.
By Fatou's lemma, it suffices to show that   
\begin{align}\label{eqn:d-SWW-new}
\int_{\Omega} \exp\bigg( c_1 \frac{|F_{\Omega,M}^{\bf b}    ({\bf x}_1,{\bf x}_2)|}{ \|S_\gamma(f)\|_{ L^\infty(\widetilde{\mathscr  M})}}\bigg) \,{\rm d} \mu_1({\bf x}_1) {\rm d} \mu_2({\bf x}_2) \leq c_2  \mu(\Omega)
\end{align} 
for any $M\geq1$, where $c_1$ and $c_2$ are  independent of $M$.

To this end, for any fixed $M\geq1$,  we will express $F_{\Omega,M}^{\bf b}  $ in terms of   Haar functions 
$$\left\{h_{u_1}^{I'} h_{u_2}^{J'}:\, I'\times J'\in \mathcal D^{[1],b_1}\times \mathcal D^{[2],b_2} ,\, u_1=1,\ldots, \mathtt M_{I'} -1,\    u_2=1,\ldots, \mathtt M_{J'} -1 \right\}.
$$ 

Recall for any $ R'=I'\times J'$, the double martingale difference of $F_{\Omega,M}^{\bf b}$ associated to $R'$ is given by 
\begin{align*}
	    &\Delta_{R'}(F_{\Omega,M}^{\bf b})\\
	    =&\sum_{I''\in Ch(I') }  \sum_{ J''\in Ch(J')}  \Big [ (F_{\Omega,M}^{\bf b})_{  I''\times   J''} \mathsf 1_{  I''\times    J''} -  (F_{\Omega,M}^{\bf b})_{   I'' \times J'}\mathsf 1_{  I''\times J'} -(F_{\Omega,M}^{\bf b})_{I'\times   J''}\mathsf 1_{I'\times   J''} +(F_{\Omega,M}^{\bf b})_{I'\times J'}\mathsf 1_{I'\times J'}  \Big ]\\
	    =&\sum_{u_1=1}^{\mathtt M_{I'}-1}\sum_{u_2=1}^{\mathtt M_{J'}-1}  \left\langle F_{\Omega,M}^{\bf b}, h_{u_1}^{I'} h_{u_2}^{J'}\right \rangle   h_{u_1}^{I'} h_{u_2}^{J'}.
\end{align*}
In particular, it follows from the cancellation that
\begin{equation}\label{eqn:cancellation-aR}
    \int_{\mathscr M_1}    a_R(f)  ({\bf x}_1,{\bf x}_2) {\rm d}\mu_1  ({\bf x}_1)=\int_{\mathscr M_2}    a_R(f)  ({\bf x}_1,{\bf x}_2){\rm d}\mu_2  ({\bf x}_2)=0
\end{equation}
and Theorem \ref{thm:Haar-base} that if $I'\in \mathcal D_i^{[1], b_1}$ for some $i<-M$ or $J'\in \mathcal D_j^{[1], b_2}$ for some $j<-M$,  we have
$$
       \Delta_{R'}(F_{\Omega,M}^{\bf b})\equiv 0.
$$
Therefore,
\begin{align*}
    F_{\Omega,M}^{\bf b}&=\sum_{   R'=I'\times J'\in \mathcal D^{[1],b_1}\times  \mathcal D^{[2],b_2}  }
	   \Delta_{R'}(F_{\Omega,M}^{\bf b})   
	=\sum_{   R'=I'\times J'\in \bigcup_{n\geq -M} \mathcal D_n^{[1],b_1}\times \bigcup_{m\geq -M}\mathcal D_m^{[2],b_2} }
  \Delta_{R'}(F_{\Omega,M}^{\bf b}).
\end{align*}

Furthermore,  for any $R'=I'\times J'$ with $R' \cap \Omega^c\neq\emptyset$, 
$$ 
      \Delta_{R'}(F_{\Omega,M}^{\bf b})=\sum_{\substack{ E(R)\in \bigcup_{i=-M}^M  \mathcal D_i^{[1],b_1} \times  \bigcup_{j=-M}^M \mathcal D_j^{[2], b_2}  \\ E( R)\subseteq \Omega  }}  \Delta_{R'} (a_R(f)).
$$
Thus for each $R$ with $E(R):=\widetilde I\times \widetilde J $ such that $E(R)\cap R'\neq\emptyset$, we have $I'\supsetneqq \widetilde I$ or $J'\supsetneqq \widetilde J$.
Hence,  by \eqref{eqn:cancellation-aR} and Theorem \ref{thm:Haar-base} again, we have 
$\Delta_{R'} (a_R(f))\equiv 0$. This, implies that 
\begin{align*}
	F_{\Omega,M}^{\bf b}
	&=\sum_{ \substack{  R'=I'\times J'\in \bigcup_{n\geq -M} \mathcal D_n^{[1],b_1}\times \bigcup_{m\geq -M}\mathcal D_m^{[2],b_2} \\   R'\subseteq \Omega  } }
  \Delta_{R'}(F_{\Omega,M}^{\bf b})\\
  &= \lim_{N\to +\infty}\sum_{ \substack{  R'=I'\times J'\in \bigcup_{n= -M}^N \mathcal D_n^{[1],b_1}\times \bigcup_{m= -M}^N \mathcal D_m^{[2],b_2} \\   R'\subseteq \Omega  } }
  \Delta_{R'}(F_{\Omega,M}^{\bf b})\\
  &= :  \lim_{N\to +\infty}   \mathbf F_{\Omega,M,N}^{\bf b}.
\end{align*}
In a similar vein, we turn to proving that 
\begin{align}\label{eqn:d-SWW-new-Haar}
    \int_{\Omega} \exp\bigg( c_1\frac{|\mathbf F_{\Omega,M,N}^{\bf b}    ({\bf x}_1,{\bf x}_2)|}{\|S_\gamma(f)\|_{ L^\infty(\widetilde{\mathscr  M})}}\bigg) \,{\rm d} \mu_1({\bf x}_1) {\rm d} \mu_2({\bf x}_2) \leq c_2 \mu(\Omega)
\end{align} 
for any $M,\, N\in \mathbb N$, where $c_1$ and $c_2$ are  independent of $M$ and $N$.

To this end, now fix $M$ and $N$. Note that  for each $R'$,  the difference $ \Delta_{R'}(F_{\Omega,M}^{\bf b})$ satisfies \textsf{(P1)}, \textsf{(P2)} and the cancellation 
$$
    \int_{\mathscr M_1}    \Delta_{R'} (F_{\Omega,M}^{\bf b})  ({\bf x}_1,{\bf x}_2)   {\rm d}\mu_1 ({\bf x}_1)=\int_{\mathscr M_2}    \Delta_{R'}(F_{\Omega,M}^{\bf b})   ({\bf x}_1,{\bf x}_2)   {\rm d}\mu_2   ({\bf x}_2)=0.
$$ 
As a consequence, $\left\{ \Delta_{R'}(F_{\Omega,M}^{\bf b})\right \}_{R'}$ satisfies \textsf{(P3)} with the constant $C_2\equiv 1$. Hence, by Theorem \ref{thm:double-newCH-dyadic} and Remark \ref{rem:global-omega}, there exist positive constants $c'$ and $C'$ independent of $M,\,N$ and $\Omega$ such that
\begin{align}\label{eqn:double-Haar-Hoeffding}
	\int_{\Omega} \exp\Bigg(c' \frac{   | \mathbf F_{\Omega,M,N}^{\bf b} ({\bf x}_1,{\bf x}_2)|   }{\|S_{\Delta}(\mathbf F_{\Omega,M,N}^{\bf b})\|_{L^\infty(\widetilde{\mathscr  M})}}\Bigg)  \,{\rm d}\mu_1(  {\bf x}_1) \, {\rm d}\mu_2({\bf  x}_2) \leq C' \mu(\Omega),
\end{align}
where $S_{\Delta}(\mathbf F_{\Omega,M,N}^{\bf b})$ is our new double square function associated to martingale differences:
$$
  S^2_{\Delta}(\mathbf F_{\Omega,M,N}^{\bf b})({\bf x}_1,{\bf x}_2)=  \sum_{ \substack{  R'=I'\times J'\in \bigcup_{n= -M}^N \mathcal D_n^{[1],b_1}\times \bigcup_{m= -M}^N \mathcal D_m^{[2],b_2} \\   R'\subseteq \Omega  } }   \left(\Delta_{R'}(F_{\Omega,M}^{\bf b})\right)^2({\bf x}_1,{\bf x}_2).
$$
This shows  the aimed \eqref{eqn:d-SWW-new-Haar}  readily, once the following estimate is proven:
  \begin{equation}\label{eqn:compare-S}
	\left\| S_{\Delta}(\mathbf F_{\Omega,M,N}^{\bf b})\right \|_{L^\infty(\widetilde{ \mathscr M})} \leq C \left \|S_\gamma(f)\right \|_{ L^\infty(\widetilde {\mathscr M})}.
\end{equation}

Firstly, we define the discrete square function
\begin{equation}\label{eqn:discrete-S}
  S_{d}(f) := \bigg( \sum_{   R=I\times J\in \mathscr D^{[1]}\times \mathscr D^{[2]}  }   \frac{|\lambda_R(f)|^2}{\mu(R)} {\mathsf{ 1}}_{E(R)}  \bigg)^{1/2},
\end{equation}
where $\lambda_R(f)$ is given in \eqref{eqn:lambdaR}. From the definitions of $R_+$ and the double volume property, we have 
\begin{align*}
    &S_{d}(f)({\bf x}_1,{\bf x}_2)\\
    \approx &\bigg(\sum_{n_1\in \mathbb Z} \sum_{n_2\in \mathbb Z}  \int_{\rho^{-1}\delta^{n_1+1}}  ^{\rho^{-1}\delta^{n_1}}
    \int_{\rho^{-1}\delta^{n_2+1}}  ^{\rho^{-1}\delta^{n_2}}  \int_{\Omega_{n_1,n_2}({\bf x}_1,{\bf x}_2)} \left| Q_{t_1^2}^{[1]} Q_{t_2^2}^{[2]} (f) ({\bf y}_1, {\bf y}_2) \right|^2 \frac{{\rm d} \mu_1({\bf y}_1) {\rm   d}t_1}{V({\bf x}_1, t_1  )   t_1}\frac{{\rm d}  \mu_2({\bf y}_2) {\rm d}t_2}{   V({\bf x}_2, t_2)  t_2}
    \bigg)^{1/2},
\end{align*}
where  $\delta$ is the parameter in Theorem \ref{thm:dyadic}, and 
$$
     \Omega_{n_1,n_2}({\bf x}_1, {\bf x}_2):=\bigcup_{R=I\times J\in \mathscr D^{[1]}_{n_1}\times \mathscr D^{[2]}_{n_2},\,  ({\bf x}_1, {\bf x}_2)\in E(R)}  R.
$$
For any fixed $\gamma=(\gamma_1,\gamma_2)$, one may verify that by combining  the definitions of $R_+$,  $E(R)$ and the double volume property again,  one may choose (and fix) $\rho=\rho\left(\min\{\gamma_1, \gamma_2\}, \delta \right )$  and constant $c>0$ such that 
 $$
     S_{d}(f)({\bf x}_1,{\bf x}_2)\leq c S_{\gamma}(f)({\bf x}_1,{\bf x}_2)\quad \text{a.e.}\     ({\bf x}_1,{\bf x}_2)\in \widetilde{\mathscr  M}.
 $$
 Consequently, the aimed \eqref{eqn:compare-S} holds if 
 \begin{equation}\label{eqn:compare-S2}
	\left\| S_{\Delta}(\mathbf F_{\Omega,M,N}^{\bf b})\right \|_{L^\infty(\widetilde{ \mathscr M})} \leq C \left \|S_d(f)\right \|_{ L^\infty(\widetilde{ \mathscr M})}.
\end{equation}

We begin to prove it. For any $R'=I'\times J'\in \bigcup_{n= -M}^N \mathcal D_n^{[1],b_1}\times \bigcup_{m= -M}^N \mathcal D_m^{[2],b_2}$ with $   R'\subseteq \Omega$, one may assume that $N>M$ without loss of generality. Observe that $\Delta_{R'}(a_R(f))\equiv 0$ whenever $\widetilde I\subsetneqq  I'$ or $\widetilde J\subsetneqq J'$, where  $E(R)=\widetilde I\times \widetilde J\in \mathcal D^{[1], b_1}\times \mathcal D^{[2], b_2}$.  Hence 
\begin{align}\label{eqn:Haar-express-atom}
	 \Delta_{R'}(F_{\Omega,M}^{\bf b}) &=\sum_{\substack{ E(R)\in \bigcup_{i=-M}^M  \mathcal D_i^{[1],b_1} \times  \bigcup_{j=-M}^M \mathcal D_j^{[2], b_2}  \\  E(R)\subseteq \Omega  }} \Delta_{R'}( a_R(f))   \nonumber \\
	 &= \sum_{\substack{ E(R)\in \bigcup_{i=-M}^M  \mathcal D_i^{[1],b_1} \times  \bigcup_{j=-M}^M \mathcal D_j^{[2], b_2}  \\  R'\subseteq E(R)   \subseteq\Omega }}   \sum_{u_1=1}^{\mathtt M_{I'}-1}\sum_{u_2=1}^{\mathtt M_{J'}-1}  \left\langle  a_R(f) , h_{u_1}^{I'} h_{u_2}^{J'}\right \rangle   h_{u_1}^{I'} h_{u_2}^{J'}.    
\end{align}
For any such $R=I\times J$ with $  R'\subseteq E(R)$,
\begin{align*}
	&\left| \left\langle  a_R(f) , h_{u_1}^{I'} h_{u_2}^{J'}\right \rangle \right| \\
	&=    c_{\psi} \Bigg|\iint_{R_+}   
	 \left[\int_{I'} K_{\psi\big(t_1\sqrt{\mathcal L_1}\big) }({\bf x}_1, {\bf y}_1)    h_{u_1}^{I'} ({\bf x}_1)   \,  {\rm d}\mu_1({\bf x}_1)\right]   
	  \left[  \int_{J'}K_{\psi\big(t_2\sqrt{\mathcal L_2}\big)}({\bf x}_2, {\bf y}_2)   h_{u_2}^{J'} ({\bf x}_2)\, {\rm d}\mu_2({\bf x}_2) \right]  
	\\
	&\qquad\qquad    \times Q_{t_1^2}^{[1]} Q_{t_2^2}^{[2]} (f) ({\bf y}_1, {\bf y}_2) \frac{   {\rm d}  \mu_1({\bf y}_1) {\rm d} t_1}{t_1}\frac{ {\rm d} \mu_2({\bf y}_2) {\rm d}t_2}{t_2}\Bigg|.
\end{align*} 
Denote the center of $I'$ (resp. $J'$) by ${\bf c_{I'}}$ (resp. $\bf c_{J'}$), it follows from Theorem \ref{thm:Haar-base} and \eqref{eqn:Holder}  that  
\begin{align*}
	&\left|\int_{I'} K_{\psi\big(t_1\sqrt{\mathcal L_1}\big) }({\bf x}_1, {\bf y}_1)    h_{u_1}^{I'} ({\bf x}_1)   \,  {\rm d}\mu_1({\bf x}_1)\right|\\
	=& \left|\int_{I'} \left(K_{\psi\big(t_1\sqrt{\mathcal L_1}\big) }({\bf x}_1, {\bf y}_1)  -K_{\psi\big(t_1\sqrt{\mathcal L_1}\big) }({\bf c_{I'}}, {\bf y}_1)   \right) h_{u_1}^{I'} ({\bf x}_1)   \,  {\rm d}\mu_1({\bf x}_1)\right|\\
	\leq & C \,\Big(\frac{\ell(I')}{t_1}\Big)^{\theta} \frac{\mu_1(I')^{1/2}}{V({\bf y}_1, t_1)}  .
\end{align*}
Similarly,
$$
   \left|\int_{J'}K_{\psi\big(t_2\sqrt{\mathcal L_2}\big)}({\bf x}_2, {\bf y}_2)   h_{u_2}^{J'} ({\bf x}_2)\, {\rm d}\mu_2({\bf x}_2) \right| \leq C \,\Big(\frac{\ell(J')}{t_2}\Big)^{\theta} \frac{\mu_2(J')^{1/2}}{V({\bf y}_2, t_2)} .
$$
Therefore,
$$
   \left| \left\langle  a_R(f) , h_{u_1}^{I'} h_{u_2}^{J'}\right \rangle \right| \leq C  \,\Big(\frac{\ell(I')}{\ell(I)}\Big)^{\theta}    \Big(\frac{\ell(J')}{\ell(J)}\Big)^{\theta} \lambda_R(f)\mu(R)^{-1/2}  \mu(R')^{1/2}
$$
and so
\begin{align*}
	    \big | \Delta_{R'}( a_R(f)) ({\bf x}_1,{\bf x}_2) \big |&\leq C  \,\Big(\frac{\ell(I')}{\ell(I)}\Big)^{\theta}    \Big(\frac{\ell(J')}{\ell(J)}\Big)^{\theta} \lambda_R(f)\mu(R)^{-1/2}\mathsf 1_{R'}({\bf x}_1,{\bf x}_2) \\
	    &\leq  C  \,\Big(\frac{\ell(I')}{\ell(\widetilde I)}\Big)^{\theta}    \Big(\frac{\ell(J')}{\ell(\widetilde J)}\Big)^{\theta} \lambda_R(f)\mu(R)^{-1/2}\mathsf 1_{R'}({\bf x}_1,{\bf x}_2)  ,
\end{align*}
where $E(R)=\widetilde I\times \widetilde J$. 

Consequently, we have
\begin{align*}
	&\left(\Delta_{R'}(F_{\Omega,M}^{\bf b})\right)^2({\bf x}_1,{\bf x}_2) \\
	\leq\,& C  \Bigg(\sum_{\substack{ E(R)\in \bigcup_{i=-M}^M  \mathcal D_i^{[1],b_1} \times  \bigcup_{j=-M}^M \mathcal D_j^{[2], b_2}  \\   R' \subseteq E(R)\subseteq \Omega}}    \Big(\frac{\ell(I')}{\ell(\widetilde I)}\Big)^{\theta}    \Big(\frac{\ell(J')}{\ell(\widetilde J)}\Big)^{\theta} \lambda_R(f)\mu(R)^{-1/2}  \Bigg)^2 \mathsf 1_{R'}({\bf x}_1,{\bf x}_2) \\
	\leq\,& C  \sum_{\substack{ E(R)\in \bigcup_{i=-M}^M  \mathcal D_i^{[1],b_1} \times  \bigcup_{j=-M}^M \mathcal D_j^{[2], b_2}  \\   R' \subseteq E(R)\subseteq \Omega}}    \Bigg(  \Big(\frac{\ell(I')}{\ell(\widetilde I)}\Big)^{\theta\over 2}    \Big(\frac{\ell(J')}{\ell(\widetilde J)}\Big)^{\theta\over 2} \lambda_R(f)\mu(R)^{-1/2}  \Bigg)^2 \sum_{i,j\geq 0}  \delta^{(i+j)\theta} \mathsf 1_{R'}({\bf x}_1,{\bf x}_2) \\
	\leq\,& C  \sum_{\substack{ E(R)\in \bigcup_{i=-M}^M  \mathcal D_i^{[1],b_1} \times  \bigcup_{j=-M}^M \mathcal D_j^{[2], b_2}  \\   R' \subseteq E(R)\subseteq \Omega}}      \Big(\frac{\ell(I')}{\ell(\widetilde I)}\Big)^{\theta }    \Big(\frac{\ell(J')}{\ell(\widetilde J)}\Big)^{\theta } \frac{|\lambda_R(f)|^2}{\mu(R)}   \mathsf 1_{R'}({\bf x}_1,{\bf x}_2).
\end{align*}
Hence,
\begin{align*}
	& S^2_{\Delta}(\mathbf F_{\Omega,M,N}^{\bf b})({\bf x}_1,{\bf x}_2)\\
	\leq \, & C   \sum_{ \substack{  R'=I'\times J'\in \bigcup_{n= -M}^N \mathcal D_n^{[1],b_1}\times \bigcup_{m= -M}^N \mathcal D_m^{[2],b_2} \\   R'\subseteq \Omega  } } 
	\sum_{\substack{ E(R)\in \bigcup_{i=-M}^M  \mathcal D_i^{[1],b_1} \times  \bigcup_{j=-M}^M \mathcal D_j^{[2], b_2}  \\   R' \subseteq E(R)\subseteq \Omega}}      \Big(\frac{\ell(I')}{\ell(\widetilde I)}\Big)^{\theta }    \Big(\frac{\ell(J')}{\ell(\widetilde J)}\Big)^{\theta } \frac{|\lambda_R(f)|^2}{\mu(R)}   \mathsf 1_{R'}({\bf x}_1,{\bf x}_2) \\
	\leq \, & C \sum_{\substack{ E(R)\in \bigcup_{i=-M}^M  \mathcal D_i^{[1],b_1} \times  \bigcup_{j=-M}^M \mathcal D_j^{[2], b_2}  \\    E(R)\subseteq \Omega}}     \frac{|\lambda_R(f)|^2}{\mu(R)}   \mathsf 1_{E(R)}({\bf x}_1,{\bf x}_2) 
	\sum_{\substack{R'=I'\times J'\in \mathcal D^{[1], b_1}\times \mathcal D^{[2], b_2}\\  R'\subseteq E(R)  :\,  ({\bf x}_1,{\bf x}_2)\in R'  } } \Big(\frac{\ell(I')}{\ell(\widetilde I)}\Big)^{\theta }    \Big(\frac{\ell(J')}{\ell(\widetilde J)}\Big)^{\theta }\\
	\leq \, & C S_d^2(f)({\bf x}_1,{\bf x}_2) ,
\end{align*}
which implies that \eqref{eqn:compare-S2} holds.

We complete the proof of Theorem \ref{thm:d-SWW}.
\end{proof}

\subsection{Proof of Theorem \ref{thm:app}}

The argument for Theorem \ref{thm:d-SWW} relies on Haar functions and the cancellation and continuity of atomic functions $a_R(f)$. We can show Theorem \ref{thm:app}, the $L^p$ lower bound estimates for the square function $S_\gamma(f)$ on  $\widetilde{\mathscr  M}=\mathscr M_1\times \cdots\times \mathscr M_n$ with $n\geq 2$ in a similar manner.

 \smallskip

\begin{proof}[Proof of of Theorem \ref{thm:app}]
Let $\mathscr D^{[j]}$ be a system of dyadic cubes in $\mathscr M_j$, and  $\left\{\mathcal D^{[j], b_j},\    b_j=1,\ldots, K_j \right\}$ be adjacent systems of dyadic cubes in $\mathscr M_j$, $j=1,\ldots, n$.

For each $f\in L^p(\widetilde{\mathscr  M})$ with $p\geq 2$, let 
$$
    \left\{\Delta_{P_1 \times \ldots \times P_n}(f)({\bf x}_1,\ldots, {\bf x}_n):\quad P_j\in \mathscr D^{[j]},\   1\leq j\leq n\right\}
$$
be the martingale differences of $f$ associated to $\mathscr D^{[1]}\times \cdots\times \mathscr D^{[n]}$. For $j=1,\ldots, n$, let
$$
    S_{\Delta, j}^2(f)({\bf x}_1,\ldots, {\bf x}_n)=\sum_{P_1 \in \mathscr D^{[1]}}    \cdots   \sum_{P_j \in \mathscr D^{[j]}}  \bigg|\sum_{P_{j+1} \in \mathscr D^{[j+1]}} \cdots   \sum_{P_n \in \mathscr D^{[n]} }  \Delta_{P_1 \times \ldots \times P_n}(f)({\bf x}_1,\ldots, {\bf x}_n)\bigg|^2.
$$
In particular, $S_{\Delta, n}(f)$ is indeed our multi-parameter square function associated to martingale differences, we denote it by $S_{\Delta}(f)$.

 In combination with Fatou's lemma, similar to Remark \ref{rem:global-omega}, the local estimates occurred in the proof of Theorems \ref{thm:newCH-dyadic} and \ref{thm:double-newCH-dyadic} can be revised to obtain the global analogs, thus,
  $$
      \|f\|_{L^p(\widetilde{\mathscr  M})}\leq C\sqrt{p} \|S_{\Delta, 1}(f)\|_{L^p(\widetilde{\mathscr  M})}
 $$ 
 and 
 $$
  \|S_{\Delta, j}(f)\|_{L^p(\widetilde{\mathscr  M})}\leq C\sqrt{p}  \|S_{\Delta, j+1}(f)\|_{L^p(\widetilde{\mathscr  M})},\quad 1\leq j\leq n-1.
 $$
 Consequently, 
\begin{equation}\label{eqn:Lp-S-Martinagle}
	  \|f\|_{L^p(\widetilde {\mathscr  M})}\leq C  p^{n\over 2}  \|S_{\Delta}(f)\|_{L^p(\widetilde {\mathscr  M})}.
\end{equation}

To show  Theorem \ref{thm:app},  let $\psi$ be as in Section \ref{sec:5.2}. By the Calder\'on reproducing formula again, 
\begin{align*}
	f({\bf x}_1,\ldots, {\bf x}_n)&=c_\psi \int_0^\infty \int_0^\infty \psi\big(t_1\sqrt{\mathcal L_1}\big) \cdots \psi\big(t_n\sqrt{\mathcal L_n}\big) Q_{t_1^2}^{[1]} \cdots Q_{t_n^2}^{[n]} (f) ({\bf x}_1,\ldots , {\bf x}_n) \frac{{\rm d} t_1}{t_1}\cdots \frac{ {\rm d} t_n}{t_n}\\
	&=\sum_{R=P_1\times \cdots\times P_n\in \mathscr D^{[1]}\times \cdots\times \mathscr D^{[n]}}   c_{\psi} \iint_{R_+} K_{\psi\big(t_1\sqrt{\mathcal L_1}\big) }({\bf x}_1, {\bf y}_1) \cdots K_{\psi\big(t_n\sqrt{\mathcal L_n}\big)}({\bf x}_n, {\bf y}_n) \\
	&\qquad \qquad \qquad \qquad \qquad \times Q_{t_1^2}^{[1]}\cdots Q_{t_n^2}^{[n]} (f) ({\bf y}_1, \ldots, {\bf y}_n) \frac{   {\rm d}  \mu_1({\bf y}_1) {\rm d} t_1}{t_1} \cdots \frac{ {\rm d} \mu_n({\bf y}_n) {\rm d}t_n}{t_n}\\
	&=:\sum_{R=P_1\times \cdots\times P_n\in \mathscr D^{[1]}\times \cdots\times \mathscr D^{[n]}}   a_R(f)({\bf x}_1,\ldots, {\bf x}_n),
\end{align*}
where 
$$
     R_+=(P_1)_+\times\cdots \times (P_n)_+\quad \text{and }\quad  (P_j)_+=P_j\times  \left[\rho^{-1}\delta \ell(P_j),\rho^{-1}\ell(P_j) \right] ,\   j=1,\ldots, n.
$$

Again, for each $R$, there exists a rectangle $E(R)\in \mathcal D^{[1], b_1}\times \cdots\times \mathcal D^{[n],b_n}$ for some $b_j\in\{1,\ldots, K_j\}$, $j=1,\ldots, n$, such that 
$$
    \supp a_R(f)\subseteq E(R)\quad \text{and}\quad \mu(E(R))\leq C \mu(R),
$$
where $C$ is a constant independent of $R$.

Now we split
$$
    f({\bf x}_1,\ldots, {\bf x}_n)=\sum_{{\bf b}=(b_1,\ldots, b_n)} f^{\bf b}({\bf x}_1,\ldots, {\bf x}_n),
$$
where 
$$
    f^{\bf b}({\bf x}_1,\ldots, {\bf x}_n)=\sum_{\substack{R=P_1\times \cdots\times P_n\in \mathscr D^{[1]}\times \cdots\times \mathscr D^{[n]} \\ E(R)\in \mathcal D^{[1], b_1}   \times \cdots\times \mathcal D^{[n],b_n}}   }   a_R(f)({\bf x}_1,\ldots, {\bf x}_n).
$$
Moreover,  similar to the argument of Theorem \ref{thm:d-SWW}, there exists $\rho=\rho(\gamma, \delta)$ such that 
 $$
     S_{d}(f)({\bf x}_1,\ldots, {\bf x}_n)\leq c S_{\gamma}(f)({\bf x}_1,\ldots, {\bf x}_n)\quad \text{a.e.}\     ({\bf x}_1,\ldots, {\bf x}_n)\in \widetilde{\mathscr  M},
 $$
 where $ S_d(f) $ is the  discrete square function for $n$ variables; see  \eqref{eqn:discrete-S} for  the case $n=2$.

Note that  one may rewrite $f^{\bf b}$ in terms of martingale differences associated to  $\mathcal D^{[1], b_1}   \times \cdots\times \mathcal D^{[n],b_n}$, that is, 
$f^{\bf b}=\sum_{R'\in \mathcal D^{[1], b_1}   \times \cdots\times \mathcal D^{[n],b_n}} \Delta_{R'}(f^{\bf b})$. Correspondingly, let 
$$
    S^2_{\Delta, \bf b}   (f^{\bf b}  )({\bf x}_1,\ldots, {\bf x}_n)=\sum_{R'\in \mathcal D^{[1], b_1}   \times \cdots\times \mathcal D^{[n],b_n}}   (\Delta_{R'}(f^{\bf b })   )^2 ({\bf x}_1,\ldots, {\bf x}_n).
$$
Similar to \eqref{eqn:Lp-S-Martinagle}, we have 
$$
      \left\|f^{\bf b}\right \|_{L^p(\widetilde{\mathscr  M})}\leq C p^{n\over 2} \left\| S_{\Delta, \bf b}   (f^{\bf b}  )\right \|_{L^p(\widetilde{\mathscr  M})}.
$$
Using the cancellation of $a_R(f)$ and H\"older's continuity of $K_{\psi(t_j\sqrt{\mathcal L_j})}$ again, one may verify that 
$$
      S_{\Delta, \bf b}   (f^{\bf b}  )({\bf x}_1,\ldots, {\bf x}_n)\leq C S_d(f^{\bf b})   ({\bf x}_1,\ldots, {\bf x}_n) \leq CS_d(f)   ({\bf x}_1,\ldots, {\bf x}_n).
$$
This deduces the aimed \eqref{eqn: Lp app}  directly, and we complete the proof of Theorem \ref{thm:app}.
\end{proof}

\subsection{No necessity of cancellation}

Let $E$ be the mapping in \eqref{eqn:injective}  and let   $E(R)=\widetilde I\times \widetilde J\in \mathcal D^{[1], b_1}\times \mathcal D^{[2],b_2}$. Recall that one crucial consequence of the cancellation
\begin{equation}\label{eqn:crucial-cancel}
    \int_{\mathscr M_1}    a_R(f)  ({\bf x}_1,{\bf x}_2) {\rm d}\mu_1  ({\bf x}_1)=\int_{\mathscr M_2}    a_R(f)  ({\bf x}_1,{\bf x}_2){\rm d}\mu_2  ({\bf x}_2)=0
\end{equation}
 is to ensure that 
$$
 \Delta_{R'} (a_R(f))\equiv 0   \   \    \text{for any}  \   \    R'=I'\times J' \text{with} \   \    I'\supsetneqq \widetilde I\  \    \text{or} \  \  J'\supsetneqq \widetilde J,
$$
 which leads to \eqref{eqn:Haar-express-atom} and ensures the entire argument for Theorem \ref{thm:d-SWW} and Theorem \ref{thm:app} process.

We now address that the cancellation \eqref{eqn:crucial-cancel}
 is actually unnecessary. 
This implies that  Theorem \ref{thm:d-SWW} and Theorem \ref{thm:app} also hold for 
 more general settings. For example, we may replace the sub-Laplacian $\mathcal L_j$ ($j=1,2$)  by Schr\"odinger operator $\mathscr L_j=\mathcal L_j+V_j$ with certain potentials $V_j\geq0$, provided that the estimates \eqref{eqn:psi-scale},     \eqref{eqn:Holder} and \eqref{eqn:al-orth} are satisfied for the kernel of $\psi(t\sqrt{\mathscr L_j})$.

We now take Theorem \ref{thm:d-SWW} as an example, to illustrate that \eqref{eqn:crucial-cancel} is unnecessary.
 
We also use the notation given aforementioned. Without cancellation, we cannot expect \eqref{eqn:crucial-cancel}. To overcome it,  denote 
$$
    {\rm Avg}_{I} \left (a_{R}(f)\right ) ({\textbf x}_1,{\textbf x}_2)= \frac{1}{|\widetilde I|}\int_{\widetilde I}  a_R(f)({\textbf y}_1,{\textbf x}_2)  {\rm d}\mu_1({\textbf y}_1)   \cdot { \mathsf 1_{E(R)}} ({\textbf x}_1,{\textbf x}_2) =  \left(a_R(f)\right)_{\widetilde I} \,  {{\mathsf 1}_{E(R)}} ({\textbf x}_1,{\textbf x}_2),
$$
$$
    {\rm Avg}_{J}(a_{R}(f)) ({\textbf x}_1, {\textbf x}_2)=  (a_R(f))_{\widetilde J}  \mathsf 1_{E(R)} ({\textbf x}_1, {\textbf x}_2) \quad \text{and}\quad   {\rm Avg}_{R}(a_{R}(f)) ({\textbf x}_1,{\textbf  x}_2)=   (a_{R}(f))_{E(R)}  \mathsf 1_{E(R)} ({\textbf x}_1, {\textbf x}_2).
$$

Then we rewrite
\begin{align*}
   a_{R}(f)   =& \left[  a_{R}(f)-   {\rm Avg}_{I}(a_{R}(f)) - {\rm Avg}_{J}(a_{R}(f))  + {\rm Avg}_{R}(a_{R}(f))\right] +  \left[   {\rm Avg}_{I}(a_{R}(f)) -  {\rm Avg}_{R}(a_{R}(f))\right] \\
   & +   \left[   {\rm Avg}_{J}(a_{R}(f)) -  {\rm Avg}_{R}(a_{R}(f))\right] + {\rm Avg}_{R}(a_{R}(f)) \\
   =:&  a_R^{(1)}(f)+  a_R^{(2)}(f)+ a_R^{(3)}(f)+ a_R^{(4)}(f).
\end{align*}

Thus, we can decompose $F_{\Omega}=\sum_{j=1}^4 F^{(j)}_{\Omega} $ accordingly and now it suffices to prove that for $j=1,2,3,4$,
$$    \int_{\Omega} \exp\bigg(c_1  \frac {|F^{(j)}_{\Omega}({\bf x}_1,{\bf x}_2)|}{ \|S_\gamma(f)\|_{ L^\infty(\widetilde{\mathscr  M})} }\bigg) \,{\rm d} \mu_1({\bf x}_1)  {\rm d}\mu_2({\bf x}_2) \leq c_2  \mu(\Omega).
$$

Note that $\int a_R^{(1)}(f)\,    {\rm d} \mu_1({\bf x}_1)=\int a_R^{(1)}(f) \,  {\rm d}\mu_2({\bf x}_2)=0 $, hence $a_R^{(1)}(f)$ plays the role of $a_R(f)$ in Theorem \ref{thm:d-SWW} and we can deal with it similarly.

Let's consider  $a_R^{(4)}(f)$. Our strategy is to show that these atomic functions $\{a_R^{(4)}(f)\}_R$ (after classification by $E(R)$) meet the framework of Theorem \ref{thm:double-newCH-dyadic}. It suffices to verify the almost-orthogonality \textsf{(P3)}.  For any given $I_0\times  J_0\in \mathcal D^{[1],b_1}\times \mathcal D^{[2],b_2}$, assume 
$$
    \mathcal G_{ I_0\times  J_0}\subseteq \left\{E(R)\in  \mathcal D^{[1],b_1}\times \mathcal D^{[2],b_2}:\,  E(R)\subseteq  I_0\times J_0\right\}.
$$ 
Recall that the left-side hand of \eqref{eqn:double-P3-revise-I} in \textsf{(P3)} can be weaken to $L^1$ norm (see Remark \ref{rem:P3}). For any fixed ${\bf x}_2\in J_0$,  
from the definition  of $a_R^{(4)}(f)$, 
\begin{align}\label{eqn: average I0}
	&\frac{1}{| I_0|} \int_{ I_0} \bigg|  \sum_I \sum_{J:\   R=I\times J, \,E(R)\in  \mathcal G_{  I_0\times  J_0}}  a_R^{(4)} (f)  ({\bf x}_1,{\bf  x}_2)\bigg| \, {\rm d} \mu_1( {\bf x}_1)\nonumber \\
	\leq\,& 	 
	\frac{1}{| I_0|} \int_{ I_0} \sum_I \sum_{J:\   R=I\times J, \,E(R)\in  \mathcal G_{I_0\times  J_0}}     {\rm sgn}(a_R^{(4)} (f))\,     (a_{R}(f))_{E(R)} {\mathsf 1_{E(R)}} ({\textbf x}_1, {\textbf x}_2)\, {\rm d} \mu_1( {\bf x}_1).
\end{align}
Note that for every fixed ${\bf x}_2$,
$$ \int_{ I_0} (a_{R}(f))_{E(R)} \mathsf 1_{E(R)}({\textbf x}_1, {\textbf x}_2)\, {\rm d} \mu_1( {\textbf x}_1)=\int_{\widetilde I}\, (a_{R}(f)({\textbf y}_1,\cdot))_{\widetilde J} \,   
{\rm d}\mu_1({\textbf y}_1)\ \ \mathsf 1_{\widetilde J}({\textbf x}_2).$$
Thus, the right-hand side of \eqref{eqn: average I0} is written as 
\begin{align}
	& 	 
	\frac{1}{| I_0|} \int_{ I_0} \sum_I \sum_{J:\   R=I\times J,\,E(R)\in  \mathcal G_{I_0\times  J_0}}     {\rm sgn}(a_R^{(4)} (f))\,     (a_{R}(f)({\textbf y}_1,\cdot))_{\widetilde J}  \mathsf 1_{E(R)} ({\textbf y}_1, {\textbf x}_2)
  \, {\rm d} \mu_1( {\bf y}_1) \nonumber\\	
   \leq \,& \bigg(	 
	\frac{1}{| I_0|} \int_{ I_0} \bigg| \sum_I \sum_{J:\   R=I\times J,\, E(R)\in  \mathcal G_{I_0\times  J_0}}     {\rm sgn}(a_R^{(4)} (f))\,     (a_{R}(f)({\textbf y}_1,\cdot))_{\widetilde J}  \mathsf 1_{E(R)} ({\textbf y}_1, {\textbf x}_2)\bigg|^2
  \, {\rm d} \mu_1( {\bf y}_1)\bigg)^{1/2} \nonumber\\	
	\leq \, &  C \bigg(\frac{1}{|I_0|} \int_{ I_0} \sum_I  \Big|  \mathcal A_{I} (f)({\bf y}_1, {\bf x}_2)  \Big|^2 \, {\rm d} \mu_1({\bf y}_1)\bigg)^{1/2},   \tag{\textsf{P3-1}}
\end{align}
 where in the last inequality we used the almost orthogonality estimate \eqref{eqn:al-orth}, and  
 $$
      \Big|  \mathcal A_{I} (f) (  {\textbf y}_1, {\textbf x}_2) \Big|^2  :=|\widetilde I|^{-1}  \int_{I_+} \bigg|  \sum_{J:\   R= I\times J,\, E(R)\in  \mathcal G_{I_0\times  J_0}}     {\rm sgn}(a_R^{(4)} (f))\,        \Big(\mathfrak a_J(f)({\textbf z}_1,\cdot,t_1)\Big)_{\widetilde J} \, \mathsf 1_{\widetilde J}  ( {\textbf x}_2)  \bigg|^2  \frac{{ \rm d} \mu_1( {\textbf  z}_1)   {\rm d}t_1}{t_1}
       \, \mathsf 1_{\widetilde I} ({\textbf y}_1),
 $$
 and 
 $$
     \mathfrak a_J(f)( {\bf z}_1, {\bf z}_2,t_1):=\int_{J_+}  K_{\psi\big(t_2\sqrt{\mathcal L_2}\big)}({\bf z}_2, {\bf w}_2) Q_{t_1^2}^{[1]} Q_{t_2^2}^{[2]} (f) ({\bf z}_1, {\bf w}_2)  \frac{ {\rm d}\mu_2( {\bf w}_2)  {\rm d}t_2}{t_2} 
 $$
 with $Q_s^{[j]}=2s{\displaystyle{\partial }\over\displaystyle \partial s}e^{-s\mathcal L_j}$, $j=1,\,2$.
 
Moreover,
 \begin{align}\label{eqn:a4}
 	&\int_{J_0}    \Big|  \mathcal A_{I} (f)({\textbf y}_1, {\textbf x}_2)  \Big|^2 \, {\rm d} \mu_2( {\textbf x}_2)\nonumber\\
 	 \lesssim \,  &     |\widetilde I|^{-1}  \int_{I_+}        \int_{J_0}     \bigg|  \sum_{J:\   R= I\times J,\, E(R)\in  \mathcal G_{ I_0\times  J_0}}     \mathfrak b_J(f) ({\textbf z}_1, {\textbf x}_2, t_1) \bigg|^2 {\rm  d} \mu_2({\textbf x}_2) \frac{{ \rm d} \mu_1({\textbf  z}_1)   {\rm d}t_1}{t_1}
        \mathsf 1_{\widetilde I} ({\textbf y}_1)\nonumber\\
        &+   |\widetilde I|^{-1}  \int_{I_+}        \int_{J_0}     \bigg|  \sum_{J:\   R= I\times J,\, E(R)\in  \mathcal G_{ I_0\times  J_0}}     {\rm sgn}(a_R^{(4)} (f))\,          \mathfrak a_J(f)({\textbf z}_1, {\textbf x}_2, t_1) \mathsf 1_{\widetilde J}  ({\textbf x}_2)  \bigg|^2  {\rm d} \mu_2({\textbf x}_2) \frac{{ \rm d}  \mu_1({\textbf z}_1)   {\rm d}t_1}{t_1}
        \mathsf 1_{\widetilde I} ({\textbf y}_1),
 \end{align}
 where 
 $$
      \mathfrak b_J(f) ({\textbf z}_1,{\textbf x}_2, t_1) =  {\rm sgn}(a_R^{(4)} (f))\,        \left[\Big(\mathfrak a_J(f)({\textbf z}_1,\cdot, t_1)\Big)_{\widetilde J} - \mathfrak a_J(f)({\textbf z}_1, {\textbf x}_2,  t_1)\right] \mathsf 1_{\widetilde J}  ({\textbf x}_2) .
 $$
 
  For any given ${\textbf z}_1$ and $t_1$, due to the cancellation of $  F_1({\textbf z}_1,{\textbf x}_2, t_1):=\sum_{J:\   R= I\times J,\, E(R)\in  \mathcal G_{ I_0\times  J_0}}     \mathfrak b_J(f) ({\textbf z}_1,{\textbf x}_2, t_1)$ on ${\textbf x}_2$, we can write it as a sum of martingale differences 
 $$\sum_{J'\in \mathcal D^{[2],b_2}}\Delta_{J'}( F_1 )=\sum_{J'\in \mathcal D^{[2],b_2}  , \,J'\subseteq J_0    }\Delta_{J'}( F_1 ) ,
 $$
 and the argument in the above subsection shows that the corresponding martingale square function satisfies  
 \begin{align*}
     \sum_{J'\in \mathcal D^{[2],b_2}  , J'\subseteq J_0    }\Delta_{J'}^2( F_1 ({\textbf z}_1,\cdot,  t_1))  ({\textbf x}_2)  \leq C
     \sum_{  J:\   R= I\times J, \, E(R)\in  \mathcal G_{ I_0\times  J_0}  }   \frac{|\lambda_R(f)|^2}{\mu(R)} {\mathsf{ 1}}_{E(R)}({\textbf z}_1,{\textbf x}_2) ,
 \end{align*}
 which is also provided by the cancellation of  $F_1$.
Thus
\begin{align*}
    &  |\widetilde I|^{-1}  \int_{I_+}        \int_{J_0}     \bigg|  \sum_{J:\   R= I\times J\in  \mathcal G_{ I_0\times  J_0}}     \mathfrak b_J(f) ({\textbf z}_1,{\textbf x}_2, t_1) \bigg|^2 {\rm  d} \mu_2({\textbf x}_2) \frac{{ \rm d} \mu_1({\textbf  z}_1)   {\rm d}t_1}{t_1}
        \mathsf 1_{\widetilde I} ({\textbf y}_1)\\
        =\, & |\widetilde I|^{-1}  \int_{I_+}        \int_{J_0}     \sum_{J'\in \mathcal D^{[2],b_2}  , J'\subseteq J_0    }\Delta_{J'}^2( F_1 ({\textbf z}_1,\cdot,  t_1))  ({\textbf x}_2) \,{\rm  d} \mu_2({\textbf x}_2) \frac{{ \rm d} \mu_1({\textbf  z}_1)   {\rm d}t_1}{t_1}
        \mathsf 1_{\widetilde I} ({\textbf y}_1)\\
        \leq \,&  C \int_{J_0}    \sum_{  J:\   R= I\times J, \, E(R)\in  \mathcal G_{ I_0\times  J_0}  } \frac{|\lambda_R(f)|^2}{\mu(R)} {\mathsf{ 1}}_{E(R)}({\textbf y}_1,{\textbf x}_2)    \, {\rm  d} \mu_2({\textbf x}_2) .
\end{align*}

Meanwhile, it follows from \eqref{eqn:al-orth} again to obtain 
\begin{align*}
      &\int_{J_0}     \bigg|  \sum_{J:\   R= I\times J,\, E(R)\in  \mathcal G_{ I_0\times  J_0}}     {\rm sgn}(a_R^{(4)} (f))\,          \mathfrak a_J(f)({\textbf z}_1,{\textbf x}_2, t_1) \mathsf 1_{\widetilde J}  ({\textbf x}_2)  \bigg|^2  {\rm d} \mu_2({\textbf x}_2)\\
      \leq\,& C\int_{J_0}     \sum_{  J:\   R= I\times J, \, E(R)\in  \mathcal G_{ I_0\times  J_0}  } \frac{|\lambda_R(f)|^2}{\mu(R)} {\mathsf{ 1}}_{E(R)}({\textbf y}_1,{\textbf x}_2)\,   {\rm  d} \mu_2({\textbf x}_2).
\end{align*}

Hence, by \eqref{eqn:a4}, for any given $I$ and ${\textbf y}_1$, 
\begin{align} 
 	\int_{J_0}    \Big|  \mathcal A_{I} (f)({\textbf y}_1, {\textbf x}_2)  \Big|^2 \, {\rm d} \mu_2( {\textbf x}_2) 	\leq \,  C \int_{J_0}  \sum_{  J:\   R= I\times J, \, E(R)\in  \mathcal G_{ I_0\times  J_0}  }  \frac{|\lambda_R(f)|^2}{\mu(R)} {\mathsf{ 1}}_{E(R)}({\textbf y}_1,{\textbf x}_2)   \,  {\rm  d} \mu_2({\textbf x}_2) .   \tag{\textsf{P3-2}}
\end{align}

Combining \eqref{eqn: average I0}, \textsf{(P3-1)} and \textsf{(P3-2)}, we 
obtain a variant of \textsf{(P3)}, by noting that $$|a_R^{(4)}(f)|\leq C \lambda_R(f)\mu(R)^{-1/2}\mathsf 1_{E(R)}.$$
This allows us to apply the  argument of Theorem \ref{thm:double-newCH-dyadic} to finish the remaining proof concerning $a_R^{(4)}(f)$.

%
 
 For  $a_R^{(2)}(f)$, the argument is similar to that of $a_R^{(4)}(f)$, and it will only exist the first term in \eqref{eqn:a4}.   The argument for   $a_R^{(3)}(f)$ is also similar by considering the integral on the second variable ${\bf x}_2$ first.  Due to the related arguments that have been discussed above, we omit the details here. 
 
Eventually,  we can derive the analogs of  Theorem \ref{thm:d-SWW} and Theorem \ref{thm:app} without assuming the cancellation.

\section{Appendix: Hoeffding inequalities  associated to  martingale differences}\label{sec:6}
\setcounter{equation}{0}
In this appendix, we  consider rougher square functions associated to martingale differences, defined by the $L^\infty$ structure.  Obviously this  square function is pointwise equal to or larger than our square function. This, in combination with the cancellation of martingale differences, will provide alternative  approaches to Hoeffding inequalities and related topics. However, we address the fact that such $L^\infty$ structure is not applicable to multi-parameter theory.

\subsection{One-parameter Hoeffding inequality associated to  martingale differences}\label{subsec:6.1}

To clarify the differences between the classical iterative approach based on generations of cubes and our approach based on the scales of square functions,  
this subsection will provide a detailed proof of the classical Hoeffding inequality for  martingale differences \eqref{eqn:classical-dyadic-Hoeffding} 
 in spaces of homogeneous type. Concretely,  let  $(X,\rho,\mu)$ be a quasi-metric space satisfying the volume doubling condition, and let $\mathcal D$ be a system of dyadic cubes in $X$. Let $f\in L^p(X,{\rm d}\mu)$ for $1<p<\infty$.
For any given $Q\in \mathcal D$ and $\alpha>0$, we have 
 \begin{equation}\label{eqn:classical-dyadic-Hoeffding-SHI}
	\int_Q \exp\Bigg (\alpha  \mathcal F_N(x)  -\frac{\alpha^2}{2} \widetilde S^2(\mathcal  F_N)(x) \Bigg) \,  {\rm d}\mu(x)\leq \mu(Q),
\end{equation}
where  
$$
   \mathcal F_N(x)=\sum_{\substack{P\in  Ch^{(n)}(Q)\\ 0\leq n\leq N}}   \Delta_{P}(f) ( x),\qquad \widetilde S^2(\mathcal F_N)(x)=\sum_{\substack{P\in  Ch^{(n)}(Q)\\ 0\leq n\leq N}}   \left \|\Delta_{P}(f)\right\|_{L^\infty(P)}^2 \mathsf 1_P(x) , 
$$
and 
$\Delta_P(f)$ the martingale difference of $f$ associated with $P$, given by 
$$
   \Delta_P(f)(x)=\sum_{R\in Ch(P)}\left(f_R-f_P\right)\mathsf{1}_R(x).
$$  
The above  \eqref{eqn:classical-dyadic-Hoeffding-SHI} is a natural generalization of the Hoeffding inequality \eqref{eqn:ESI} in $[0,1]$ and its extension \eqref{eqn:ESI higher}  (alternatively,  the equivalent expression \eqref{eqn:classical-dyadic-Hoeffding}) for higher dimensional Euclidean spaces.

Rubin's method is based on the following Hoeffding's lemma. 

\begin{lemma}  \label{lem:Hoeffding}
Let $(X, \nu)$ be a probability space and let $g: X \mapsto \mathbb R$ be a random variable (i.e., a measurable function) such that $|g(x)| \leq \mathfrak t$ for all $x \in X$ and
$$
    \int_{X} g(x) \, {\rm d} \nu(x)=0 .
$$
Then
$$
    \int_{X} \exp (g(x))\, {\rm d} \nu(x) \leq \cosh (\mathfrak t)\leq \exp\left(\frac{\mathfrak t^2}{2}\right).
$$
\end{lemma}

\begin{proof}
Here we provide a new proof.
We write
\begin{equation}\label{eqn:g-a}
   g(x)=\alpha(x)\cdot \mathfrak t+(1-\alpha(x))\cdot (-\mathfrak t) 
\end{equation}
for each $x\in X$. That is, $\displaystyle \alpha(x)=\frac{g(x)+ \mathfrak t }{ 2 \mathfrak t}$ is a measurable function on $X$, and it follows from the cancellation of $g$ to see
$$
    \int_X \alpha(x) \, {\rm d} \nu(x)=\int_X \frac{g(x)+\mathfrak t}{ 2 \mathfrak t}\, {\rm d} \nu(x)=\frac{1}{2}.
$$

Note that $ \alpha(x)\in [0,1]$ for each $x\in X$, then by \eqref{eqn:g-a}, one may apply the convexity  of the exponential function to obtain
\begin{align*}
	\int_X \exp(g(x)) \, {\rm d} \nu(x)&\leq  \int_X  \left[ \alpha (x)  e^{\mathfrak t}  +(1-\alpha(x)) e^{-\mathfrak t}\right]  {\rm d} \nu(x)\\
	&=  e^{\mathfrak t}      \int_X    \alpha (x)   \, {\rm d} \nu(x)  +    e^{- \mathfrak t}   \int_X    \left[1-\alpha (x) \right]   {\rm d} \nu(x)\\
	&=\frac{e^{\mathfrak t} +e^{-\mathfrak t}}{2},
\end{align*}
as desired.
\end{proof}

\begin{remark}
The above Hoeffding's lemma is sharp. For example, when $X=[0,1]$, $\nu$ is the Lebesgue measure on $[0,1]$, and taking
$$
    g(x)=
    \begin{cases}
    	-\mathfrak t,\quad  & x\in \big [0,\frac{1}{2}\big ],\\[6pt]
    	\mathfrak t,\quad & x\in \big( \frac{1}{2},1\big].
    \end{cases}
$$
Then 
$$
    \int_{X} \exp (g(x))\,{\rm  d} \nu(x) =\frac{e^{\mathfrak t} +e^{-\mathfrak t}}{2}= \cosh (\mathfrak t).
$$
\end{remark}

\smallskip

\begin{proof}[Proof of \eqref{eqn:classical-dyadic-Hoeffding-SHI}]
Rubin's argument is mainly the following iteration on generations of cubes: for any $k\geq 1$,
\begin{align}\label{eqn:iteration-Rubin}
	\int_Q \exp\left (\alpha \mathcal F_{k}( x)  -\frac{\alpha^2}{2} \widetilde {S}^2(\mathcal F_k)(x)  \right )  \,  {\rm d}\mu(x)
	\leq  \int_Q \exp\left (\alpha \mathcal F_{k-1}( x)  -\frac{\alpha^2}{2} \widetilde {S}^2(\mathcal F_{k-1})(x)  \right )  \,  {\rm d}\mu(x).
\end{align}

To this end, note that for each $P\in \mathcal D$, the martingale difference $\Delta_P(f)$ satisfies the following properties:
\begin{itemize}
	\item $\Delta_P(f)$ is supported on $P$;
	\item $\Delta_P(f)$ is constant on each $R\in Ch(P)$;
	\item $\int_P   \Delta_P(f)(x)\,{\rm d}\mu(x)=0$.
\end{itemize}
Hence $\mathcal F_{k-1}\equiv \mathcal F_{k-1}(c_P)$ and $\widetilde {S}^2(\mathcal F_{k-1})\equiv \widetilde {S}^2(\mathcal F_{k-1})(c_P)$  on each $P\in Ch^{(k)}(Q)$, where $c_P$ is the center of $P$, and it follows from Lemma \ref{lem:Hoeffding} that
$$
    \int_P \exp\left(\alpha \Delta_P(f)(x)  -\frac{\alpha^2}{2}\|\Delta_P(f)\|_{L^\infty(P)}^2 \right) {\rm d}\mu(x)\leq \mu(P).
$$
As a consequence,
\begin{align*}
	&\int_Q \exp\left (\alpha \mathcal  F_{k}( x)  -\frac{\alpha^2}{2} \widetilde {S}^2(\mathcal F_k)(x)  \right )  \,  {\rm d}\mu(x)\\
	 = &  \int_Q \exp\left (\alpha \mathcal F_{k-1}( x)  -\frac{\alpha^2}{2} \widetilde {S}^2(\mathcal F_{k-1})(x)   +\alpha \sum_{ P\in  Ch^{(k)}(Q)} \Delta_{P}(f)(x)-\frac{\alpha^2}{2}    \sum_{ P\in  Ch^{(k)}(Q)} \|\Delta_{P}(f)\|_{L^\infty(P)}^2 \mathsf 1_P(x)  \right )  \,  {\rm d}\mu(x)\\
	 =& \sum_{P\in Ch^{(k)} (Q)} \exp\left (\alpha \mathcal F_{k-1}( c_P)  -\frac{\alpha^2}{2} \widetilde {S}^2(\mathcal F_{k-1})(c_P)  \right )  \int_P  \exp\left ( \alpha  \Delta_{P}(f)(x)-\frac{\alpha^2}{2}     \|\Delta_{P}(f)\|_{L^\infty(P)}^2 \mathsf 1_P(x)  \right )  \,  {\rm d}\mu(x)\\
	 \leq &  \sum_{P\in Ch^{(k)} (Q)} \exp\left (\alpha\mathcal  F_{k-1}( c_P)  -\frac{\alpha^2}{2} \widetilde {S}^2(\mathcal F_{k-1})(c_P)  \right )  \mu(P)\\
	 =& \int_Q \exp\left (\alpha \mathcal F_{k-1}( x)  -\frac{\alpha^2}{2} \widetilde {S}^2(\mathcal F_{k-1})(x)  \right )  \,  {\rm d}\mu(x).
\end{align*}
That is, \eqref{eqn:iteration-Rubin} holds.  

Therefore, for any $N$, one may apply this iteration on generations of cubes to see
\begin{align*}
   \text{LHS \   of \   } \eqref{eqn:classical-dyadic-Hoeffding-SHI} 
  \leq\,  & \int_Q \exp\Bigg (\alpha  \sum_{\substack{P\in  Ch^{(n)}(Q)\\ 0\leq n\leq N-1}}   \Delta_{P}(f) ( x)  -\frac{\alpha^2}{2} \sum_{\substack{P\in  Ch^{(n)}(Q)\\ 0\leq n\leq N-1}}   \|\Delta_{P}(f)\|_{L^\infty(P)}^2 \mathsf 1_P(x)  \Bigg)  \,  {\rm d}\mu(x)\\
  \leq\, & \cdots\cdots\\
  \leq \,&  \int_Q \exp\Bigg (\alpha \Delta_Q(f)(x)-\frac{\alpha^2}{2}\|\Delta_Q(f)\|_{L^\infty(Q)}^2 \mathsf 1_Q(x)\Bigg)\,  {\rm d}\mu(x)\\
  \leq\, &  \text{RHS \   of \   } \eqref{eqn:classical-dyadic-Hoeffding-SHI},
\end{align*}
as desired.
\end{proof}

\smallskip

\subsection{Bi-parameter Hoeffding inequality associated to  martingale differences}\label{subsec:6.2}
For  $i=1,2$, let  $(X_i,\rho_i,\mu_i)$ be as in Subsection \ref{subsec:6.1}, 
and let $\mathcal D^{[i]}$ be a system of dyadic cubes on $X_i$. 
Let $\{\Delta_{P_1\times P_2}(f)\}$ be double martingale differences of some function $f$  on a product space $X_1\times X_2$ of homogeneous type,  the argument in \cite[Lemma 2.2]{Pi86}  based on the cancellation property of $\Delta_{P_1\times P_2}(f)$
 can deduce that there exist positive constant $c$ and $C$ such that for any $\alpha>0$ and $x_1\in Q_1$,
\begin{align}\label{eqn:one-shift-bi-I-weak-martingale} 
		\int_{Q_2} \exp\left(   \alpha S_1( \mathcal F_{N,N})     (x_1,x_2) -c \alpha^2 \widetilde{S}^2(\mathcal F_{N,N})(x_1,x_2)\right){\rm d}\mu_2(x_2)\leq C\mu_2(Q_2),
\end{align}
where 
$$
    \mathcal F_{N,N}(x_1,x_2):=\sum_{\substack{P_1\in  Ch^{(n)}(Q_1)\\ 0\leq n\leq N}} \sum_{\substack{P_2\in  Ch^{(n)}(Q_2)\\ 0\leq n\leq N}} \Delta_{P_1\times P_2}(f)(x_1,x_2)
$$
and 
$$
     \widetilde{S}^2(\mathcal F_{N,N})(x_1,x_2):=\sum_{\substack{P_1\in  Ch^{(n)}(Q_1)\\ 0\leq n\leq N}} 
     \sum_{\substack{P_2\in  Ch^{(n)}(Q_2)\\ 0\leq n\leq N}} \left \|\Delta_{P_1\times P_2}(f)\right\|_{L^\infty(P_1\times P_2)}^2\mathsf 1_{P_1\times P_2}(x_1,x_2).
$$
It's clear that $\widetilde{S}(\mathcal  F_{N,N})\geq S(\mathcal F_{N,N})$ pointwise.

Observe that the lacunary version also holds by using the same argument. In particular, for any  $\mathcal G_{Q_1, Q_2} \subseteq  \bigcup_{0\leq n\leq N}Ch^{(n)}(Q_1)\times \bigcup_{0\leq n\leq N}Ch^{(n)}(Q_2)$, $x_1\in Q_1$ and $\alpha>0$,
\begin{align}\label{eqn:one-shift-bi-lacu}
	&\int_{Q_2} \exp\Bigg[\alpha \bigg(\sum_{P_1}  \bigg|\sum_{P_2:\, P_1\times P_2\in \mathcal G_{Q_1,Q_2}} \Delta_{P_1\times P_2}(f)(x_1, x_2) \bigg|^2  \, \bigg)^{1/2}\Bigg.   \nonumber\\
	&\qquad\qquad    \Bigg. -c\alpha^2 \sum_{P_1\times P_2\in \mathcal G_{Q_1,Q_2}}   \left\|\Delta_{P_1\times P_2} (f)  \right\|_{L^\infty(P_1\times P_2)}^2 \Bigg]\,  {\rm d}\mu_2(x_2)\leq C\mu_2(Q_2).
\end{align}

Based on it, we can establish the following improved good-$\lambda$ inequality, the corresponding one-parameter counterpart is well-known.  The special case when $X_1\times X_2=[0,1]^2$ was previously addressed in \cite{Pi86} by the second author. However, due to the nonlinear property of the operator $S_1$ and the fact that $|\Delta_{P_1\times P_2}(f)|$ may take different values on each child of $P_1\times P_2$, the subsequent argument is more intricate.

\begin{proposition}\label{coro:good-lambda}
There exist positive  constants  $C$ and $C'$ and sufficiently small $\epsilon$,
such that for any $x_1\in Q_1$ and $\lambda>0$,
\begin{align}\label{eqn:good-lambda}
	&\mu_2\Big(\Big \{x_2\in Q_2: S_1^*(\mathcal  F_{N,N})(x_1,x_2)>2\lambda, \  \   \widetilde S(\mathcal  F_{N,N})(x_1,x_2)<\epsilon \lambda \Big \}\Big)\nonumber\\
	\leq \,& C \exp\left(-\frac{C'}  {\epsilon^2}\right) \mu_2\Big(\Big\{x_2\in Q_2: S_1^*(\mathcal F_{N,N})(x_1,x_2)>\lambda  \Big \}\Big),
\end{align}
where
$$
  S_1^*(\mathcal  F_{N,N})(x_1,x_2):=\sup_{0\leq k\leq N} S_1(\mathcal  F_{N,k})(x_1,x_2)=\sup_{0\leq k\leq N}\Bigg(\sum_{\substack{P_1\in  Ch^{(n)}(Q_1)\\ 0\leq n\leq N}} \Big|\sum_{\substack{P_2\in  Ch^{(n)}(Q_2)\\ 0\leq n\leq k}}  \Delta_{P_1\times P_2}(f) ( x_1, x_2)\Big|^2\, \Bigg)^{1/2}.
$$
\end{proposition}

\begin{proof}
Fixed $x_1\in Q_1$ and $\lambda>0$, one may assume that $
    \Big\{x_2\in Q_2: S_1^*(\mathcal  F_{N,N})(x_1,x_2)>\lambda  \Big \}\neq \emptyset
$. 

Moreover, without loss of generality, 
it suffices to assume 
\begin{equation}\label{eqn:small-assume}
	 \max_{x_2\in Q_2} S_1(\mathcal  F_{N,0})(x_1,x_2)<\lambda.
\end{equation}
Indeed, $S_1(\mathcal F_{N,0})(x_1,\cdot)$  is constant on each  $P_2\in Ch(Q_2)$. If $S_1(\mathcal  F_{N,0})(x_1,\cdot)\geq \lambda$ on some $P_2\in Ch(Q_2)$, then 
$$
       \widetilde S(\mathcal  F_{N,N})(x_1,\cdot)\geq \widetilde S(\mathcal  F_{N,0}) (x_1,\cdot)\geq S(\mathcal  F_{N,0}) (x_1,\cdot)=
       S_1(\mathcal  F_{N,0})(x_1,\cdot)\geq \lambda \   \   \text{on}\   \   P_2.
$$
Note that $\epsilon>0$ in \eqref{eqn:good-lambda} is sufficiently small, it turns to prove that  \eqref{eqn:good-lambda} holds on $Q_2\setminus P_2$.

Denote $\mathcal G_2=\bigcup_{0\leq n\leq N} Ch^{(n)}(Q_2)$ and let
$$
    \mathsf B:=\left\{R_2\in Ch^{(N)}(Q_2):\     \left \|S_1^*(\mathcal F_{N,N})(x_1,\cdot)\right \|_{L^\infty (R_2)}^2 >\lambda^2\right\},
$$
clearly $ \displaystyle\Big\{x_2\in Q_2: S_1^*(\mathcal F_{N,N})(x_1,x_2)>\lambda  \Big \}\subseteq \bigcup_{R\in \mathsf B}R$.

We now provide a constructive algorithm to obtain a precise decomposition of the level set\\ $ \displaystyle\Big\{x_2\in Q_2: S_1^*(\mathcal  F_{N,N})(x_1,x_2)>\lambda  \Big \}$
into pairwise disjoint cubes.

\medskip

\noindent {\bf Algorithm:}

\medskip

For each  $R_2\in \mathsf B$ and  the  corresponding chain $\text{Chain}(R_2)$ 
$$
     Q_2=P_{2,0}\supseteq P_{2,1}\supseteq \cdots\supseteq P_{2,k}\supseteq P_{2,k+1}\supseteq\cdots\supseteq P_{2,N}=R_2
$$
generated by $R_2$  in $\mathcal G_2$, there must be finite (at least one) critical cubes
$$
         P_{2,C_1}\supseteq P_{2,C_2}\supseteq \cdots \supseteq  P_{2, C_K}
$$
in this chain, such that 
\begin{equation}\label{eqn:C2-new}
	   \left \|S_1(\mathcal  F_{N, C_k})(x_1,\cdot)\right \|_{L^\infty(P_{2, C_k})} ^2 >\lambda^2 ,\quad 1\leq k\leq K.
\end{equation}
Otherwise, $S_1(\mathcal  F_{N, k})(x_1,\cdot)\leq \lambda $ on $P_{2,k}$ for each $1\leq k\leq N$, consequently,   $S_1(\mathcal  F_{N, k})(x_1,\cdot)\leq \lambda $ on $R_2$ for each $k$, and thus $S_1^*(\mathcal  F_{N, N})(x_1,\cdot)\leq \lambda $ on $R_2$, this contradicts  $R_2\in \mathsf B$.

Note that $P_{2,C_1}$ is the first cube in $\text{Chain}(R_2)$ satisfying \eqref{eqn:C2-new} and $C_1>0$ due to the assumption \eqref{eqn:small-assume}, this implies the small-scale property that 
$$
     \left \|S_1(\mathcal F_{N, k})(x_1,\cdot)\right \|_{L^\infty(P_{2, k})} ^2 \leq \lambda^2 ,\quad 0\leq k\leq C_1-1.
$$
Similarly, for any two adjacent critical cubes $P_{2,C_i}$ and $P_{2, C_{i+1}}$, if $C_{i+1}> C_i+1$, then all cubes $P_{2,k}$ in $\text{Chain}(R_2)$ with $C_i<k<C_{i+1}$
  satisfy the above small-scale property.

Now we label  $C_{R_2}^I$ and $C_{R_2}^{II}$  associated to $R_2$ from the following setting, which will be crucial for our subsequent argument.
\begin{itemize}
	\item Verify from the first critical cube $P_{2,C_1}$, 
	\begin{itemize}
	\item {\it Case 1}.  if $ S_{1}^2( \mathcal  F_{N, C_1})(x_1,\cdot)  >\lambda ^{2}$ on $P_{2, C_1}$, then let 
	$$
	   C_{R_2}^I=P_{2, C_1}\quad \text{and}\quad C_{R_2}^{II}=P_{2, C_1}.
	$$
	\item {\it Case 2}. if $ \displaystyle \min_{x_2\in P_{2, C_1}}S_{1}^2( \mathcal  F_{N, C_1})(x_1,x_2)     \leq \lambda ^{2}$ while  $ S_{1}^2(\mathcal   F_{N, C_1}) (x_1,\cdot)  >\lambda ^{2}$ on $P_{2,C_1+1}$, then 
	let 
	$$
	    C_{R_2}^I=P_{2, C_1}\quad \text{and}\quad C_{R_2}^{II}=P_{2, C_1+1}.
	$$
	
	In particular, if $P_{2,C_1}=R_2$ in this case, precisely, that is,  $\displaystyle \max_{x_2\in R_2}S_{1}^2( \mathcal  F_{N, C_1})(x_1,x_2)     > \lambda ^{2}$ while $\displaystyle \min_{x_2\in R_2}S_{1}^2( \mathcal  F_{N, C_1})(x_1,x_2)     \leq \lambda ^{2}$. It means  $S_{1}^2( \mathcal F_{N, N}) (x_1,\cdot)   >\lambda ^{2}$ on at least one but not all children of $R_2$, let 
	$$
	    C_{R_2}^I=R_2\quad \text{and}\quad C_{R_2}^{II}=\bigcup_{\left\{W\in Ch(R_2):\,   S_{1}^2(  \mathcal  F_{N, N}) (x_1,\cdot)  >\lambda ^{2}    \    \text{on}\   W\right\}} W.
	$$
	Recall that now $C_{R_2}^{II}\subsetneqq R_2$ with $C_{R_2}^{II}\in Ch(R_2)$ or $C_{R_2}^{II}$ is the union of certain (at least two but not all) children of $R_2$.
	
	 \item {\it Case 3}. Otherwise, $ \displaystyle \min_{x_2\in P_{2, C_1}}S_{1}^2( \mathcal  F_{N, C_1})(x_1,x_2)     \leq \lambda ^{2}$ and   $ S_{1}^2( \mathcal    F_{N, C_1}) (x_1,\cdot)  >\lambda ^{2}$ on some (at least one) $P'\in Ch(P_{2,C_1})$, but $P'$ is not in the $\text{Chain}(R_2)$. Such $P'$ must be in the chain generated by another $R'_2\in \mathsf B$, and it will be verified by the same method when handling $\text{Chain}(R'_2)$ to label  $C_{R'_2}^I$ and $C_{R'_2}^{II}$.
	 
	 In this case, we will turn to verify the next critical cube $P_{2,C_2}$.  Indeed,  by contradiction, such $P_{2,C_2}$ must be exist due to $R_2\in \mathsf B$. 
\end{itemize}
\smallskip

\item For any $k\geq 2$, once {\it Case 3} for $P_{2,C_{k-1}}$ happens, we will verify 
$P_{2,C_k}$  in the same manner. That is, 
\begin{itemize}
\item if  the above {\it Case 1} happens  (replacing $C_1$ therein by $C_k$), let
$$
     C_{R_2}^I=P_{2, C_k}\quad \text{and}\quad C_{R_2}^{II}=P_{2, C_k}.
$$

\item If   the above {\it Case 2} happens  (replacing $C_1$ therein by $C_k$), 
then let
$$
    C_{R_2}^I=P_{2, C_k}\quad \text{and}\quad C_{R_2}^{II}=P_{2, C_k+1}
$$
when $C_k<N$, and let
$$
    C_{R_2}^I=R_2\quad \text{and}\quad C_{R_2}^{II}=\bigcup_{\left\{W\in Ch(R_2):\,   S_{1}^2( \mathcal  F_{N, N}) (x_1,\cdot)  >\lambda ^{2}    \    \text{on}\   W\right\}}W
$$
when $C_k=N$.

\item If  the above {\it Case 3} happens  (replacing $C_1$ therein by $C_k$), we turn to verify the next critical cube $P_{2,C_{k+1}}$.
\end{itemize}

\item 
Eventually, combining the fact that the critical cubes $\{P_{2,C_k}\}_{1\leq k\leq K}$ are finite and $R_2\in \mathsf B$, we can conclude by a standard contradiction argument that it's impossible for 
{\it Case 3}  to occur for each $P_{2,C_k}$, thus we can  obtain $C_R^{I}$ and $C_R^{II}$ clearly in a finite number of steps. See Figure \ref{fig:2} for instance.
\end{itemize}

\noindent {\bf End of Algorithm}

\medskip

\begin{figure}
\tikzset{every picture/.style={line width=0.75pt}}  
\begin{tikzpicture}[x=0.75pt,y=0.75pt,yscale=-1,xscale=1]

\draw   (1.67,3) -- (112,3) -- (112,113.33) -- (1.67,113.33) -- cycle ;
 
\draw   (5.67,123) -- (101,123) -- (101,218.33) -- (5.67,218.33) -- cycle ;
 
\draw    [fill={rgb, 255:yellow, 255; green, 0; blue, 0}  ,fill opacity=1 ] (16.33,355.5) -- (91,355.5) -- (91,430.17) -- (16.33,430.17) -- cycle ;
 
\draw (21,440) -- (86,440) -- (86,505) -- (21,505) -- cycle ;
 
\draw   (31,610) -- (77,610) -- (77,656) -- (31,656) -- cycle ;
 
\draw   (38,702) -- (65,702) -- (65,729) -- (38,729) -- cycle ;
 
\draw [line width=1.5]  [dash pattern={on 5.63pt off 4.5pt}]  (53,223) -- (53,238.5) -- (53,256.5) ;
 
\draw [line width=1.5]  [dash pattern={on 5.63pt off 4.5pt}]  (52,662) -- (52,677.5) -- (52,695.5) ;
 
\draw [line width=1.5]  [dash pattern={on 5.63pt off 4.5pt}]  (53,508) -- (53,523.5) -- (53,541.5) ;
 
\draw     [fill={rgb, 255:yellow, 255; green, 0; blue, 0}  ,fill opacity=1 ] (27,546) -- (81,546) -- (81,600) -- (27,600) -- cycle ;
 
\draw   (12.33,264.5) -- (95,264.5) -- (95,347.17) -- (12.33,347.17) -- cycle ;
 
\draw   (120,305.84) .. controls (124.67,305.9) and (127.03,303.6) .. (127.08,298.93) -- (128.37,190.43) .. controls (128.45,183.76) and (130.82,180.46) .. (135.49,180.51) .. controls (130.82,180.46) and (128.53,177.1) .. (128.61,170.43)(128.57,173.43) -- (129.9,61.93) .. controls (129.95,57.26) and (127.65,54.9) .. (122.98,54.84) ;
 
\draw    (288,384) -- (325.62,344.3) ;
\draw [shift={(327,342.84)}, rotate = 133.46] [color={rgb, 255:red, 0; green, 0; blue, 0 }  ][line width=0.75]    (10.93,-3.29) .. controls (6.95,-1.4) and (3.31,-0.3) .. (0,0) .. controls (3.31,0.3) and (6.95,1.4) .. (10.93,3.29)   ;
 
\draw    (288,391.84) -- (324.59,428.43) ;
\draw [shift={(326,429.84)}, rotate = 225] [color={rgb, 255:red, 0; green, 0; blue, 0 }  ][line width=0.75]    (10.93,-3.29) .. controls (6.95,-1.4) and (3.31,-0.3) .. (0,0) .. controls (3.31,0.3) and (6.95,1.4) .. (10.93,3.29)   ;
 
\draw    (147,408.34) -- (147,555.84) ;
\draw [shift={(147,557.84)}, rotate = 270] [color={rgb, 255:red, 0; green, 0; blue, 0 }  ][line width=0.75]    (10.93,-3.29) .. controls (6.95,-1.4) and (3.31,-0.3) .. (0,0) .. controls (3.31,0.3) and (6.95,1.4) .. (10.93,3.29)   ;
 
\draw    (278,575) -- (315.62,535.3) ;
\draw [shift={(317,533.84)}, rotate = 133.46] [color={rgb, 255:red, 0; green, 0; blue, 0 }  ][line width=0.75]    (10.93,-3.29) .. controls (6.95,-1.4) and (3.31,-0.3) .. (0,0) .. controls (3.31,0.3) and (6.95,1.4) .. (10.93,3.29)   ;
 
\draw    (278,582.84) -- (314.59,619.43) ;
\draw [shift={(316,620.84)}, rotate = 225] [color={rgb, 255:red, 0; green, 0; blue, 0 }  ][line width=0.75]    (10.93,-3.29) .. controls (6.95,-1.4) and (3.31,-0.3) .. (0,0) .. controls (3.31,0.3) and (6.95,1.4) .. (10.93,3.29)   ;
 
\draw    (148,596.34) -- (148,660.84) ;
\draw [shift={(148,662.84)}, rotate = 270] [color={rgb, 255:red, 0; green, 0; blue, 0 }  ][line width=0.75]    (10.93,-3.29) .. controls (6.95,-1.4) and (3.31,-0.3) .. (0,0) .. controls (3.31,0.3) and (6.95,1.4) .. (10.93,3.29)   ;

\draw (58,171) node  [font=\small] [align=left] {\begin{minipage}[lt]{25.84pt}\setlength\topsep{0pt}
$ P_{2,1}$
\end{minipage}};
 
\draw (52.5,476.5) node  [font=\small] [align=left] {\begin{minipage}[lt]{31.96pt}\setlength\topsep{0pt}
$  P_{2,C_1 +1}$
\end{minipage}};
 
\draw (54,633) node  [font=\small] [align=left] {\begin{minipage}[lt]{31.28pt}\setlength\topsep{0pt}
$  P_{2,C_2 +1}$
\end{minipage}};
 
\draw (70,739.75) node  [font=\footnotesize] [align=left] {\begin{minipage}[lt]{96.56pt}\setlength\topsep{0pt}
$  P_{2,N} =R_2 \in Ch^{( N)}( Q_2)$
\end{minipage}};
 
\draw (62.5,60.5) node  [font=\small] [align=left] {\begin{minipage}[lt]{48.28pt}\setlength\topsep{0pt}
$  Q_{2} =P_{2,0}$
\end{minipage}};
 
\draw (54.5,395.5) node  [font=\small] [align=left] {\begin{minipage}[lt]{25.16pt}\setlength\topsep{0pt}
$  P_{2, C_1}$
\end{minipage}};
 
\draw (53.5,574) node  [font=\small] [align=left] {\begin{minipage}[lt]{23.8pt}\setlength\topsep{0pt}
$  P_{2,C_2}$
\end{minipage}};
 
\draw (57,303.5) node  [font=\small] [align=left] {\begin{minipage}[lt]{34pt}\setlength\topsep{0pt}
$  P_{2,C_1 -1}$
\end{minipage}};
 
\draw (346,179.42) node   [align=left] {\begin{minipage}[lt]{273.36pt}\setlength\topsep{0pt}
$\displaystyle \forall \ 0\leq k\leq C_{1} -1, \    \     \| S_1( \mathcal F_{N,k})( x_{1} ,\cdot ) \| _{L^{\infty }( P_{2,k})}^{2} \leq \lambda ^{2}$
\end{minipage}};
 
\draw (199,387) node   [align=left] {\begin{minipage}[lt]{125.8pt}\setlength\topsep{0pt}
$\displaystyle \| S_{1}( \mathcal F_{N, C_1}) \| _{L^{\infty }( P_{2, C_1})}^{2}  >\lambda ^{2}$
\end{minipage}};
 
\draw (481,350) node  [font=\small] [align=left] {\begin{minipage}[lt]{194.48pt}\setlength\topsep{0pt}
{\small Case 1. }$\displaystyle S_1^2( \mathcal F_{N, C_1})( x_1,\cdot ) >\lambda ^{2}$ on $ P_{2, C_1}$\\[4pt]{\small \color{blue}
Let $C_{R_2}^{I} =P_{2, C_1}$ and $C_{R_2}^{II}=P_{2, C_1}$
}
\end{minipage}};
 
\draw (477.5,429.92) node  [font=\small] [align=left] {\begin{minipage}[lt]{187pt}\setlength\topsep{0pt}
{\small Case 2. }$\displaystyle \min  S_1^{2}( \mathcal F_{N, C_1})( x_1 ,\cdot ) \leq \lambda ^2$ on $P_{2,C_1}$ \\  while $\displaystyle S_1^{2}(\mathcal  F_{N, C_1})( x_1 ,\cdot ) >\lambda^2$ on $  P_{2,C_1 +1}$\\[4pt]{\small \color{blue}
Let $C_{R_2}^{I} =P_{2, C_1}$ and $C_{R_{2}}^{II}=P_{2,C_1 +1}$
}
\end{minipage}};
 
\draw (199,579) node   [align=left] {\begin{minipage}[lt]{125.8pt}\setlength\topsep{0pt}
$\displaystyle \| S_{1}( \mathcal F_{N, C_2}) \| _{L^{\infty }( P_{2, C_2})}^{2}  >\lambda ^{2}$
\end{minipage}};
 
\draw (240.5,481.59) node  [font=\small] [align=left] {\begin{minipage}[lt]{117.64pt}\setlength\topsep{0pt}
{\small Case 3. Otherwise}\\   {\small \color{blue}Consider the next critical cube $  P_{2,C_{2}}$ similarly}\\
\end{minipage}};
 
\draw (471,541) node  [font=\small] [align=left] {\begin{minipage}[lt]{194.48pt}\setlength\topsep{0pt}
{\small Case 1. }$\displaystyle S_1^2( \mathcal F_{N, C_2})( x_{1} ,\cdot ) > \lambda ^{2}$ on $\displaystyle P_{2,C_2}$\\[4pt]{\small 
\color{blue} Let $ C_{R_{2}}^{I} =P_{2, C_2}$ and $ C_{R_{2}}^{II}=P_{2, C_{2}}$
}
\end{minipage}};
 
\draw (467.5,620) node  [font=\small] [align=left] {\begin{minipage}[lt]{187pt}\setlength\topsep{0pt}
{\small Case 2. }$\displaystyle \min S_1^2( F_{N, C_2})( x_1,\cdot ) \leq \lambda ^2$ on $ P_{2,C_2}$\\  while $\displaystyle S_1^2( F_{N, C_2})( x_1 ,\cdot ) > \lambda ^2$ on $ P_{2,C_2 +1}$\\[4pt]{\small \color{blue}
Let $C_{R_2}^{I} =P_{2, C_2}$ and $C_{R_2}^{II}=P_{2,C_2+1}$
}
\end{minipage}};
 
\draw (241.5,635.24) node  [font=\small] [align=left] {\begin{minipage}[lt]{117.64pt}\setlength\topsep{0pt}
{\small Case 3. Otherwise}\\  {\color{blue}Continue similarly}\\
\end{minipage}};
 
\draw (320,695) node   
[align=left] {\begin{minipage}[lt]{330pt}\setlength\topsep{0pt}
{\small \color{blue}
since the critical cubes $\displaystyle \{P_{2, C_k}\}_{k}$ with $\displaystyle \| S_{1}( \mathcal F_{N, C_k}) \| _{L^{\infty }( P_{2, C_k})}^{2}  >\lambda ^{2}$ are finite, and due to $\displaystyle \| S_1^*(\mathcal  F_{N,N}) \| _{L^{\infty }( R_{2})}^{2}  >\lambda ^{2}$, eventually we can choose   $\displaystyle C_{R_{2}}^{I}$ and $\displaystyle C_{R_{2}}^{II}$ }
\end{minipage}};
\end{tikzpicture}
\caption{ The selection of  $C_{R_2}^I$ and $C_{R_2}^{II}$ in the chain generated by $R_2\in Ch^{(N)}(Q_2)$ }
  \label{fig:2}
\end{figure}
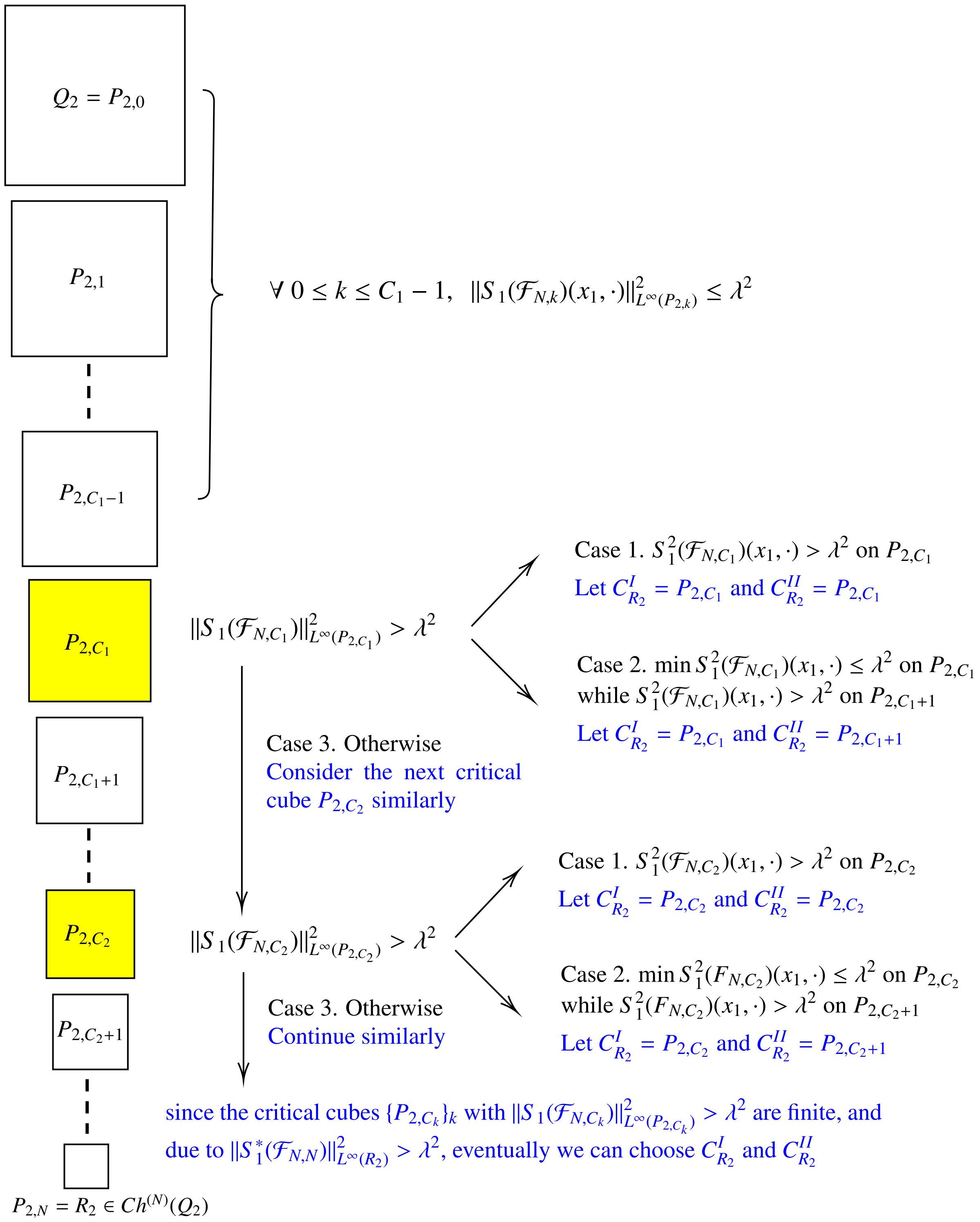

Now we define 
$$
  \mathcal G_2^{(I)}:=\left\{  C_{R_2}^I:\,   R_2\in \mathsf B\right\} \quad \text{and}\quad
     \mathcal G_2^{(II)}:=\left\{  C_{R_2}^{II}:\,   R_2\in \mathsf B  \right\}.
$$
Based on the algorithm above, we now decompose 
\begin{align}\label{precise decomposition}
    \Big\{x_2\in Q_2: S_1^*(\mathcal  F_{N,N})(x_1,x_2)>\lambda  \Big \}=\bigcup_{P_2\in  \mathcal G_2^{(II)}} P_2.
\end{align}
We highlight that for two distinct cubes $R_2,R'_2\in \mathsf B$, either   $C_{R_2}^{II}=C_{R'_2}^{II}$ or $C_{R_2}^{II}\cap C_{R'_2}^{II}=\emptyset$. Hence any two distinct elements in $\mathcal G_2^{(II)}$ are mutually disjoint. Thus, \eqref{precise decomposition} is a precise decomposition.

Besides, the role of $\mathcal{G}_2^{(I)}$ can be illustrated by noting that for each $P_2\in \mathcal G_2^{(II)}$, there exists a unique $\widetilde{P}_2\in \mathcal G_2^{(I)}$ such that one of the following situations happens:
\begin{itemize}
	\item $P_2=\widetilde{P}_2$;
	\item $P_2\in Ch(\widetilde{P}_2)$;
	\item $\widetilde{P}_2\in \mathsf  B$ and $P_2$ is the union of certain (at least two but not all) children of $\widetilde{P}_2$.
\end{itemize}

In particular,  the above algorithm  implies that for each $P_2\in \mathcal G_2^{(II)}$ associated to some $R_2\in Ch^{(N)}(Q)$,
$$
        \sum_{\substack{P_1\in  Ch^{(n)}(Q_1)\\ 0\leq n\leq N}} \Bigg| \sum_{P_2'\in \mathcal G_2:\, P_2'\supsetneqq \widetilde P_2}  \Delta_{P_1\times P_2'}(f) (x_1,\cdot)\Bigg|^2\leq  \lambda^2       \ \     \text{on}\  \    P_2.    \leqno{(\mathsf C1)}
$$
Otherwise, due to \textsf{(P2)}, we further have $\displaystyle \sum_{\substack{P_1\in  Ch^{(n)}(Q_1)\\ 0\leq n\leq N}} \Bigg| \sum_{P_2'\in \mathcal G_2:\, P_2'\supsetneqq \widetilde P_2}  \Delta_{P_1\times P_2'}(f) (x_1,\cdot)\Bigg|^2>  \lambda^2$ on  $\widetilde P_2$. Then we must have $C_{R_2}^I\supsetneqq \widetilde P_2$, which leads to a contradiction. On the other hand, it's clear that 
$$
        \sum_{\substack{P_1\in  Ch^{(n)}(Q_1)\\ 0\leq n\leq N}} \Bigg| \sum_{P_2'\in \mathcal G_2:\, P_2'\supseteq \widetilde P_2}  \Delta_{P_1\times P_2'} (f)(x_1,\cdot)\Bigg|^2>  \lambda^2       \ \     \text{on}\  \    P_2.    \leqno{(\mathsf C2)}
$$

Additionally, for any $P_2\in  {\mathcal G}_2^{(II)}$, if $P_2\in Ch^{(N)}(Q_2)$ or even $P_2\subsetneqq R_2$ for some $R_2\in Ch^{(N)}(Q_2)$, then such $P_2$ will have no contribution  
to the improved good-$\lambda$ inequality \eqref{eqn:good-lambda}, since now
\begin{equation}\label{eqn:TR}
	T_{P_2}:=\Big \{x_2\in P_2: S_1^*(\mathcal F_{N,N})(x_1,x_2)>2\lambda, \  \   \widetilde S(\mathcal F_{N,N})(x_1,x_2)\leq \epsilon \lambda \Big \}=\emptyset.
\end{equation}
Concretely,
if  $P_2\in Ch^{(N)}(Q_2)$ (the argument for $P_2\subsetneqq R_2$ for some $R_2\in Ch^{(N)}(Q_2)$ is similar and we skip it), then there exists a unique  $\widetilde{ P}_2\in \mathcal G_2^{(I)}$ satisfying one of the following cases.

\noindent{\it Case I. $\widetilde {P}_2=P_2\in Ch^{(N)}(Q_2)$}.  
By the algorithm above, we have 
$$
       S_1^2(\mathcal  F_{N, N-1})(x_1,\cdot) \leq \lambda^2 \  \   \text{on}\   \   P_2,
$$
and  
$$
   S_1^*(\mathcal  F_{N,N})(x_1,\cdot)= S_1(\mathcal  F_{N,N})(x_1,\cdot) >\lambda\  \   \text{on}\   \   P_2.
$$

Moreover,  if $T_{P_2}\neq \emptyset$, then 
$S_1^*(\mathcal  F_{N,N})(x_1,\cdot)=S_1(\mathcal  F_{N,N})(x_1,\cdot )>2\lambda$ on $T_{P_2}$. Combining $\widetilde S(\mathcal  F_{N,N})(x_1,x_2)\leq \epsilon \lambda $ for $x_2\in T_{P_2}$, we have for $x_2\in T_{P_2}$
\begin{align*}
	4\lambda^2 \leq &\sum_{\substack{P_1\in  Ch^{(n)}(Q_1)\\ 0\leq n\leq N}} \Bigg| \sum_{P_2'\in \mathcal G,\, P_2'\supsetneqq P_2}  \Delta_{P_1\times P_2'}(f) (x_1,x_2) \Bigg|^2 +\sum_{\substack{P_1\in  Ch^{(n)}(Q_1)\\ 0\leq n\leq N}} \left|\Delta_{P_1\times P_2}(f) (x_1,x_2) \right |^2 \\
	&  +2\sum_{\substack{P_1\in  Ch^{(n)}(Q_1)\\ 0\leq n\leq N}} \left[  \sum_{P_2'\in \mathcal G_2,\, P_2'\supsetneqq P_2}  \Delta_{P_1\times P_2'} (f)  (x_1,x_2) \cdot \Delta_{P_1 \times P_2} (f) (x_1,x_2) \right]\\
	\leq & S_1^2(\mathcal F_{N, N-1})(x_1,x_2)+ \widetilde S^2(\mathcal F_{N,N})(x_1,x_2)+ 2 S_1(\mathcal F_{N,N+1})(x_1,x_2)\cdot \widetilde S(\mathcal F_{N,N})(x_1,x_2)\\
	\leq & \lambda^2 +\epsilon^2\lambda^2  +2\lambda \cdot \epsilon\lambda,
\end{align*}
which leads to a contradiction whenever $\epsilon<1/2$.

\smallskip

\noindent{\it Case II. $P_2\in Ch(\widetilde P_2)$. } Then $\widetilde P_2\in Ch^{(N-1)}(Q_2)$, and by the  above algorithm again,
$$
       S_1^2(\mathcal  F_{N, N-2})(x_1,\cdot) \leq \lambda^2 \  \   \text{on}\   \   \widetilde P_2,
$$
and  for each  $x_2\in   P_2$,
$$
   \lambda<S_1^*(\mathcal  F_{N,N})(x_1,x_2) = S_1(\mathcal  F_{N,N-1})(x_1,x_2) \  \text{or}\  =S_1(\mathcal  F_{N,N})(x_1,x_2) .
$$

Given $x_2\in T_{P_2}$, we have the following two cases:

\noindent $\bullet$ If $\displaystyle S_1(\mathcal  F_{N,N-1})(x_1,x_2)>\frac{3\lambda}{2}$, then similarly,
\begin{align*}
	\Big(\frac{3\lambda}{2}\Big)^2 < \, & S_1^2 (\mathcal  F_{N,N-1})(x_1,x_2) \\
	\leq\,&  S_1^2(\mathcal  F_{N,N-2}) (x_1,x_2)+ \widetilde S^2(\mathcal  F_{N,N})(x_1,x_2) +2 S_1(\mathcal  F_{N,N-2}) (x_1,x_2)\cdot \widetilde S(\mathcal F_{N,N})(x_1,x_2) \\
	\leq \,&  \lambda^2 +\epsilon^2\lambda^2  +2\lambda \cdot \epsilon\lambda,
\end{align*}
which contradicts to the assumption that $\epsilon$ is sufficiently small.

\noindent $\bullet$ Otherwise, if $\displaystyle S_1(\mathcal  F_{N,N-1})(x_1,x_2)\leq \frac{3\lambda}{2}$, then $S_1^*(\mathcal  F_{N,N})(x_1,x_2)= S_1(\mathcal  F_{N,N})(x_1,x_2)>2\lambda $. Similarly, estimates on $S_1(\mathcal  F_{N,N})(x_1,x_2)$ yield a contradiction again.

Therefore, in the following discussion, we always assume that the cubes and atomic functions involved exist, and the above meaningless situations do not occur. That is, we may assume that 
$$
    \mathcal G_2^{(II)}=\left\{C_{R_2}^{II}:\,   R_2\in \mathsf B\quad \text{and}\quad C_{R_2}^{II} \supsetneqq R_2\right\}.
$$

Recall that for any $P_2\in  \mathcal G_2^{(II)}$ and the corresponding $\widetilde P_{2}\in  \mathcal G_2^{(I)}$, we have  
$$
   P_2'\subseteq P_2 \  \    \text{ for\    each }   \     
      P_2'\subsetneqq \widetilde P_{2}.
$$ 
 This, in combination with 
 \eqref{eqn:one-shift-bi-I-weak-martingale}, yields that for any $P_2\in  \mathcal G_2^{(II)}$ and $\alpha>0$,
\begin{align}\label{eqn:one-shift-bi-iteration-II}
	&\int_{P_2} \exp \Bigg(\alpha\sqrt{\sum_{\substack{P_1\in  Ch^{(n)}(Q_1)\\ 0\leq n\leq N}} \bigg|\sum_{\substack{P_2'\in  \mathcal G_2\\ P_2'\subsetneqq \widetilde P_2}} \Delta_{P_1\times P_2'}(f)(x_1,x_2)\bigg|^2}  \Bigg.\\
	&\qquad\qquad  \Bigg.   - c\alpha^2  \sum_{\substack{P_1\in  Ch^{(n)}(Q_1)\\ 0\leq n\leq N}}  \sum_{\substack{P_2'\in  \mathcal G_2\\ P_2'\subsetneqq \widetilde P_2}} \left\|\Delta_{P_1\times P_2'}(f)\right\|_{L^\infty}^2 \Bigg) \, {\rm d}\mu_2(x_2) 
	\leq \,  C\mu_2(P_2).\nonumber
\end{align}

Similarly, using the algorithm above, we can continue to decompose for  $P_2\in  \mathcal G_2^{(II)}$ that 
\begin{equation}\label{eqn:Rj}
	 \Big\{x_2\in Q_2: S_1^*(\mathcal F_{N,N})(x_1,x_2)>2\lambda  \Big \}\cap P_2=\bigcup_j P_{2,j} \   \    \text{(disjoint)}
\end{equation}
such that for each $P_{2,j}$,  
there exists $\widetilde P_{2,j}$ with $\widetilde P_{2,j}=  P_{2,j}$ or $ P_{2,j}\in Ch(\widetilde P_{2,j})$, satisfying
$$
        \sum_{\substack{P_1\in  Ch^{(n)}(Q_1)\\ 0\leq n\leq N}} \Bigg| \sum_{P_2'\in \mathcal G_2:\, P_2'\supsetneqq \widetilde P_{2,j}}  \Delta_{P_1\times P_2'}(f) (x_1,\cdot )\Bigg|^2\leq (2\lambda)^2       \ \     \text{on}\  \    P_{2,j},    \leqno{(\mathsf C1')}
$$
and 
$$
        \sum_{\substack{P_1\in  Ch^{(n)}(Q_1)\\ 0\leq n\leq N}} \Bigg| \sum_{P_2'\in \mathcal G_2:\, P_2'\supseteq \widetilde P_{2,j}}  \Delta_{P_1\times P_2'}(f) (x_1,\cdot)\Bigg|^2 > (2\lambda)^2       \ \     \text{on}\  \    P_{2,j}. \leqno{(\mathsf C2')}
$$

It suffices to consider the case $\widetilde P_{2,j}\subsetneqq \widetilde P_2$ for each $j$, by noting that by  a contradiction similar to the proof of \eqref{eqn:TR}, one may verify that 
 $  \left\{  x_2\in P_{2,j}: S(\mathcal F_{N,N})(x_1,x_2)\leq \epsilon\lambda\right\}=\emptyset $ if $\widetilde P_{2,j}=\widetilde P_2$. This ensures that 
 $$
    \left\{ P_2':\  \widetilde P_{2,j}\subseteq P_2'\subsetneqq \widetilde P_2\   \text{for\  some\  }j\right\}\neq \emptyset.
 $$
Note that the  lacunary type Hoeffding's inequality  \eqref{eqn:one-shift-bi-lacu} also holds. In particular, one may  
constrain the set of $P_2'$'s in \eqref{eqn:one-shift-bi-iteration-II}  as the above $ \left\{ P_2':\  \widetilde P_{2,j}\subseteq P_2'\subsetneqq \widetilde P_2\   \text{for\  some\  }j\right\}$, then as a consequence of \eqref{eqn:one-shift-bi-iteration-II},
\begin{align} \label{eqn:one-shift-bi-iteration-III}
	&\sum_j \int_{P_{2,j}\cap \left\{\widetilde S(\mathcal F_{N,N})(x_1,\cdot )\leq \epsilon \lambda\right\}} \exp \bigg(\alpha \sqrt{\sum_{\substack{P_1\in  Ch^{(n)}(Q_1)\\ 0\leq n\leq N}} \bigg|\sum_{ P_2':\, \widetilde P_{2,j}\subseteq P_2'\subsetneqq \widetilde P_2} \Delta_{P_1\times P_2'}(f)(x_1,x_2)\bigg|^2}  \bigg. \nonumber\\
	&\qquad \qquad \qquad \qquad \qquad \qquad  \bigg.  - c\alpha^2  \sum_{\substack{P_1\in  Ch^{(n)}(Q_1)\\ 0\leq n\leq N}}  \sum_{P_2':\, \widetilde P_{2,j}\subseteq P_2'\subsetneqq \widetilde P_2} \left\|\Delta_{P_1\times P_2'}(f) \right\|_{L^\infty}^2 \bigg) \, {\rm d}\mu(x_2)\\
	\leq \,&  C\mu_2(P_2).\nonumber
\end{align}
  
For each $j$ and $x_2\in P_{2,j}\cap \left\{\widetilde S(\mathcal F_{N,N})(x_1,\cdot)\leq \epsilon \lambda\right\}$, combining ($\mathsf C$1) and ($\mathsf C$2'), we have 
\begin{align*}
	2\lambda <&\sqrt{\sum_{\substack{P_1\in  Ch^{(n)}(Q_1)\\ 0\leq n\leq N}} \bigg|\sum_{ P_2'\in \mathcal G_2:\,   P_2'\supseteq \widetilde P_{2,j}} \Delta_{P_1\times P_2'}(f)(x_1,x_2)\bigg|^2}\\
	\leq &
	\sqrt{\sum_{\substack{P_1\in  Ch^{(n)}(Q_1)\\ 0\leq n\leq N}} \bigg|\sum_{ P_2':\, \widetilde P_{2,j}\subseteq P_2'\subsetneqq \widetilde P_2} \Delta_{P_1\times P_2'} (f)\bigg|^2}+\sqrt{\sum_{\substack{P_1\in  Ch^{(n)}(Q_1)\\ 0\leq n\leq N}} \big| \Delta_{P_1\times \widetilde P_2}(f) \big|^2}+\sqrt{\sum_{\substack{P_1\in  Ch^{(n)}(Q_1)\\ 0\leq n\leq N}} \bigg|\sum_{ P_2':\,  P_2'\supsetneqq \widetilde P_2} \Delta_{P_1\times P_2'}(f) \bigg|^2}\\
	\leq & \sqrt{\sum_{\substack{P_1\in  Ch^{(n)}(Q_1)\\ 0\leq n\leq N}} \bigg|\sum_{ P_2':\, \widetilde P_{2,j}\subseteq P_2'\subsetneqq \widetilde P_2} \Delta_{P_1\times P_2'} (f)\bigg|^2}+ \epsilon \lambda +\lambda.
\end{align*}
 Hence
 $$
      \sqrt{\sum_{\substack{P_1\in  Ch^{(n)}(Q_1)\\ 0\leq n\leq N}} \bigg|\sum_{ P_2':\, \widetilde P_{2,j}\subseteq P_2'\subsetneqq \widetilde P_2} \Delta_{P_1\times P_2'}(f) \bigg|^2}\geq (1-\epsilon)\lambda
 $$
 on $P_{2,j}\cap \left\{\widetilde S(\mathcal  F_{N,N})(x_1,\cdot)\leq \epsilon \lambda\right\}$. Plugging this into \eqref{eqn:one-shift-bi-iteration-III}, we have that for any $\alpha>0$,
 \begin{align*}
 	\sum_j \exp\left(\alpha (1-\epsilon)\lambda -c\alpha^2 \epsilon^2\lambda^2\right)\mu_2 \left(  \left\{ P_{2,j}\cap \left\{ \widetilde S(\mathcal F_{N,N})(x_1,\cdot )\leq \epsilon \lambda\right\}\right\}\right)\leq  C\mu_2(P_2).
 \end{align*}
 Then following \eqref{eqn:Rj}, we have
 $$
     	\exp\left(\alpha (1-\epsilon)\lambda -c\alpha^2 \epsilon^2\lambda^2\right)\mu_2 \left(  \left\{P_{2}\cap \left\{S_1^*(\mathcal  F_{N,N})(x_1,\cdot)>2\lambda ,\  \  \widetilde S(\mathcal F_{N,N})(x_1,\cdot )\leq \epsilon \lambda\right\}\right\}\right)\leq  C\mu_2(P_2).
 $$
 Summing over these $P_2\in \mathcal G_2^{(II)}$, for any given $\lambda$ and $\epsilon$ sufficiently small (one may assume $\epsilon<1/2$),  we have that for any $\alpha>0$,
 \begin{align*}
	 &\mu_2\Big(\Big \{x_2\in Q_2: S_1^*(\mathcal F_{N,N})(x_1,x_2)>2\lambda, \  \   \widetilde S(\mathcal  F_{N,N})(x_1,x_2)<\epsilon \lambda \Big \}\Big)  \\
	\leq \,& C\exp\left(-\alpha (1-\epsilon)\lambda +c\alpha^2 \epsilon^2\lambda^2\right) \mu_2\Big(\Big\{x_2\in Q_2: S_1^*(\mathcal  F_{N,N})(x_1,x_2)>\lambda  \Big \}\Big).
\end{align*}
 Note that 
 $$
     -\alpha (1-\epsilon)\lambda +c\alpha^2 \epsilon^2\lambda^2=c\epsilon^2\lambda^2\left[\left(\alpha-\frac{(1-\epsilon)\lambda}{2c\epsilon^2\lambda^2}\right)^2 -\frac{(1-\epsilon)\lambda^2}{4c^2\epsilon^4\lambda^4}\right]\geq -\frac{(1-\epsilon)^2}{4c\epsilon^2}.
 $$
 Taking $\displaystyle \alpha=\frac{(1-\epsilon)\lambda}{2c\epsilon^2\lambda^2}$ and assuming $\epsilon<1/2$,  we obtain 
 \begin{align*}
	 &\mu_2\Big(\Big \{x_2\in Q_2: S_1^*(\mathcal  F_{N,N})(x_1,x_2)>2\lambda, \   \  \widetilde S(\mathcal F_{N,N})(x_1,x_2)<\epsilon \lambda \Big \}\Big)  \\
	\leq \,& C\exp\left(-\frac{1}{8c\epsilon^2} \right) \mu_2\Big(\Big\{x_2\in Q_2: S_1^*(\mathcal  F_{N,N})(x_1,x_2)>\lambda  \Big \}\Big)
\end{align*}
as desired.
Hence we complete the proof of Proposition \ref{coro:good-lambda}.
\end{proof}

\smallskip

As a consequence, we can obtain the sharp $L^p$ lower bound for the square function $\widetilde S$ associated to dyadic martingale differences that
\begin{equation}\label{eqn:part-II-Lp-CH-tilde}
	   \left\|S_1(\mathcal F_{N,N})(x_1, \cdot)\right\|_{L^p(Q_2,{\rm d}\mu_2)}\leq C\sqrt{p} \left\|\widetilde S(\mathcal F_{N,N})(x_1,\cdot)\right\|_{L^p(Q_2,{\rm d}\mu_2)}\quad \text{for}\   p\geq 2,
\end{equation}
without using the $A_1$ weight approach as in Sections \ref{sec:3} and \ref{sec:4}.  

\medskip

{\bf Proof of   \eqref{eqn:part-II-Lp-CH-tilde}:} 

For any $\lambda >0$ and any fixed $x_1\in Q_1$,  we split  
\begin{align*}
	&\Big\{x_2\in Q_2:\,S_1^*(\mathcal F_{N,N})(x_1,x_2)>2\lambda\Big \}\\
	=&\Big \{S_1^*(\mathcal F_{N,N})>2\lambda,\  \  \widetilde S(\mathcal F_{N,N})>\epsilon\lambda\Big\} \bigcup  \Big\{S_1^*(\mathcal  F_{N,N})>2\lambda,\  \  \widetilde S(\mathcal F_{N,N})\leq \epsilon\lambda\Big\} \\
	\subseteq & \Big \{ \widetilde S(\mathcal F_{N,N})>\epsilon\lambda\Big\} \bigcup  \Big\{S_1^*(\mathcal F_{N,N})>2\lambda, \   \   \widetilde S(\mathcal F_{N,N})\leq \epsilon\lambda\Big\} .
\end{align*}
Using Proposition \ref{coro:good-lambda}, 
\begin{align*}
	&\int_{Q_2} \Big|S_1^* (\mathcal F_{N,N})(x_1,x_2)\Big|^p {\rm d}\mu_2(x_2)\\
	=& p \int_0^\infty (2\lambda)^{p-1}\mu_2 \Big  (
	\Big  \{   x_2\in Q_2:\,    S_1^*(\mathcal F_{N,N}) (x_1, x_2)>2\lambda \Big  \} 
	\Big ) \, {\rm d}(2\lambda)\\
	\leq & \frac{2^p }{\epsilon^p}\int_{Q_2} \big|\widetilde S(\mathcal F_{N,N}) (x_1,x_2) \big|^p \, {\rm d}\mu_2(x_2) + 2^pC \exp\left(-\frac{C'}  {\epsilon^2}\right)   \int_{Q_2} \big|S_1^*(\mathcal F_{N,N}) (x_1,x_2)\big|^p  \, {\rm d}\mu_2(x_2) .
\end{align*}
Taking $\epsilon$ such that $2^pC \exp\left(-\frac{C'}  {\epsilon^2}\right) \approx 1/2$, that is, $\epsilon\approx p^{-1/2}$, we obtain \eqref{eqn:part-II-Lp-CH-tilde} readily.
 \hfill $\Box$

Alternatively, \eqref{eqn:part-II-Lp-CH-tilde} follows from Proposition \ref{prop:double-Lp-CH} directly by $S(\mathcal F_{N,N})(x_1,x_2)\leq \widetilde S(\mathcal F_{N,N})(x_1,x_2)$.
We highlight again that such square function $\widetilde S$ is hard to use in multi-parameter settings.

\section*{Subsequent applications}

Since we completed this paper in 2024, both the one-parameter Hoeffding inequality and the intermediate square functions have found further applications.

The almost orthogonal atoms in Theorem~\ref{thm:newCH-dyadic} may have nonzero mean. Wu~\cite{WuJGA} used this inequality to prove the classical exponential-square estimate of Chang, Wilson and Wolff~\cite{CWW} and its operator version due to Chen, Duong, Wu and Yan~\cite{CDWY} by the same argument. The proof does not require the heat semigroup to preserve constants.

For Dunkl area integrals, the orbit distance governing heat kernel decay agrees on each Weyl chamber with the Euclidean distance used in the regularity estimate. The Dunkl heat semigroup preserves constants on the whole space, so the reconstruction atoms have integral zero. Their integrals over individual chambers, however, may be nonzero. Chang et al.~\cite{DunklArea} used the Hoeffding inequality on each chamber to prove the reverse $L^p$ estimate for Dunkl area integrals with the sharp factor $\sqrt p$, $p>2$.

Li and Wu~\cite{LWZygmund} adapted the intermediate square function of Section~\ref{sec:4} to the coupled scales $s,t,st$ of Zygmund dilations. At a fixed aperture, the sharp reverse $L^p$ bounds have order $p$, $p^{3/2}$ or $p^2$, depending on the class of kernels. For a $k$-step Euclidean flag, Han, Li, Wick and Wu~\cite{HLWWFlag} used decay in the scale ratios to sum the terms that lack cancellation in a coordinate. On any product $R$ of coordinate cubes, the flag Poisson area on $R$ controls the oscillation left after subtracting each coordinate mean in turn. The sharp factor in this local $L^p$ estimate is $p^{k/2}$.

For products of $n$ Dunkl heat area integrals, the time scales are independent. Within each factor, reflections couple the values of the function on different chambers. By estimating the nonzero chamber means, Han, Li, Wick and Wu~\cite{HLWWDunkl} were able to iterate the intermediate square function estimates and prove local exponential integrability of power $2/n$ for the reconstructed sums.
\vskip   1cm

\noindent{\bf Acknowledgments.} Li and Wu are supported by ARC DP220100285. Wu is  also supported by  NNSF of China \# 12201002, Anhui NSF of China \# 2208085QA03 and Excellent University Research and Innovation Team in Anhui Province \# 2024AH010002.

\vskip 0.5cm


\end{document}